\pdfoutput=1
\RequirePackage{ifpdf}
\ifpdf % We are running pdfTeX in pdf mode
\documentclass[pdftex]{sigma}
\else
\documentclass{sigma}
\fi

\usepackage[mathscr]{eucal}
\usepackage{amsfonts,dsfont}
\usepackage{stmaryrd}

\numberwithin{equation}{section}

\newtheorem{Theorem}{Theorem}[section]
\newtheorem{Corollary}[Theorem]{Corollary}
\newtheorem{Lemma}[Theorem]{Lemma}
\newtheorem{Proposition}[Theorem]{Proposition}

{\theoremstyle{definition}
\newtheorem{Definition}[Theorem]{Definition}
\newtheorem{Example}[Theorem]{Example}
\newtheorem{Remark}[Theorem]{Remark}
}

\newcommand{\n}{^{(n)}}
\newcommand{\ii}{^{(\infty )}}
\newcommand{\vv}{\mathbf{v}}
\newcommand{\ww}{\mathbf{w}}
\newcommand{\VV}{\mathbf{V}}

\newcommand{\bb}{\pmb{\beta}}
\newcommand{\bl}{\pmb{\lambda}}
\newcommand{\bt}{\pmb{\tau}}

\newcommand{\bs}{\pmb{\sigma}}

\newcommand{\g}{\mathfrak{g}}

\newcommand{\D}{\mathcal{D}}

\newcommand{\G}{\mathcal{G}}

\newcommand{\bH}{\mathbf{H}}
\newcommand{\Jn}{\J^n}

\newcommand{\bO}{\mathbf{\Omega}}
\newcommand{\bo}{\pmb{\omega}}
\newcommand{\bp}{\mathbf{p}}
\newcommand{\bS}{\pmb{\Sigma}}
\newcommand{\bT}{\mathbf{T}}

\def\dom{\mathop{\rm dom}\nolimits}
\def\dim{\mathop{\rm dim}\nolimits}

\newcommand{\pp}[2]{\frac{\partial #1}{\partial #2}}

\def\comp{\raise 1pt \hbox{$\,\scriptstyle\circ\,$}}

\def\ro#1{{\rm #1}}

\def\j {\ro j}
\def\J{\ro J}
\def\interior{\mathbin{\hbox{\hbox{{\vbox
    {\hrule height.4pt width6pt depth0pt}}}\vrule height6pt width.4pt
depth0pt}\,}}

\begin{document}

\allowdisplaybreaks

\renewcommand{\thefootnote}{$\star$}

\renewcommand{\PaperNumber}{029}

\FirstPageHeading

\ShortArticleName{Solving Local Equivalence Problems
with the Equivariant Moving Frame Method}

\ArticleName{Solving Local Equivalence Problems\\
with the Equivariant Moving Frame Method\footnote{This paper
is a~contribution to the Special Issue ``Symmetries of Dif\/ferential Equations: Frames, Invariants and~Applications''.
The full collection is available at
\href{http://www.emis.de/journals/SIGMA/SDE2012.html}{http://www.emis.de/journals/SIGMA/SDE2012.html}}}

\Author{Francis VALIQUETTE}

\AuthorNameForHeading{F.~Valiquette}

\Address{Department of Mathematics and Statistics, Dalhousie University,\\
Halifax, Nova Scotia, B3H 3J5, Canada}

\Email{\href{mailto:francisv@mathstat.dal.ca}{francisv@mathstat.dal.ca}}
\URLaddress{\url{http://www.mathstat.dal.ca/~francisv}}

\ArticleDates{Received July 21, 2012, in f\/inal form March 31, 2013; Published online April 05, 2013}

\Abstract{Given a~Lie pseudo-group action, an equivariant moving frame exists in the neighborhood of a~submanifold jet
provided the action is free and regular.
For local equiva\-len\-ce problems the freeness requirement cannot always be satisf\/ied and in this paper we show that,
with the appropriate modif\/ications and assumptions, the equivariant moving frame constructions extend to submanifold
jets where the pseudo-group does not act freely at any order.
Once this is done, we review the solution to the local equivalence problem of submanifolds within the equivariant
moving frame framework.
This of\/fers an alternative approach to Cartan's equivalence method based on the theory of $G$-structures.}


\Classification{53A55; 58A15}
\vspace{-3mm}

\begin{flushright}
\begin{minipage}{67mm}
\it In honor of Peter Olver's $60^{\text{th}}$ birthday.
Thank you for your guidance.
\end{minipage}
\end{flushright}

\vspace{-3mm}

\renewcommand{\thefootnote}{\arabic{footnote}}
\setcounter{footnote}{0}

\section{Introduction}

\looseness=-1
First introduced by the Estonian mathematician Martin Bartels and primarily developed by \'Elie Cartan,~\cite{AR-1993},
the method of moving frames is a~powerful tool for studying geometric properties of submanifolds under the action of
a~(pseudo-)group of transformations.
In 1999, Fels and Olver proposed in~\cite{FO-1999} a~new theoretical foundation to the method of moving frames.
For a~Lie group~$G$ acting on the $n^\text{th}$ order jet space $\J^n(M,p)$ of $p$-dimensional submanifolds of~$M$,
a~moving frame is a~$G$-equivariant section of the trivial bundle $\J^n(M,p)\times G\to\J^n(M,p)$.
This new framework to moving frames, now known as the \emph{equivariant moving frame} method, possesses some appealing
features.
First, it decouples the moving frame theory from reliance on any form of frame bundle or connection and can thereby be
applied to almost any type of group action.
Secondly, every equivariant moving frame comes with an invariantization map that sends dif\/ferential functions,
dif\/ferential forms, and vector f\/ields to their invariant counterparts yielding a~complete collection of (local)
dif\/ferential invariants, invariant dif\/ferential forms, and invariant vector f\/ields.
In general, the invariantization map and the exterior dif\/ferential do not commute, and this lack of commutativity is
encapsulated in the \emph{universal recurrence formula} which is at the heart of many new results in the f\/ield.
For example, using this fundamental formula, Kogan and Olver were able to obtained in~\cite{KO-2003} a~general
group-invariant formula for the Euler--Lagrange equations of an invariant variational problem, while in~\cite{TV-2010}
the same formula was used to show that the cohomology of the invariant Euler--Lagrange complex is isomorphic to the Lie
algebra cohomology of its symmetry group.
But more importantly, the universal recurrence formula is the key that unveils the structure of the algebra of
dif\/ferential invariants,~\cite{FO-1999,H-2009,OP-2009-1}, essential to the solution of local equivalence problems.

\looseness=1
Recently, the theory of equivariant moving frames was successfully extended to inf\/inite-dimensional Lie pseudo-group
actions in~\cite{OP-2009-1, OP-2005,OP-2008}; opening the way to many new applications.
The f\/irst application appeared in~\cite{COP-2008} where the algebra of dif\/ferential invariants of the
inf\/inite-dimensional symmetry group of the Kadomtsev--Petviashvili equation was completely characterized.
An application to the classif\/ication of Laplace invariants and the factorization of linear partial dif\/ferential
operator can be found in~\cite{SM-2008}, and an adaptation of Vessiot's group foliation method using moving frames was
proposed in~\cite{TV} (see~\cite{P-2008} for an alternative construction based on exterior dif\/ferential systems).
As a~further application, in this paper we revisit the solution of the local equivalence problem of submanifolds under
an inf\/inite-dimensional Lie pseudo-group action using the equivariant moving frame formalism.
Following Cartan, the solution is essentially obtained by constructing suf\/f\/iciently many invariants so that one can
distinguish inequivalent submanifolds.
With the equivariant moving frame method, these invariants are easily constructed using the invariantization map and
their signature manifold is completely characterized by the universal recurrence formula.
Since the algorithms do not rely on the theory of exterior dif\/ferential systems and
$G$-structures,~\cite{BCGGG-1991,G-1989,K-1989,L-1924,O-1995}, the solution based on the equivariant moving frame
method of\/fers an interesting alternative to Cartan's equivalence method of coframes.

\looseness=1
To construct an equivariant moving frame in the neighborhood of a~submanifold jet, the pseudo-group action must be
(locally) free and regular.
Unfortunately, given an equivalence problem, the freeness requirement cannot always be satisf\/ied, and more often than
not many interesting results occur at the submanifold jets where the action cannot be made free by prolongation.
For example, it is well-known that a~second order ordinary dif\/ferential equation $u_{xx}=Q(x,u,u_x)$ is equivalent to
$u_{xx}=0$ under a~point transformation if and only if it admits an eight-dimensional symmetry group isomorphic to
$\text{SL}(3)$,~\cite{C-1955,G-1989,O-1995,T-1896}.
For such a~dif\/ferential equation, the pseudo-group of point transformations cannot act freely and the freeness
assumption must be relaxed if one wants to obtain a~complete solution of the local equivalence problem using the
equivariant moving frame method.
As one might expect, the idea is to modify the standard moving frame algorithms by incorporating the isotropy group
into the constructions to obtain what we call a~\emph{partial equivariant moving frame}.
Using (partial) moving frames we can solve a~wide range of local equivalence problems, which includes equivalence
problems between coframes.
To illustrate the method we consider the local equivalence problem of second order ordinary dif\/ferential equations
under point transformations and contact transformations, and the simultaneous equivalence of a~two-form and a~vector
f\/ield on $\mathbb{R}^3$.
By revisiting these standard examples, our goal is to highlight some links between the (partial) equivariant moving
frame approach and Cartan's method.

The solution of a~local equivalence problem relies on the \emph{fundamental basis theorem} (also known as the
Lie--Tresse theorem) which states that, under appropriate regularity assumptions, the algebra of dif\/ferential
invariants on $\J^\infty(M,p)$ is locally generated by a~f\/inite set of dif\/ferential invariants and exactly $p$
linearly independent invariant total derivative operators.
Under the assumption that a~pseudo-group action is regular and locally free on a~dense open subbundle of~$\J^\infty(M,p)$, a~constructive proof of the fundamental basis theorem based on the equivariant mo\-ving frame method
was recently given in~\cite{OP-2009-1}.
In Section~\ref{algebra of differential invariants - section}, we adapt the algebraic constructions introduced
in~\cite{OP-2009-1} to cover pseudo-groups acting regularly and freely on invariant regular subbundles of
$\J^\infty(M,p)$ and also consider the case of regular pseudo-groups acting non-freely on invariant regular subbundles
of $\J^\infty(M,p)$.
These adaptations are necessary to give a~complete solution to the local equivalence problem of submanifolds.

\begin{Remark}
The theory of inf\/inite-dimensional Lie pseudo-groups relies on the Cartan--K\"ahler
theorem,~\cite{BCGGG-1991,O-1995}, which requires analyticity.
For this reason, all our constructions and results hold in the analytic category.
Implicitly, all manifolds, maps, dif\/ferential forms and vector f\/ields are thus assumed to be analytic.
For Lie pseudo-groups of f\/inite type, in other words for local Lie group actions, analyticity can be replaced by
smoothness.
\end{Remark}

\begin{Remark}
%\label{local}
Following the global notation convention used in~\cite{FO-1999,OP-2009-1,OP-2005,OP-2008}, given a~map $\varphi\colon
M\to N$ between two manifolds $M$ and $N$ we allow the domain of the map to be a~proper open subset of the manifold
$M$: dom $\varphi\subset M$.
Hence, while we use global notation throughout the exposition, all moving frame constructions and results should be
understood to hold micro-locally, i.e.\
on open subsets of the submanifold jet bundle $\J^\infty(M,p)$.
Similarly, dif\/ferential forms and vector f\/ields on $\J^\infty(M,p)$ are assumed to be def\/ined micro-locally.
\end{Remark}

\section{Structure equations}
%\label{Section Structure Theory}

Following~\cite{COP-2005,OP-2005} we begin by recalling how the structure equations of a~Lie pseudo-group are obtained
from its inf\/initesimal data.
As we will see in Section~\ref{Section Equivalence Problems}, the structure equations of an equivalence pseudo-group
provide the link between the equivariant moving frame method and Cartan's moving frame approach.

\subsection{Dif\/feomorphism pseudo-group}

Let $M$ be an $m$-dimensional manifold.
We denote by $\D=\D(M)$ the pseudo-group of all local dif\/feomorphisms of $M$.
For each integer $0\leq n\leq\infty$, let $\D\n$ denote the bundle formed by their $n^\text{th}$ order jets.
For $k\geq n$, let $\widetilde{\pi}^k_n\colon\D^{(k)}\to\D\n$ denote the standard projection.
Following Cartan,~\cite{C-1953-2, C-1953-1}, and the recent work of Olver and
Pohjanpelto,~\cite{OP-2009-1, OP-2005,OP-2008}, in some local coordinate system we use lower case letters, $z, x,u, \ldots$ for the \emph{source coordinates} \mbox{$\widetilde{\bs}(\varphi)\!=\!z\!=\!(z^1,\ldots,z^m)\in M$} of a~local
dif\/feomorphism $Z=\varphi(z)$ and corresponding upper case letters $Z, X, U, \ldots$ for the \emph{target
coordinates} \mbox{$\widetilde{\bt}(\varphi)=Z=(Z^1,\ldots,Z^m)\in M$}.
The local coordinates of the $n$-jet of a~local dif\/feomorphism $\varphi$ are then given by $\j_n\varphi=(z,Z\n)$,
where $z$ are the source coordinates and~$Z\n$ denotes the collection of derivatives
$Z^a_B=\partial^k\varphi^a/\partial z^{b^1}\cdots\partial z^{b^k}$ with $1\leq a,b^1,\ldots,b^k\leq m$ and $0\leq
k=\#B\leq n$.

The dif\/feomorphism jet bundle $\D\ii$ has the structure of a~groupoid,~\cite{M-1987}.
The groupoid multiplication follows from the composition of local dif\/feomorphisms.
For $g\ii|_z=\j_\infty\varphi|_z$ and $h\ii|_Z=\j_\infty\psi|_Z$ with
$Z=\widetilde{\bt}\ii(\j_\infty\varphi|_z)=\widetilde{\bs}\ii(\j_\infty\psi|_Z)$, we have that $(h\ii\cdot
g\ii)|_z=\j_\infty(\psi\comp\varphi)|_z$.
Throughout the paper, the dif\/feomorphism pseudo-group $\D$ acts on $\D\ii$ by right multiplication:
\begin{gather}
\label{right action}
R_\psi(\j_\infty\varphi|_z)=\j_\infty\big(\varphi\comp\psi^{-1}\big)|_{\psi(z)},
\end{gather}
whenever def\/ined.

The cotangent space $T^*\D^{(\infty)}$ naturally splits into horizontal and contact (groupoid) components.
The \emph{horizontal subbundle} is spanned by the right-invariant one-forms
\begin{gather*}
%\label{invariant horizontal coframe}
\sigma^a=d_M Z^a=\sum_{b=1}^m Z^a_b dz^b,
\qquad
a=1,\ldots,m,
\end{gather*}
while the \emph{contact subbundle} is spanned by the (right-invariant) \emph{Maurer--Cartan forms}
\begin{gather}
\label{mc forms}
\mu^a_B,
\qquad
a=1,\ldots,m,
\qquad
\#B\geq0.
\end{gather}
Their coordinate expressions are obtained by repeatedly applying the total derivative operators
\begin{gather}
\label{DZ}
\mathbb{D}_{Z^b}=\sum_{a=1}^m W^a_b \mathbb{D}_{z^a},
\qquad
\big(W^b_a\big)=\big(Z^b_a\big)^{-1}
\end{gather}
to the order zero Maurer--Cartan forms
\begin{gather*}
\mu^a=d_G Z^a=dZ^a-\sum_{b=1}^m Z^a_b dz^b
\end{gather*}
so that
\begin{gather*}
\mu^a_B=\mathbb{D}_Z^B\mu^a=\mathbb{D}_{Z^{b^1}}\cdots\mathbb{D}_{Z^{b^k}}\mu^a,
\qquad
k=\#B.
\end{gather*}
We refer to~\cite{OP-2005} for more details.

To concisely express the structure equations of the invariant coframe
$\{ \ldots \sigma^a \ldots \mu^a_B \ldots \}$, the vector-valued Maurer--Cartan formal power series
$\mu\llbracket H\rrbracket=(\mu^1\llbracket H\rrbracket,\ldots,\mu^m\llbracket H\rrbracket)^T$ with components
\begin{gather}
\label{mc power series}
\mu^a\llbracket H\rrbracket=\sum_{\#B\geq0} \frac{1}{B!} \mu_B^a H^B,
\qquad
a=1,\ldots,m,
\end{gather}
is introduced.
In the above formula, $H=(H^1,\ldots,H^m)$ are formal power series parameters while
$B!=\widetilde{b}^1! \widetilde{b}^2!\cdots\widetilde{b}^m!$ with $\widetilde{b}^a$ standing for the number of
occurrences of the integer $1\leq a \leq m$ in $B$.
The structure equations for the Maurer--Cartan forms are then obtained by comparing the coef\/f\/icients of the various
powers of $H$ in the power series identity
\begin{subequations}\label{diff structure eq}
\begin{gather}%\label{diff structure eq1}
d\mu\llbracket H\rrbracket=\nabla\mu\llbracket H\rrbracket\wedge(\mu\llbracket H\rrbracket-dZ),
\end{gather}
where $dZ=(dZ^1,\ldots,dZ^m)^T$ and $\nabla\mu\llbracket H\rrbracket=\left({\partial\mu^b\llbracket
H\rrbracket}/{\partial H^a}\right)$ denotes the $m\times m$ Jacobian matrix obtained by formally dif\/ferentiating the
power series~\eqref{mc power series} with respect to $H=(H^1,\ldots,H^m)$.
On the other hand, the structure equations for the invariant horizontal one-forms
$\sigma=(\sigma^1,\ldots$, $\sigma^m)^T$
are
\begin{gather}%\label{diff structure eq2}
d\sigma=\nabla\mu\llbracket0\rrbracket\wedge\sigma.
\end{gather}
\end{subequations}

\subsection{Lie pseudo-groups}

Several variants of the technical def\/inition of a~Lie pseudo-group exist in the
literature,~\cite{C-1953-1,GS-1966,J-1962,K-1975,K-1959,SS-1965}.
In the analytic category, Lie pseudo-groups can be def\/ined as follows.
\begin{Definition}
\label{Lie pseudo-group definition}
A pseudo-group $\G\subset\D$ is called a~\emph{Lie pseudo-group} of order $n_\star\geq1$ if, for all f\/inite $n\geq
n_\star\colon$
\begin{itemize}\itemsep=0pt
\itemsep=0pt \item $\G\n\subset\D\n$ forms a~smooth embedded subbundle, \item the projection
$\widetilde{\pi}^{n+1}_n\colon\G^{(n+1)}\to\G\n$ is a~f\/ibration, \item every local dif\/feomorphism $\phi\in\D$
satisfying $\j_{n_\star}\phi\subset\G^{(n_\star)}$ belongs to $\G$, \item $\G\n=\text{pr}^{(n-n_\star)}\G^{(n_\star)}$
is obtained by prolongation.
\end{itemize}
\end{Definition}

For $n\geq n_\star$, Def\/inition~\ref{Lie pseudo-group definition} implies that the pseudo-group jet subbundle
$\G^{(n)}\subset\D^{(n)}$ is characterized by a~formally integrable system of $n^\text{th}$ order dif\/ferential
equations
\begin{gather}
\label{determining system}
F^{(n)}\big(z,Z^{(n)}\big)=0,
\end{gather}
called the ($n^\text{th}$ order) \emph{determining system} of $\G\n$.

At the inf\/initesimal level, let
\begin{gather}
\label{vector field}
\vv=\sum_{a=1}^m\zeta^a(z)\pp{}{z^a}
\in
TM
\end{gather}
denote a~local vector f\/ield on $M$.
For $0\leq n\leq\infty$, let $\J^n TM$ denote the bundle of $n^\text{th}$ order jets of sections of $TM$ with local
coordinates
\begin{gather*}
(z,\zeta\n)=\big(\ldots z^a\ldots\zeta^a_B\ldots\big),
\end{gather*}
where $\zeta^a_B$ denotes the derivative coordinates of order $0\leq\#B\leq n$.
Let $\g$ denote the (local) Lie algebra of $\G$ consisting of all local vector f\/ields on $M$ tangent to the
pseudo-group orbits.
A vector f\/ield~\eqref{vector field} is in $\mathfrak{g}$ if and only if its jets satisfy the $n^{\text{th}}$ order
(formally integrable) \emph{infinitesimal determining system}
\begin{gather}
\label{infinitesimal determining system}
L\n(z,\zeta\n)=\sum_{a=1}^m\sum_{\#B\leq n} h^{B}_{a;\upsilon}(z) \zeta^a_B=0,
\qquad
\upsilon=1,\ldots,k,
\qquad
n\geq n_\star,
\end{gather}
obtained by linearizing the determining system~\eqref{determining system} at the identity jet $\mathds{1}\n$.
\begin{Theorem}
For each $n\geq n_\star$, the linear relations among the $($restricted$)$ Maurer--Cartan forms $\mu\n|_\G$ are obtained by
making the formal substitution or ``lift'' $($see~\eqref{lift of submanifold jets}, \eqref{lift of a vector field jet}
below$)$
\begin{gather*}
%\label{formal substitution}
z^a\longrightarrow Z^a,
\qquad
\zeta^a_B\longrightarrow\mu^a_B
\end{gather*}
in the infinitesimal determining equations~\eqref{infinitesimal determining system}:
\begin{gather}
\label{lifted infinitesimal determining equations}
L\n\big(Z,\mu\n\big)=0.
\end{gather}
\end{Theorem}

\begin{Corollary}
%\label{corollary MC structure eq}
The structure equations of a~Lie pseudo-group $\G$ are obtained by restricting the diffeomorphism structure
equations~\eqref{diff structure eq} to the solution space of~\eqref{lifted infinitesimal determining equations}:
\begin{gather}
\label{Lie pseudo-group structure equations}
\big(d\sigma=\nabla\mu\llbracket0\rrbracket\wedge\sigma,
\;
d\mu\llbracket H\rrbracket=
\nabla\mu\llbracket H\rrbracket\wedge(\mu\llbracket H\rrbracket-dZ)\big)\big|_{L\ii(Z,\mu\ii)=0}.
\end{gather}
\end{Corollary}

\begin{Example}
%\label{structure ODE example}
Let $M=\J^2(\mathbb{R}^2,1)$ be the second order jet bundle of curves in the plane with local coordinates
\begin{gather*}
x,
\qquad
u,
\qquad
p=u_x,
\qquad
q=u_{xx}.
\end{gather*}
To illustrate the constructions occurring in this paper we will consider the equivalence problem of second order
ordinary dif\/ferential equations
\begin{gather}
\label{ODE}
q=F(x,u,p)
\end{gather}
under the Lie pseudo-group of \emph{contact transformations}
\begin{gather}\label{contact}
X=\chi(x,u,p),
\qquad
U=\psi(x,u,p),
\qquad
P=\beta(x,u,p),
\qquad
Q=\frac{\beta_x+p\beta_u+q\beta_p}{\chi_x+p\chi_u+q\chi_p},
\end{gather}
where the functions $(\chi,\psi,\beta)\in\D(\mathbb{R}^3)$ satisfy the contact conditions
\begin{gather*}
%\label{contact transformations determining system}
\psi_p=\beta \chi_p,
\qquad
\beta(\chi_x+p\chi_u)=\psi_x+p \psi_u,
\end{gather*}
and the Lie pseudo-group of \emph{point transformations}
\begin{gather}
\label{point}
X=\chi(x,u),
\qquad
U=\psi(x,u),
\qquad
P=\frac{\widehat{D}\psi}{\widehat{D}\chi},
\qquad
Q=\frac{\widehat{D}^2\psi\cdot\widehat{D}\chi-\widehat{D}\psi\cdot\widehat{D}^2\chi}{(\widehat{D}\chi)^3},
\end{gather}
with $(\chi,\psi)\in\D(\mathbb{R}^2)$ and
\begin{gather*}
\widehat{D}=\pp{}{x}+p\pp{}{u}+q\pp{}{p}.
\end{gather*}

We now compute the low order structure equations for the contact pseudo-group~\eqref{contact}.
Let
\begin{gather}
\label{ODE vector field}
\vv=\xi(x,u,p,q)\pp{}{x}+\eta(x,u,p,q)\pp{}{u}+\alpha(x,u,p,q)\pp{}{p}+\gamma(x,u,p,q)\pp{}{q}
\end{gather}
denote a~local vector f\/ield on $M=\J^2(\mathbb{R}^2,1)$.
The vector f\/ield~\eqref{ODE vector field} is tangent to the orbits of the pseudo-group action~\eqref{contact} if
and only if
\begin{gather*}
\xi=\xi(x,u,p),
\qquad
\eta=\eta(x,u,p),
\qquad
\alpha=\alpha(x,u,p),
\qquad
\gamma=\widehat{D}\alpha-q \widehat{D}\xi,
\end{gather*}
and
\begin{gather*}
\eta_p=p \xi_p,
\qquad
\alpha=\eta_x+p(\eta_u-\xi_x)+p^2\xi_u.
\end{gather*}
Hence, the coef\/f\/icients of the inf\/initesimal generator~\eqref{ODE vector field} satisfy the determining system
\begin{gather}
\label{contact infinitesimal determining system}
\xi_q=\eta_q=\alpha_q=0,
\quad\;
\eta_p=p \xi_p,
\quad\;
\alpha=\eta_x+p(\eta_u-\xi_x)+p^2\xi_u,
\quad\;
\gamma=\widehat{D}\alpha-q \widehat{D}\xi.
\end{gather}
Under the replacement
\begin{gather*}
\xi_A\rightarrow\mu^x_A,
\qquad
\eta_A\rightarrow\mu^u_A,
\qquad
\alpha_A\rightarrow\mu^p_A,
\qquad
\gamma_A\rightarrow\mu^q_A,
\qquad
(x,u,p,q)\rightarrow(X,U,P,Q),
\end{gather*}
the inf\/initesimal determining equations~\eqref{contact infinitesimal determining system} yield the linear
dependencies
\begin{gather}
\mu^x_Q=\mu^u_Q=\mu^p_Q=0,
\qquad
\mu^u_P=P\mu^x_P,
\qquad
\mu^p=\mu^u_X+P\big(\mu^u_U-\mu^x_X\big)+P^2\mu^x_U,
\nonumber
\\
\mu^q=\mu^p_X+P\mu^p_U+Q\mu^p_P-Q\big(\mu^x_X+P\mu^x_U+Q\mu^x_P\big),
\label{contact lifted equations}
\end{gather}
among the Maurer--Cartan forms of order $\leq1$.
Dif\/ferentiating~\eqref{contact lifted equations} with respect to $\mathbb{D}_X$, $\mathbb{D}_Y$, $\mathbb{D}_P$,
$\mathbb{D}_Q$ as def\/ined in~\eqref{DZ}, we obtain the linear relations among the higher order Maurer--Cartan forms.
It follows that
\begin{gather*}
%\label{contact mc basis}
\mu_{X^iU^jP^k}=\mu^x_{X^i U^j P^k},
\qquad
\nu_{X^iU^j}=\mu^u_{X^iU^j},
\qquad
i,j,k\geq0,
\end{gather*}
is a~basis of Maurer--Cartan forms.
Restricting the structure equations of the dif\/feomorphism pseudo-group $\D(\mathbb{R}^4)$ to~\eqref{contact lifted
equations} and its prolongations we obtain the structure equations
\begin{gather}
d\sigma^x=-d\mu=\mu^{}_X\wedge\sigma^x+\mu^{}_U\wedge\sigma^u+\mu^{}_P\wedge\sigma^p,\nonumber
\\
d\sigma^u=-d\nu=\nu^{}_X\wedge\sigma^x+\nu^{}_U\wedge\sigma^u+P\mu_P\wedge\sigma^p,\nonumber
\\
d\sigma^p=
[\nu^{}_{XX}+P(\nu^{}_{UX}-\mu^{}_{XX})+P^2\mu^{}_{UX}]\wedge\sigma^x+[\nu^{}_{UX}+P(\nu^{}_{UU}-\mu^{}_{UX})+P^2\mu^{}_{UU}]\wedge\sigma^u\nonumber
\\
\phantom{d\sigma^p=}
{}+[\nu^{}_U-\mu^{}_X+P(2\mu^{}_U+\mu_{XP})+P(2P\mu_{UP}-\mu_{XP})]\wedge\sigma^p,\nonumber
\\
d\sigma^q=\mu^q_X\wedge\sigma^x+\mu^q_U\wedge\sigma^u+\mu^q_P\wedge\sigma^p+\mu^q_Q\wedge\sigma^q,
\label{horizontal contact structure equations}
\end{gather}
for the horizontal coframe and the order 0 Maurer--Cartan forms $\mu$, $\nu$.
We do not write the structure equations for the higher order Maurer--Cartan forms as these are not needed subsequently.

For the pseudo-group of point transformations~\eqref{point}, it suf\/f\/ices to add the constraints
\begin{gather*}
\xi_p=\eta_p=0
\end{gather*}
to~\eqref{contact infinitesimal determining system} to obtain the inf\/initesimal determining equations of its Lie
algebra:
\begin{gather}
\xi_p=\eta_p=0,
\qquad
\xi_q=\eta_q=\alpha_q=0,
\nonumber
\\
\alpha=\eta_x+p(\eta_u-\xi_x)+p^2\xi_u,
\qquad
\gamma=\widehat{D}\alpha-q(\xi_x+p\xi_u).
\label{point infinitesimal determining system}
\end{gather}
Taking the lift of~\eqref{point infinitesimal determining system} we obtain the linear relations
\begin{gather*}
\mu^x_P=\mu^u_P=0,
\qquad
\mu^x_Q=\mu^u_Q=\mu^p_Q=0,
\qquad
\mu^p=\mu^u_X+P\big(\mu^u_U-\mu^x_X\big)+P^2\mu^x_U,
\\
\mu^q=\mu^p_X+P\mu^p_U+Q\mu^p_P-Q\big(\mu^x_X+P\mu^x_U\big),
\end{gather*}
among the Maurer--Cartan forms of order $\leq1$.
A basis of Maurer--Cartan forms is thus given~by
\begin{gather}
\label{point basis mc forms}
\mu_{X^iU^j}=\mu^x_{X^iU^j},
\qquad
\nu_{X^iU^j}=\mu^u_{X^i U^j}.
\end{gather}
By setting $\mu_P=\mu_{PX}=\mu_{PU}=\cdots=0$ in~\eqref{horizontal contact structure equations} we obtain the structure
equations for the horizontal coframe $\{\sigma^x,\sigma^u,\sigma^p,\sigma^q\}$.
On the other hand, the structure equations for the Maurer--Cartan forms~\eqref{point basis mc forms} are given by the
structure equations of the dif\/feomorphism pseudo-group $\D(\mathbb{R}^2)$,~\cite{OP-2005}:
\begin{gather}
d\mu=\sigma^x\wedge\mu_X+\sigma^u\wedge\mu_U,\nonumber
\\
d\mu_X=\sigma^x\wedge\mu_{XX}+\sigma^u\wedge\mu_{XU}+\mu_U\wedge\nu_X,\nonumber
\\
d\mu_U=\sigma^x\wedge\mu_{XU}+\sigma^u\wedge\mu_{UU}+\mu_X\wedge\mu_U+\mu_U\wedge\nu_U,\nonumber
\\
%\hskip0.2cm\vdots
\cdots\cdots\cdots\cdots\cdots\cdots\cdots\cdots\cdots\cdots\cdots\cdots\cdots\cdots\cdots\cdots\cdots
%\cdots\cdots\cdots\cdots\cdots\cdots\cdots\cdots\cdots\cdots\cdots\cdots\cdots\cdots\cdots\cdots
\nonumber
\\
d\nu=\sigma^x\wedge\nu_X+\sigma^u\wedge\nu_U,\nonumber
\\
d\nu_X=\sigma^x\wedge\nu_{XX}+\sigma^u\wedge\nu_{XU}+\nu_X\wedge\mu_X+\nu_U\wedge\nu_X,\nonumber
\\
d\nu_U=\sigma^x\wedge\nu_{XU}+\sigma^u\wedge\nu_{UU}+\nu_X\wedge\mu_U,\nonumber
\\
d\nu_{UU}=
\sigma^x\wedge\nu_{XUU}+\sigma^u\wedge\nu_{UUU}+2\nu_{XU}\wedge\nu_U+\nu_X\wedge\mu_{UU}+\nu_{UU}\wedge\nu_U,\nonumber
\\
d\nu_{XU}=
\sigma^x\!\wedge\nu_{XXU}+\sigma^u\!\wedge\nu_{XUU}\!+\nu_{XU}\wedge\mu_X+\nu_{XX}\wedge\mu_U+\nu_X\wedge\mu_{XU}+\nu_{UU}\wedge\nu_X,\nonumber
\\
%\hskip0.2cm\vdots
\cdots\cdots\cdots\cdots\cdots\cdots\cdots\cdots\cdots\cdots\cdots\cdots\cdots\cdots\cdots\cdots\cdots
\label{mc structure equations}
\end{gather}
\end{Example}

\section{Equivariant moving frames}
\label{Section Moving Frame}

As in the previous section, let $\G$ be a~Lie pseudo-group acting on an $m$-dimensional manifold~$M$.
We are now interested in the induced action of~$\G$ on $p$-dimensional submanifolds $S\subset M$, where $1\leq p<m$.
For each integer $0\leq n\leq\infty$, let $\J^n=\J^n(M,p)$ denote the $n^\text{th}$ order \emph{submanifold jet
bundle}, def\/ined as the set of equivalence classes under the equivalence relation of $n^\text{th}$ order
contact,~\cite{O-1995}.
For $k\geq n$, we use $\pi^k_n\colon\J^k\to\J^n$ to denote the canonical projection.
We introduce local coordinates $z=(x,u)=(x^1,\ldots,x^p,u^1,\ldots,u^q)$ on $M$ so that submanifolds that are
transverse to the vertical f\/ibers $\{x=x_0\}$ are (locally) given as graphs of smooth functions $u=f(x)$.
(Submanifolds with vertical tangent spaces are handled by a~dif\/ferent choice of local coordinates.)
In this adapted system of coordinates, the standard coordinates on $\Jn$ are
\begin{gather*}
z\n=\big(x,u\n\big)=\big( \ldots x^i \ldots u^\alpha_J \ldots \big),
\end{gather*}
where $u^\alpha_J$ denotes the derivative coordinates of order $0\leq\#J\leq n$.

Let $\mathcal{B}\n\to\J^n$ denote the $n^{\text{th}}$ order \emph{lifted bundle} obtained by pulling back $\G\n\to M$
via the projection $\pi^n_0\colon\J^n\to M$.
Local coordinates on $\mathcal{B}\n$ are given by $\big(z\n,g\n\big)$, where the base coordinates $z\n=\big(x,u\n\big)\in\J^n$ are the
submanifold jet coordinates and the f\/iber coordinates are the pseudo-group parameters $g\n$.
The bundle $\mathcal{B}\n$ carries the structure of a~groupoid with \emph{source map} $\bs\n\big(z\n,g\n\big)=z\n$ and
\emph{target map} $\bt\n\big(z\n,g\n\big)=Z\n=g\n\cdot z\n$ given by the \emph{prolonged action}.
To compute the prolonged action we introduce on $\mathcal{B}\ii$ the \emph{lifted horizontal coframe}
\begin{gather*}
%\label{lifted horizontal coframe}
\omega^i=d_H X^i=\sum_{j=1}^p \big(D_{x^j}X^i\big) dx^j,
\qquad
i=1,\ldots,p,
\end{gather*}
where $D_{x^1},\ldots,D_{x^p}$ are the total derivative operators.
The \emph{lifted total derivative operators} $D_{X^1},\ldots,D_{X^p}$ are then def\/ined by the formula
\begin{gather*}
d_H F\big(z\n\big)=\sum_{i=1}^p (D_{x^i}F) dx^i=\sum_{i=1}^p (D_{X^i}F) \omega^i.
\end{gather*}
More explicitly,
\begin{gather}
\label{lifted total differential operators}
D_{X^i}=\sum_{j=1}^p \widehat{W}^j_i D_{x^j},
\qquad
\text{with}
\qquad
\big(\widehat{W}^j_i\big)=\big(D_{x^i}X^j\big)^{-1}.
\end{gather}
Dif\/ferentiating the target dependent variables $U^\alpha$ with respect to the lifted total derivative
operators~\eqref{lifted total differential operators} we obtain the explicit expressions for the prolonged action
$Z\n=g\n\cdot z\n$:
\begin{gather}
\label{prolonged action}
X^i,
\qquad
U^\alpha_J=D^J_X U^\alpha.
\end{gather}

A local dif\/feomorphism $h\in\G$ acts on $\{\big(z\n,g\n\big) \in \mathcal{B}\n \, |\,\pi^n_0\big(z\n\big)\in\dom h\}$ by right
multi\-pli\-cation:
\begin{gather}
\label{lifted action}
R_h\big(z\n,g\n\big)=\big(h\n\cdot z\n,g\n\cdot\big(h\n\big)^{-1}\big).
\end{gather}
The ($n^\text{th}$ order) \emph{lifted action} action~\eqref{lifted action} is given by the concatenation of the
prolonged action on submanifold jets with the restriction of the right action~\eqref{right action} to $\G\n$.
The expressions~\eqref{prolonged action} are invariant under the lifted action~\eqref{lifted action} and these
functions are called \emph{lifted invariants}.

\subsection{Regular submanifold jets}
\label{subsection regular jets}

The existence of a~moving frame requires the prolonged pseudo-group action on submanifold jets to be (locally) free and
regular,~\cite{OP-2008}.
\begin{Definition}
%\label{regular action}
A pseudo-group acts regularly in the neighborhood of a~point $z\n\in\J^n$ if the pseudo-group orbits have the same
dimension and there are arbitrarily small neighborhoods whose intersection with each orbit is a~connected subset
thereof.
\end{Definition}
\begin{Definition}
%\label{definition free action}
The \emph{isotropy subgroup} of a~submanifold jet $z\n\in\J^n$ is def\/ined as
\begin{gather*}
\G_{z\n}\n=\big\{g\n\in\G\n|_z\,:\,g\n\cdot z\n=z\n\big\},
\end{gather*}
where $\pi^n_0\big(z\n\big)=z$.
The pseudo-group is said to act \emph{freely} at $z\n$ if $\G\n_{z\n}=\{\mathds{1}\n|_z\}$.
The pseudo-group acts \emph{locally freely} at $z\n$ if $\G\n_{z\n}$ is discrete.
\end{Definition}
\begin{Definition}
A submanifold jet $z\ii$ is said to be \emph{regular} if there exists a~f\/inite $n\geq1$ such that the pseudo-group
$\G$ acts freely at $z\n=\pi^\infty_n\big(z\ii\big)$.
Let $\mathcal{V}^\infty\subset\J^\infty$ denote the subset of all regular submanifold jets.
\end{Definition}

Following the foundational papers~\cite{OP-2009-1, OP-2008} we, for the moment, assume that the pseudo-group~$\G$ acts
regularly on $\mathcal{V}^\infty$ and that this set is a~dense open subbundle of~$\J^\infty$.
In Example~\ref{regular moving frame example}, we will see that, in general, these assumptions need to be relaxed.
\begin{Definition}
%\label{moving frame definition}
Let $\mathcal{V}^n=\pi^\infty_n(\mathcal{V}^\infty)$ denote the truncation of the regular submanifold jets to order $n$.
A (right) \emph{moving frame} of \emph{order} $n$ is a~(local) $\G$-equivariant section
$\widehat{\rho}^{ (n)}\colon\mathcal{V}^n\to\mathcal{B}\n$.
\end{Definition}

\noindent In local coordinates we use the notation
\begin{gather*}
\widehat{\rho}^{ (n)}\big(z\n\big)=\big(z\n,\rho\n\big(z\n\big)\big)
\end{gather*}
to denote a~moving frame.
Right equivariance means that
\begin{gather*}
R_g\widehat{\rho}^{ (n)}\big(z\n\big)=\widehat{\rho}^{ (n)}\big(g\n\cdot z\n\big),
\qquad
g\in\G,
\qquad
\text{when def\/ined.}
\end{gather*}
\begin{Theorem}
%\label{existence moving frame}
Suppose $\G$ acts $($locally$)$ freely on $\mathcal{V}^n\subset\J^n$, with its orbits forming a~regular foliation, then an
$n^\text{th}$ order moving frame exists in a~neighborhood of every $z\n\in\mathcal{V}^n$.
\end{Theorem}

Once a~pseudo-group action is free, a~result known as the \emph{persistence of freeness},~\cite{OP-2009-1,OP-2009-2},
guarantees that the action remains free under prolongation.
\begin{Theorem}
%\label{persistence of freeness theorem}
If a~Lie pseudo-group $\G$ acts $($locally$)$ freely at $z\n$ then it acts $($locally$)$ freely at any $z^{(k)}\in\J^k$, $k\geq
n$, with $\pi^k_n(z^{(k)})=z\n$.
The minimal $n$ such that the action becomes free is called the \emph{order of freeness} and is denoted by $n^\star$.
\end{Theorem}

An order $n\geq n^\star$ moving frame is constructed through a~normalization procedure based on the choice of
a~cross-section $\mathcal{K}^n\subset\mathcal{V}^n$ to the pseudo-group orbits.
The associated (locally def\/ined) right moving frame section
$\widehat{\rho}^{ (n)}\colon\mathcal{V}^n\to\mathcal{B}\n$ is uniquely characterized by the condition that
$\bt\n\big(\widehat{\rho}^{ (n)}\big(z\n\big)\big)=\rho\n\big(z\n\big)\cdot z\n\in\mathcal{K}^n$.
For simplicity, we assume that $\mathcal{K}^n=\{z_{i_1}=c_1,\ldots,z_{i_{r_n}}=c_{r_n}: r_n=\dim \G\n|_z\}$ is
a~coordinate cross-section.
Then, the moving frame $\widehat{\rho}^{ (n)}$ is obtained by solving the \emph{normalization equations}
\begin{gather}
\label{normalization equations}
Z_{i_1}\big(x,u\n,g\n\big)=c_1,
\qquad
\ldots
\qquad
Z_{i_{r_n}}\big(x,u\n,g\n\big)=c_{r_n},
\end{gather}
for the pseudo-group parameters $g\n=\rho\n\big(x,u\n\big)$.
The invariants appearing on the left-hand side of the normalization equations~\eqref{normalization equations} are
called \emph{phantom invariants}.
As one increases the order from $n$ to $k>n$, a~new cross-section must be selected.
We require that these cross-sections be compatible in the sense that $\pi^k_n(\mathcal{K}^k)=\mathcal{K}^n$       %  \big)  ?
for all $k>n$.
This in turn, implies the compatibility of the moving frames:
$\pi^k_n\big(\widehat{\rho}^{ (k)}\big(z^{(k)}\big)\big)=\widehat{\rho}^{ (n)}\big(\pi^k_n\big(z^{(k)}\big)\big)$.
A compatible sequence of moving frames is simply called a~moving frame and is denoted by
$\widehat{\rho}\colon\mathcal{V}^\infty\to\mathcal{B}^{(\infty)}$.
Finally, we require the compatible cross-sections to be of \emph{minimal order},~\cite{OP-2009-1}.
Intuitively, this is equivalent to requiring that the pseudo-group parameters be normalized as soon as possible during
the normalization procedure.
\begin{figure}[!h]\centering
\includegraphics[scale=0.9]{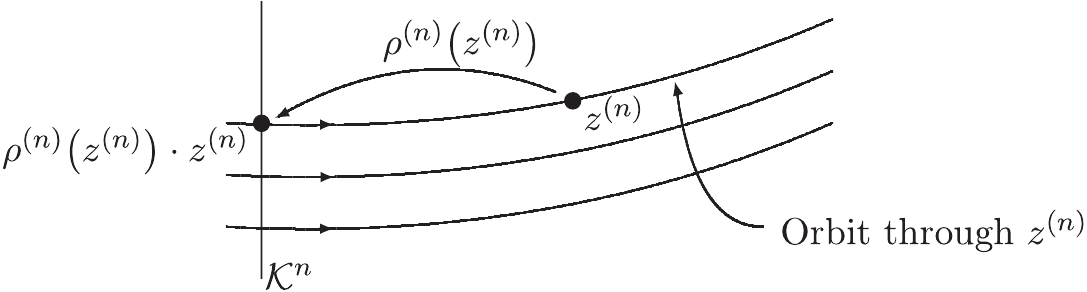}
%\begin{center}
%\begin{picture}
%(200,80) \qbezier(0,35)(95,30)(175,65) \qbezier(0,20)(95,15)(175,50) \qbezier(0,50)(95,45)(175,80)
%\put(30,34.4){\vector(1,0){1}} \put(30,49.6){\vector(1,0){1}} \put(30,19.4){\vector(1,0){1}} \put(10,5){\line(0,1){80}}
%\put(11,2){$\mathcal{K}^n$} \put(100,56.5){\circle*{5}} \put(103,48){$z\n$} \qbezier(130,60)(135,20)(155,20)
%\put(130.5,58){\vector(-1,4){1}} \put(160,15){Orbit through $z\n$} \put(-65,39){$\rho\n\big(z\n\big)\cdot z\n$}
%\put(10,50){\circle*{5}} \qbezier(15,52)(55,75)(95,59) \put(15,52){\vector(-2,-1){1}} \put(45,70){$\rho\n\big(z\n\big)$}
%\end{picture}
\caption{Moving frame $\widehat{\rho}^{ (n)}\big(z\n\big)=\big(z\n,\rho\n\big(z\n\big)\big)$.}
%\end{center}
\end{figure}

We now introduce the invariantization map associated with a~moving frame.
First, we note that the space of dif\/ferential forms on $\mathcal{B}^{(\infty)}$ splits into
\begin{gather*}
\bO^*=\bigoplus_{k,l}\bO^{k,l}=\bigoplus_{i,j,l}\bO^{i,j,l},
\end{gather*}
where $l$ indicates the number of Maurer--Cartan forms~\eqref{mc forms}, $k=i+j$ the number of \emph{jet forms}, with
$i$ indicating the number of horizontal forms $dx^i$ and $j$ the number of basic submanifold jet contact forms
\begin{gather*}
%\label{basic contact forms}
\theta^\alpha_J=du^\alpha_J-\sum_{i=1}^p u^\alpha_{J,i} dx^i,
\qquad
\alpha=1,\ldots,q,
\qquad
\#J\geq0.
\end{gather*}
Next, let
\begin{gather*}
%\label{jet differential splitting}
\bO^*_J=\bigoplus_k\bO^{k,0}=\bigoplus_{i,j}\bO^{i,j,0}
\end{gather*}
denote the subspace of jet forms consisting of those dif\/ferential forms containing no Maurer--Cartan forms.
Then, we introduce the projection $\pi_J\colon\bO^*\to\bO^*_J$ which takes a~dif\/ferential form~$\Omega$ on~$\mathcal{B}^{(\infty)}$ to its jet component~$\pi_J(\Omega)$ obtained by annihilating the Maurer--Cartan forms in~$\Omega$.
\begin{Definition}
%\label{definition lift}
Let $\Omega$ be a~dif\/ferential form def\/ined on $\J^\infty$.
Its \emph{lift} is the invariant jet form
\begin{gather}
\label{lift}
\bl(\Omega)=\pi_J\big[\big(\bt\ii\big)^*(\Omega)\big]
\end{gather}
def\/ined on the lifted bundle $\mathcal{B}\ii$.
\end{Definition}

In particular, setting $\Omega$ in~\eqref{lift} to be each of the submanifold jet coordinates $x^i$, $u^\alpha_J$, the
lift map~\eqref{lift} coincides with the prolonged action~\eqref{prolonged action}:
\begin{gather}
\label{lift of submanifold jets}
\bl\big(x^i\big)=X^i,
\qquad
\bl\big(u^\alpha_J\big)=U^\alpha_J.
\end{gather}
Also, we note that the lift of the horizontal forms $dx^1, \ldots, dx^p$
\begin{gather*}
%\label{lifted horizontal forms}
\bl\big(dx^i\big)=\sum_{i=1}^p \big(D_{x^j}X^i\big) dx^j+\sum_{\alpha=1}^q X^i_{u^\alpha} \theta^\alpha=
\omega^i+\sum_{\alpha=1}^q X^i_{u^\alpha} \theta^\alpha
\end{gather*}
are invariant horizontal forms if and only if the pseudo-group action is projectable meaning that $X^i_{u^\alpha}=0$.
On the other hand, the lift of a~contact form is always a~contact form.
In the following, we ignore contact forms and introduce the equivalence relation $\equiv$ to indicate equality modulo
contact forms.
\begin{Definition}
Let $\widehat{\rho}\colon\mathcal{V}^\infty\to\mathcal{B}^{(\infty)}$ be a~moving frame.
If $\Omega$ is a~dif\/ferential form on $\mathcal{V}^\infty$, then its \emph{invariantization} is the invariant
dif\/ferential form
\begin{gather}
\label{invariantization map}
\iota(\Omega)=\widehat{\rho}^{\,*}[\bl(\Omega)].
\end{gather}
\end{Definition}

In the following, to lighten the notation, we will usually omit writing moving frame pull-backs:
\begin{gather*}
\omega^j=\widehat{\rho}^{\,*}(\omega^j)\equiv\widehat{\rho}^{\,*}\big[\bl\big(dx^j\big)\big]=\iota\big(dx^j\big),
\\
\big(X,U\ii\big)=\widehat{\rho}^{\,*}\big(X,U\ii\big)=\widehat{\rho}^{\,*}\big[\bl\big(x,u\ii\big)\big]=\iota\big(x,u\ii\big).
\end{gather*}
\begin{Proposition}
%\label{normalized invariants proposition}
The normalized differential invariants
\begin{gather*}
%\label{normalized invariants}
X^i=\iota\big(x^i\big),
\qquad
U^\alpha_J=\iota\big(u^\alpha_J\big),
\end{gather*}
contain a~complete set of functionally independent differential invariants.
\end{Proposition}

One of the most important results in the theory of equivariant moving frames is the \emph{universal recurrence formula}
for lifted/invariantized dif\/ferential forms,~\cite{OP-2008}.
To write down the formula we must extend the lift map~\eqref{lift} to vector f\/ield jet coordinates.
\begin{Definition}
%\label{definition lift mc forms}
The lift of a~vector jet coordinate $\zeta^a_B$ is def\/ined to be the Maurer--Cartan form $\mu^a_B$:
\begin{gather}
\label{lift of a vector field jet}
\bl\big(\zeta^a_B\big)=\mu^a_B,
\qquad
\text{for}
\quad
a=1,\ldots,m,
\quad
\#B\geq0.
\end{gather}
\end{Definition}

\noindent More generally, the lift of any f\/inite linear combination of vector f\/ield jet coordinates is
\begin{gather*}
\bl\left[\sum_{a=1}^m\sum_{\#B\geq0} P^B_a\big(x,u\n\big) \zeta^a_B\right]=\sum_{a=1}^m\sum_{\#B\geq0} P^B_a\big(X,U\n\big) \mu^a_B.
\end{gather*}

Recall that if
\begin{gather}
\label{xu infinitesimal generator}
\vv=\sum_{i=1}^p \xi^i(x,u)\pp{}{x^i}+\sum_{\alpha=1}^q \phi_\alpha(x,u)\pp{}{u^\alpha} \in \g
\end{gather}
is an inf\/initesimal generator of the pseudo-group action, then its \emph{prolongation} is the vector f\/ield
\begin{gather}
\vv^{(\infty)}=\sum_{i=1}^p \xi^i(x,u)\pp{}{x^i}+\sum_{\alpha=
1}^q\sum_{\#J\geq0} \phi^J_\alpha\pp{}{u^\alpha_J} \in \mathfrak{g}^{(\infty)},
\label{prolonged vector field}
\end{gather}
where
\begin{gather*}
%\label{prolonged coefficients}
\phi_\alpha^J=D_J Q^\alpha+\sum_{i=1}^p \xi^i u^\alpha_{J,i},
\qquad
\text{and}
\qquad
Q^\alpha=\phi_\alpha-\sum_{i=1}^p\xi^i u^\alpha_i
\end{gather*}
are the \emph{characteristic} components of the vector f\/ield~\eqref{xu infinitesimal generator}.
Then, for $z\ii\in\J^\infty|_z$ the prolongation formula~\eqref{prolonged vector field} def\/ines the
\emph{prolongation map}
\begin{gather}
\label{prolongation map}
\bp=\bp_{z\ii}\ii\colon \ \J^\infty TM|_z\to T\J^\infty|_{z\ii},
\qquad
\bp\ii_{z\ii}(\j_\infty\vv|_z)=\vv\ii|_{z\ii}.
\end{gather}
\begin{Theorem}
Let $\Omega$ be a~differential form on $\J^\infty$.
Then
\begin{gather}
\label{lifted recurrence}
d[\bl(\Omega)]=\bl\big[d\Omega+\vv^{(\infty)}(\Omega)\big],
\qquad
\text{where}
\qquad
\vv\ii\in\g\ii
\end{gather}
and $\vv\ii(\Omega)$ denotes the Lie derivative of $\Omega$ along $\vv\ii$.
\end{Theorem}
We refer to~\cite{OP-2008} for a~proof of~\eqref{lifted recurrence}.
In particular, the identity~\eqref{lifted recurrence} applies to the lifted dif\/ferential invariants $X^i$,
$U^\alpha_J$:
\begin{gather}
dX^i\equiv \omega^i+\mu^i,i=1,\ldots,p,
\nonumber
\\
d U^\alpha_J\equiv\sum_{j=1}^p U^\alpha_{J,j} \omega^j+\bl(\phi^J_\alpha),
\qquad
\alpha=1,\ldots,q,
\qquad
\#J\geq0.
\label{jet coordinate recurrence relations}
\end{gather}
\begin{Corollary}
Let $\widehat{\rho}\colon\mathcal{V}^\infty\to\mathcal{B}^{(\infty)}$ be a~moving frame and $\Omega$ a~differential
form on $\mathcal{V}^\infty$, then
\begin{gather}
\label{universal recurrence relations}
d[\iota(\Omega)]=\iota\big[d\Omega+\vv^{(\infty)}(\Omega)\big],
\qquad
\text{where}
\qquad
\vv\ii\in\g\ii.
\end{gather}
\end{Corollary}

Of particular interest to us is when $\Omega$ is one of the submanifold jet coordinate func\-tions~$x^i$,~$u^\alpha_J$.
Equation~\eqref{universal recurrence relations} then produces the recurrence relations
\begin{gather}
\label{horizontal recurrence relations}
\D_j X^i=\delta_j^i+N^i_j,
\qquad
\D_j U^\alpha_J=U^\alpha_{J,j}+M^\alpha_{J,j},
\end{gather}
where $N^i_j$ and $M^\alpha_{J,j}$ are correction terms and
\begin{gather}
\label{invariant total derivative operators}
\D_k=\sum_{j=1}^p \widehat{\rho}^{\,*}\big(W^j_k\big) D_{x^j},
\qquad
k=1,\ldots,p,
\end{gather}
are invariant total derivative operators obtained by normalizing the lifted total derivative opera\-tors~\eqref{lifted total differential operators}.

The commutation relations among the invariant total derivative operators $\D_1,\ldots,\D_p$ can be deduced from the
universal recurrence relation~\eqref{universal recurrence relations}.
By setting $\Omega=dx^i$ in~\eqref{universal recurrence relations} we obtain the equations
\begin{gather*}
%\label{structure equation invariant coframe}
d\omega^i\equiv-\sum_{1\leq j<k\leq p}Y^i_{jk} \omega^j\wedge\omega^k,
\qquad
i=1,\ldots,p.
\end{gather*}
Since the operators $\D_i$ are dual to the contact invariant horizontal forms $\omega^i$, it follows that
\begin{gather}
\label{commutation relations}
[\D_i,\D_j]=\sum_{k=1}^p Y^k_{ij} \D_k,
\qquad
1\leq i,j\leq p,
\end{gather}
and the invariant functions $Y^i_{jk}$ are called \emph{commutator invariants},~\cite{OP-2009-1}.
An important feature of the recurrence formula~\eqref{universal recurrence relations} (or~\eqref{horizontal recurrence
relations}) is that the coordinate expressions for the invariants $(X,U^{(\infty)})$, the contact invariant horizontal
1-forms $\omega^i$, the Maurer--Cartan forms $\mu^a_B$ and the moving frame $\widehat{\rho}$ are not required to
compute these equations.
One only needs to know the cross-section $\mathcal{K}^\infty\subset\mathcal{V}^\infty$ def\/ining $\widehat{\rho}$ and
the expression of $\vv^{(\infty)}\in\mathfrak{g}^{(\infty)}$.
The key observation is that the pulled-back Maurer--Cartan forms $\widehat{\rho}^{\,*}(\mu^a_B)$ can be obtained from
the phantom invariant recurrence relations.
We refer the reader to~\cite{COP-2008,OP-2008,OP-2009-2} for concrete examples of the moving frame implementation.

\begin{Example}
\label{regular moving frame example}
As mentioned at the beginning of this section, the theory introduced above assumes the Lie pseudo-group action to be
free and regular on a~dense open subset $\mathcal{V}^\infty\subset\J^\infty$.
Using the local equivalence problem of second order ordinary dif\/ferential equations under the pseudo-group of point
transformations~\eqref{point}, we now show that, in general, these assumptions should be relaxed.
Working symbolically, we use the recurrence relations~\eqref{jet coordinate recurrence relations} to f\/ind disjoint
sets of regular submanifold jets each admitting their own moving frame.

The f\/irst step consists of determining the ``universal normalizations'' which hold for any second order ordinary
dif\/ferential equation.
Beginning with the order 0 recurrence relations
\begin{gather}
dX\equiv\omega^x+\mu,
\nonumber
\\
dU\equiv\omega^u+\nu,
\nonumber
\\
dP\equiv\omega^p+\nu_X+P(2\nu_Y-\mu_X)-P^2\mu_U,
\nonumber
\\
dQ\equiv Q_P \omega^p+Q_U \omega^u+Q_X \omega^x+\nu_{XX}+Q(\nu_U-2\mu_X)+P(2\nu_{XU}-\mu_{XX})
\nonumber
\\
\phantom{dQ\equiv}
{}-3PQ\mu_U+P^2(\nu_{UU}-2\mu_{XU})-P^3\mu_{UU},
\label{order 0 recurrence relations}
\end{gather}
the lone appearance of the linearly independent Maurer--Cartan forms $\mu$, $\nu$, $\nu_X$, $\nu_{XX}$ in the group
dif\/ferential component of the recurrence relations~\eqref{order 0 recurrence relations} implies that we can normalize
the lifted invariants
\begin{gather}
\label{order 0 normalizations}
X=U=P=Q=0
\end{gather}
to zero, independently of the dif\/ferential equations.
Substituting~\eqref{order 0 normalizations} into~\eqref{order 0 recurrence relations} we obtain a~system of equations
that can be solved for the (partially) normalized Maurer--Cartan forms
\begin{gather*}
\mu\equiv-\omega^x,
\qquad
\nu\equiv-\omega^u,
\qquad
\nu_X\equiv-\omega^p,
\qquad
\nu_{XX}\equiv-\big(Q_P \omega^p+Q_U \omega^u+Q_X \omega^x\big).
\end{gather*}
Continuing the normalization procedure, order by order, we come to the conclusion that it is always possible to
normalize the lifted invariants
\begin{gather}
X=U=P=0,
\nonumber
\\
Q_{U^jX^k}=Q_{P U^j X^k}=Q_{P^2U^j}=Q_{P^2U^j X}=Q_{P^3U^j}=Q_{P^3U^jX}=0,
\qquad
j,k\geq0,
\label{intermediate normalizations}
\end{gather}
to zero (see~\cite[Appendix~B]{MV} for more details).
The normalizations~\eqref{intermediate normalizations} lead to the normali\-za\-tion of all the Maurer--Cartan
forms~\eqref{point basis mc forms} except for
\begin{gather}
\label{unnormalized mc forms}
\mu_X,
\qquad
\mu_U,
\qquad
\nu_U,
\qquad
\nu_{UU},
\qquad
\nu_{XU}.
\end{gather}
To proceed further, the value of the remaining (partially normalized) lifted invariants
\begin{gather}
\label{intermediate remaining invariants}
Q_{P^{k+4}U^j X^i},
\qquad
Q_{P^3U^j X^{i+2}},
\qquad
Q_{P^2U^j X^{i+2}},
\qquad
i,j,k\geq0,
\end{gather}
must be carefully analyzed.
More explicitly, the invariants~\eqref{intermediate remaining invariants} of order $\leq6$ are
\begin{gather}
n=4:
\quad
Q_{P^4},\ Q_{P^2X^2},
\nonumber
\\
n=5:
\quad
Q_{P^5},\ Q_{P^4U},\ Q_{P^4X},\ Q_{P^3X^2},\ Q_{P^2U X^2},\ Q_{P^2X^3},
\nonumber
\\
n=6:
\quad
Q_{P^6},\ Q_{P^5U},\ Q_{P^5X},\ Q_{P^4U^2},\ Q_{P^4U X},\ Q_{P^4X^2},\ Q_{P^3U X^2},
\nonumber
\\
\phantom{n=6:} \
\quad
Q_{P^3X^3},\ Q_{P^2U^2X^2},\ Q_{P^2U X^3},\ Q_{P^2X^4}.
\label{456 invariants}
\end{gather}
Writing the recurrence relations for the invariants~\eqref{456 invariants} of order $\leq5$, taking into account the
normalizations~\eqref{intermediate normalizations}, we obtain
\begin{gather}
dQ_{P^4}\equiv Q_{P^5} \omega^p+Q_{P^4U} \omega^u+Q_{P^4X} \omega^x+Q_{P^4}(2\mu_X-3\nu_U),\nonumber \\
dQ_{P^2X^2}\equiv Q_{P^3X^2} \omega^p+Q_{P^2U X^2} \omega^u+Q_{P^2X^3} \omega^x-Q_{P^2X^2}(\nu_U+2\mu_X),\nonumber \\
dQ_{P^5}\equiv Q_{P^6} \omega^p+Q_{P^5U} \omega^u+Q_{P^5X} \omega^x+5Q_{P^4} \mu_U+Q_{P^5}(3\mu_X-4\nu_U),\nonumber \\
dQ_{P^4X}\equiv(Q_{P^5X}+Q_{P^4U}) \omega^p+Q_{P^4U X} \omega^u+Q_{P^4X^2} \omega^x+Q_{P^4} \nu_{UX}
 +Q_{P^4X}(\mu_X-3\nu_U),\nonumber \\
dQ_{P^4U}\equiv Q_{P^5U} \omega^p+Q_{P^4U^2} \omega^u+Q_{P^4U X} \omega^x-2Q_{P^4} \nu_{UU}-Q_{P^5} \nu_{UX}-Q_{P^4X} \mu_U\nonumber \\
\phantom{dQ_{P^4U}\equiv}
{} +Q_{P^4U}(2\mu_X-4\nu_U),\nonumber \\
dQ_{P^3X^2}\equiv Q_{P^4X^2} \omega^p+Q_{P^3U X^2} \omega^u+(Q_{P^3X^3}-2Q_{P^2U X^2}) \omega^x-Q_{P^2X^2} \mu_U\nonumber \\
\phantom{dQ_{P^3X^2}\equiv}
{} -Q_{P^3X^2}(2\nu_U+\mu_X),\nonumber \\
dQ_{P^2U X^2}\equiv Q_{P^3U X^2} \omega^p+Q_{P^2U^2X^2} \omega^u+Q_{P^2U X^3} \omega^x-2Q_{P^2X^2} \nu_{UU}-Q_{P^3X^2} \nu_{UX}\nonumber \\
\phantom{dQ_{P^2U X^2}\equiv}
{} -Q_{P^2X^3} \mu_U-2Q_{P^2U X^2}(\nu_U+\mu_X),\nonumber \\
dQ_{P^2X^3}=(Q_{P^3X^3}-Q_{P^2UX^2}) \omega^p+Q_{P^2U X^3} \omega^u+Q_{P^2X^4} \omega^x-5Q_{P^2X^2} \nu_{UX}\nonumber \\
\phantom{dQ_{P^2X^2}\equiv}
{} -Q_{P^2X^3}(\nu_U+3\mu_X).
\label{45 recurrence relations}
\end{gather}
At this juncture, the normalization procedure splits into dif\/ferent branches depending on the value of the
fourth-order lifted invariants
\begin{gather*}
Q_{P^4}=\frac{\chi_x^2}{\psi_u^3}q_{pppp},
\qquad
Q_{P^2X^2}=
\frac{\widehat{D}^2(q_{pp})-4\widehat{D}(q_{u p})-q_{p}\widehat{D}(q_{p p})+6q_{uu}-3q_u q_{p p}+4q_{p}q_{u p}}{\psi_u\chi_x^2}
\end{gather*}
obtained by implementing the moving frame construction.
There are 4 branches to consider\footnote{To distinguish between lifted invariants that are set equal to a~constant by
normalization from those that are identically constant, we use the notation $=$ and $\equiv$ respectively.}
\begin{gather*}
\hphantom{I}{\bf I)}~Q_{P^4}\not\equiv0 \ \text{and} \ Q_{P^2X^2}\not\equiv0, \qquad {\bf III)}~Q_{P^4}\not\equiv0 \ \text{and} \ Q_{P^2X^2}\equiv0,
\\
{\bf II)}~Q_{P^4}\equiv0 \ \text{and} \  Q_{P^2X^2}\not\equiv0, \qquad {\bf IV)}~Q_{P^4}\equiv0 \ \text{and} \ Q_{P^2X^2}\equiv0.
\end{gather*}
In Case {\bf IV)} the pseudo-group action is not free and we dif\/fer this branch of the
equivalence problem to the next section.
For the other three cases, let
\begin{gather*}
%\label{intermediate cross-section}
\mathcal{K}^\infty=\{x=u=p=q_{u^j x^k}=q_{pu^jx^k}=q_{p^2u^j}=q_{p^2u^jx}=q_{p^3u^j}=q_{p^3u^jx}=0\,|\, j,k\geq0\}
\end{gather*}
denote the ``intermediate cross-section'' corresponding to the universal normalizations~\eqref{intermediate
normalizations}.
In Case~{\bf I)}, the action is free on the subbundle of regular submanifold jets
\begin{gather*}
\mathcal{V}^\infty_1=\big\{z\ii\in\J^\infty\,|\,Q_{P^4}\not\equiv0,\; Q_{P^2X^2}\not\equiv0\big\}.
\end{gather*}
From the recurrence relations~\eqref{45 recurrence relations}, we f\/ind that a~cross-section is given by
\begin{gather}
\label{K1}
\mathcal{K}^\infty_1=\mathcal{K}^\infty\cup\{q_{p^4}=q_{p^2x^2}=1,\;q_{p^5}=q_{p^4u}=q_{p^4x}=
0\}\subset\mathcal{V}^\infty_1
\end{gather}
leading to a~moving frame $\widehat{\rho}_1\colon\mathcal{V}^\infty_1\to\mathcal{B}\ii$.
In Case {\bf II)}, it can be shown, see~\cite{MV} for more details, that when $Q_{P^2X^2}\not\equiv0$ the (partially
normalized) lifted invariants $Q_{P^3X^3}$ and $Q_{P^2X^4}$ cannot be simultaneously equal to zero.
There are then 2 subbundles of regular submanifold jets
\begin{gather*}
\mathcal{V}^\infty_2=\big\{z\ii\in\J^\infty\,|\,Q_{P^2X^2}\not\equiv0,\;D_JQ_{P^4}\equiv0,\;Q_{P^2X^4}\not\equiv0\big\},
\\
\mathcal{V}^\infty_3=\big\{z\ii\in\J^\infty\,|\, Q_{P^2X^2}\not\equiv0,\;D_JQ_{P^4}\equiv0,\;Q_{P^3X^3}\not\equiv0\big\}.
\end{gather*}
In the above formula, the total derivative operator $D_J=D_{j^1}\cdots D_{j^k}$ in the variables $x^1=x$, $x^2=u$,
$x^3=p$ ranges over the multi-indices of length $\#J\geq0$.
Admissible cross-sections are given by
\begin{gather*}
\mathcal{K}^\infty_{2}=\mathcal{K}^\infty\cup\{q_{p^2x^2}=1,\; q_{p^3x^2}=q_{p^2ux^2}=q_{p^2x^3}=0,\;q_{p^2x^4}=
1\}\subset\mathcal{V}^\infty_2,
\\
\mathcal{K}^\infty_{3}=\mathcal{K}^\infty\cup\{q_{p^2x^2}=1,\;q_{p^3x^2}=q_{p^2ux^2}=q_{p^2x^3}=0,\;q_{p^3x^3}=
1\}\subset\mathcal{V}^\infty_3.
\end{gather*}
Similarly, in Case {\bf III)} the non-degeneracy condition $Q_{P^4}\not\equiv0$ implies that $Q_{P^6}$ and $Q_{P^5X}$
cannot both be equal to zero so that on
\begin{gather*}
\mathcal{V}^\infty_4=\big\{z\ii\in\J^\infty\,|\,D_JQ_{P^2X^2}\equiv0,\;Q_{P^4}\not\equiv0,\;Q_{P^6}\not\equiv0\big\},
\\
\mathcal{V}^\infty_5=\big\{z\ii\in\J^\infty\,|\,D_JQ_{P^2X^2}\equiv0,\;Q_{P^4}\not\equiv0,\;Q_{P^5X}\not\equiv0\big\},
\end{gather*}
moving frames can be constructed.
On $\mathcal{V}^\infty_4$ and $\mathcal{V}^\infty_5$ possible cross-sections are given by
\begin{gather*}
\mathcal{K}^\infty_{4}=\mathcal{K}^\infty\cup\{q_{p^4}=1,\;q_{p^5}=q_{p^4u}=q_{p^4x}=0,\;q_{p^6}=
1\}\subset\mathcal{V}^\infty_4,
\\
\mathcal{K}^\infty_{5}=\mathcal{K}^\infty\cup\{q_{p^4}=1,\;q_{p^5}=q_{p^4u}=q_{p^4x}=0,\;q_{p^5x}=
1\}\subset\mathcal{V}^\infty_5,
\end{gather*}
respectively.
\end{Example}

As illustrated by the above computations, the assumption that $\mathcal{V}^\infty$ is a~dense open subset of~$\J^\infty$ on which the pseudo-group acts regularly (and freely) is too restrictive.
The generic Case~{\bf I)} corresponding the set of regular submanifold jets $\mathcal{V}^\infty_1$ is the only branch
of the equivalence problem satisfying this assumption.
To encompass the other cases, we observe that, apart from a~f\/inite number of non-degeneracy conditions,
$\mathcal{V}^\infty_2, \ldots, \mathcal{V}^\infty_5$ are characterized by $\G$-invariant systems of dif\/ferential
equations which we formalize in the following def\/inition.
\begin{Definition}
\label{regular subbundle}
Let $\G$ be a~Lie pseudo-group acting on $\J^\infty$.
A \emph{$\G$-invariant} subbundle $S^\infty\subset\J^\infty$ is said to be \emph{regular} of order $\overline{n}\geq1$
if, for all f\/inite $n\geq\overline{n}\colon$
\begin{itemize}\itemsep=0pt
\item $S^n=\pi^\infty_n(S^\infty)\subset\J^n$ forms a~smooth embedded subbundle,
\item the projection
$\pi^{n+1}_n\colon S^{n+1}\to S^n$ is a~f\/ibration,
\item $S^n=\text{pr}^{(n-\overline{n})}S^{\overline{n}}$ is
obtained by prolongation, \item $S^n$ is $\G$-invariant.
\end{itemize}
\end{Definition}

In light of the above comments, we allow the set of regular submanifold jets to be the f\/inite union
\begin{gather*}
\mathcal{V}^\infty=\bigcup_{i=1}^k \mathcal{V}^\infty_i
\end{gather*}
of dense open subsets of $\G$-invariant regular subbundles of $\J^\infty$.
In other words, the subbundles~$\mathcal{V}^\infty_i$ are $\G$-invariant Zariski open subsets of $\J^\infty$.
Def\/inition~\ref{regular subbundle} implies that for each subbundle~$\mathcal{V}^\infty_i$ there exists~$n_i\geq1$
such that for $n\geq n_i$ the subbundle $\mathcal{V}^{n}_i=\pi^\infty_n(\mathcal{V}^\infty_i)\subset\J^n$ is
characterized by a~$\G$-invariant system of formally integrable dif\/ferential equations
\begin{gather}
\label{i-th determining system}
E^{(n)}_i(x,u^{(n)})=0
\end{gather}
plus, possibly, a~f\/inite number of $\G$-invariant non-degeneracy conditions.

\begin{Example}
The determining equations~\eqref{i-th determining system} and the non-degeneracy conditions naturally occur as one
tries to normalize the parameters of an equivalence pseudo-group in the moving frame algorithm.
In Example~\ref{regular moving frame example}, the determining system for~$\mathcal{V}_1^\infty$ is trivial while the
non-degeneracy conditions are given by
\begin{gather*}
Q_{P^4}\not\equiv0,
\qquad
Q_{P^2X^2}\not\equiv0.
\end{gather*}
For $\mathcal{V}_2^\infty$, the order of regularity is $n_2=4$, and for $n\geq n_2$ the determining system of
$\mathcal{V}^n_2$ is
\begin{gather*}
D_J\big(Q_{P^4}\big)\equiv0,
\qquad
0\leq\#J\leq n-n_2,
\end{gather*}
to which we add the non-degeneracy conditions $Q_{P^2X^2}\not\equiv0$, $Q_{P^2X^4}\not\equiv0$.
\end{Example}
\begin{Remark}
While the pseudo-group action will, in general, not be regular on the whole set of regular submanifold jets
$\mathcal{V}^\infty$, in order to construct a~moving frame on each invariant subbundle $\mathcal{V}^\infty_i$, we
require the restriction of pseudo-group action to $\mathcal{V}^\infty_i$ to be regular in the subset topology.
\end{Remark}

Assuming the action is free and regular on each invariant regular subbundle $\mathcal{V}^\infty_i$, we can construct
a~moving frame $\widehat{\rho}_i\colon\mathcal{V}^\infty_i\to\mathcal{B}\ii_i$
($\mathcal{B}\ii_i=\mathcal{B}\ii|_{\mathcal{V}^\infty_i}$) by choosing a~cross-section
\mbox{$\mathcal{K}^\infty_i\subset\mathcal{V}^\infty_i$}.
Each moving frame will have its own order of freeness $n_i^\star\geq n_i$ and on each subbundle $\mathcal{V}^\infty_i$
the recurrence relations~\eqref{horizontal recurrence relations} completely determine the algebra of dif\/ferential
invariants.

\subsection{Singular submanifold jets}
%\label{subsection totally singular jets}

We now would like to extend the moving frame method to submanifold jets $z\ii$ where a~pseudo-group does not act freely.
\begin{Example}
%\label{singular jet example}
Any pseudo-group satisfying
\begin{gather}
\label{large pseudo-groups}
r_n=\dim \G\n|_z>\dim \J^n|_z=q\binom{p+n}{p}
\qquad
\text{for all}
\quad
n\geq1
\end{gather}
cannot act freely.
Indeed, the inequality~\eqref{large pseudo-groups} implies that for all $n\geq1$ the isotropy group $\G_{z\n}\n$ is
non-trivial since the dimension of the pseudo-group jet $\G\n|_z$ is larger than the dimension of the jet space $\Jn$
on which it acts.
The contact pseudo-group~\eqref{contact} is an example of such pseudo-group.
\end{Example}
\begin{Example}
\label{partial moving frame example}
A second example is given by Case~{\bf IV)} of Example~\ref{regular moving frame example}.
When $Q_{P^4}\equiv Q_{P^2X^2}\equiv0$ it follows from the recurrence relations~\eqref{45 recurrence relations} that
the f\/ifth and sixth order (partially normalized) lifted invariants~\eqref{456 invariants} are identically equal to
zero.
This combined with the higher order recurrence relations then implies that all higher order
invariants~\eqref{intermediate remaining invariants} are also equal to zero.
Hence, in Case~{\bf IV)} there are no further lifted invariants available to normalize the Maurer--Cartan
forms~\eqref{unnormalized mc forms}.
Said dif\/ferently, on the subbundle
\begin{gather*}
%\label{point singular jets}
\mathscr{S}^\infty=\big\{z\ii\in\J^\infty\,\big|\,D_JQ_{P^4}\equiv D_JQ_{P^2X^2}\equiv0,
\;
\#J\geq0\big\}.
\end{gather*}
the pseudo-group parameters
\begin{gather*}
%\label{isotropy parameters}
\chi_x,
\qquad
\chi_u,
\qquad
\psi_u,
\qquad
\psi_{uu},
\qquad
\psi_{xu}
\end{gather*}
cannot be normalized and these pseudo-group jets parametrize the 5-dimensional isotropy group $\G\ii_{z\ii}$ of
a~submanifold jet $z\ii\in\mathscr{S}^\infty$.
\end{Example}
\begin{Definition}
%\label{definition totally singular}
A submanifold jet $z^{(\infty)}\in\J^\infty$ is said to be \emph{singular} if its isotropy group is non-trivial:
\begin{gather*}
\G\ii_{z\ii}=\big\{g\ii\in\G\ii\big|_z:\;g\ii\cdot z\ii=z\ii\big\}\neq\big\{\mathds{1}^{(\infty)}\big|_z\big\}.
\end{gather*}
The set of singular submanifold jets is denoted by $\mathscr{S}^\infty$.
\end{Definition}

As with regular submanifold jets, we allow the set of singular submanifold jets to be a~f\/inite union of open dense
subsets of $\G$-invariant regular subbundles
\begin{gather*}
\mathscr{S}^\infty=\bigcup_{i=1}^\ell \mathscr{S}^\infty_i
\end{gather*}
so that for $n\geq n_i\geq1$, the subbundle $\mathscr{S}^n_i=\pi^\infty_n(\mathscr{S}^\infty_i)$ is characterized by
a~formally integrable system of $\G$-invariant dif\/ferential equations (with possibly f\/initely many invariant
non-degeneracy conditions).
\begin{Remark}
Though the pseudo-group action on $\mathscr{S}^\infty_i$ is not free, we still require the restriction of the action to
$\mathscr{S}^\infty_i$ to be regular in the subset topology.
\end{Remark}

Let $z\ii\in\mathscr{S}^\infty_i$ be a~f\/ixed submanifold jet and $z\n=\pi^\infty_n\big(z\ii\big)\in\mathscr{S}^n_i$ its
$n^\text{th}$ order truncation.
Then $\G\ii_{z\ii}$ is the projective limit
\begin{gather*}
\G\ii_{z\ii}=\varprojlim \G\n_{z\n}.
\end{gather*}
Apart from the order 0 constraint $g\cdot z=z$, we observe that the isotropy requirement $g\n\cdot z\n=z\n$ gives
a~system of dif\/ferential constraints for the pseudo-group parameters $g\n\in\G\n_{z\n}$, and the prolongation of this
system is given by $g^{(n+1)}\cdot z^{(n+1)}=z^{(n+1)}$.
To see this, consider an analytic submanifold $(x,u(x))$ with $(x,\j_\infty u(x))=z\ii$.
Then, given
\begin{gather*}
X=\chi(x,u),\, U=\psi(x,u) \in \G_{z},
\end{gather*}
the isotropy condition $g^{(n)}\cdot z^{(n)}=z^{(n)}$ is obtained by dif\/ferentiating
\begin{gather*}
u(\chi(x,u(x)))=\psi(x,u(x))
\end{gather*}
with respect to the independent variables $x=(x^1,\ldots,x^p)$.
By Cartan--Kuranishi's prolongation theorem,~\cite{BCGGG-1991,S-2009}, we conclude that, generically, there exists
$n^\star_i\geq n_i$ such that the system of dif\/ferential equations $g^{(n^\star)}\cdot z^{(n^\star)}=z^{(n^\star)}$
is formally integrable (and eventually involutive for some $n\geq n^\star_i$).
In the following, $n^\star_i$ is assumed to be independent of the submanifold jet $z\ii\in\mathscr{S}^\infty_i$ and
$n^\star_i$ is called the \emph{order of partial freeness}.
For $n\geq n^\star_i$ the pseudo-group $\G$ is said to act \emph{partially freely} on $\mathscr{S}^n_i$.
\begin{Remark}
If $z\ii\in\mathcal{V}^\infty_i$ is a~regular submanifold jet, the order of partial freeness def\/ined above
corresponds to the usual order of freeness.
In this case the only solution to $g\n\cdot z\n=z\n$ is $\G\n_{z\n}=\{\mathds{1}\n|_z\}$ for $n\geq n^\star_i$.
\end{Remark}
\begin{Example}
\label{isotropy example}
To illustrate the above discussion we consider the Lie pseudo-group
\begin{gather}
\label{pj pseudo-group}
X=f(x),
\qquad
Y=e(x,y)=f^\prime(x) y+g(x),
\qquad
U=u+\frac{e_x(x,y)}{f^\prime(x)}
\end{gather}
acting on graphs of functions $(x,y,u(x,y))$.
This is one of the main examples used in~\cite{OP-2009-1, OP-2008} to illustrate the moving frame method.
The pseudo-group~\eqref{pj pseudo-group} acts locally freely on the set of regular submanifold jets
$\mathcal{V}^\infty_-=\J^\infty\cap\{u_{yy}<0\}$ and $\mathcal{V}^\infty_+=\J^\infty\cap\{u_{yy}>0\}$.
On the other hand, when restricted to the jet space of linear functions $u(x,y)=a(x) y+b(x)$ given by
\begin{gather}
\label{pj singular jets}
\mathscr{S}^\infty=\J^\infty\cap\{u_{x^iy^{j+2}}=0:\, i,j\geq0\},
\end{gather}
the pseudo-group action is not free.
To obtain the determining equations of the isotropy pseudo-group at a~submanifold jet
$z\ii=(x,y,u\ii)\in\mathscr{S}^\infty$ we dif\/ferentiate the equality
\begin{gather}\label{isotropy condition}
u(f(x),e(x,y))=u(x,y)+\frac{e_x(x,y)}{f_x(x)}
\end{gather}
with respect to $x$ and $y$.
Up to order 2, taking into account that $u_{yy}=0$ on $\mathscr{S}^\infty$, we obtain the constraints
\begin{gather*}
f=x,
\qquad
e=y,
\qquad
u_x f_x+u_y e_x=u_x+\frac{e_{xx}}{f_x}-\frac{e_x f_{xx}}{f_x^2},
\qquad
u_y f_x=u_y+\frac{f_{xx}}{f_x},
\\
u_{xx} f_x+u_x f_{xx}+u_{xy} f_x e_x+u_y e_{xx}=
u_{xx}+\frac{e_{xxx}}{f_x}-2\frac{e_{xx} f_{xx}}{f_x^2}-\frac{e_x f_{xxx}}{f_x^2}+2\frac{e_x f_{xx}^2}{f_x^3},
\\
u_{xy} f_x^2+u_y f_{xx}=u_{xy}+\frac{f_{xxx}}{f_x}-\frac{f_{xx}^2}{f_x^2},
\end{gather*}
on the jets of $f(x)$ and $e(x,y)$.
Solving the equations~\eqref{isotropy condition} (and their prolongations) for the pseudo-group parameters we conclude
that $f_x$ parametrizes the isotropy group $\G\ii_{z\ii}$.
\end{Example}

\subsubsection{Partial moving frames}
%\label{partial moving frame section}

Though it is not possible to construct a~moving frame on $\mathscr{S}^\infty_i$ as in Section~\ref{subsection regular
jets}, it is nevertheless possible to introduce the notion of a~\emph{partial moving frame}.
For $n^\star_i\leq n\leq\infty$, where $n^\star_i$ is the order of partial freeness, we introduce the $n^\text{th}$
order \emph{prolonged bundle}
\begin{gather*}
\mathcal{P}\n_i=\big\{\big(z\n,g\n\big):z\n\in\mathscr{S}^n_i
\
\text{and}
\
g\n\in\G\n_{z\n}\big\}.
\end{gather*}
A local dif\/feomorphism $h\in\G$ acts on the set
\begin{gather*}
\big\{\big(z\n,g\n\big)\in\mathcal{P}\n_i\,\big|\, \pi^n_0\big(z\n\big)\in\dom h\big\}
\end{gather*}
by
\begin{gather}
\label{action on prolonged bundle}
h\n\cdot\big(z\n,g\n\big)=\big(h\n\cdot z\n,K_{h\n}\big(g\n\big)\big),
\end{gather}
where $K_{h\n}\big(g\n\big)=h\n\cdot g\n\cdot\big(h^{-1}\big)\n$ is the conjugation action.

\begin{Definition}
For $n\geq n^\star_i$, an $n^\text{th}$ order (right) \emph{partial moving frame} on $\mathscr{S}^n_i$ is
a~$\G$-equivariant bundle map
\begin{gather*}
\widehat{\rho}^{ (n)}_i\colon \ \mathcal{P}\n_i\to\mathcal{B}\n_i,
\qquad
\text{where}
\qquad
\mathcal{B}\n_i=\mathcal{B}\n\big|_{\mathscr{S}^n_i}.
\end{gather*}
\end{Definition}

In terms of the action~\eqref{action on prolonged bundle}, right $\G$-equivariance means that
\begin{gather*}
R_h\widehat{\rho}^{ (n)}_i\big(z\n,g\n\big)=\widehat{\rho}^{ (n)}_i\big(h\n\cdot\big(z\n,g\n\big)\big).
\end{gather*}

\begin{figure}[!h]\centering
\includegraphics[scale=0.9]{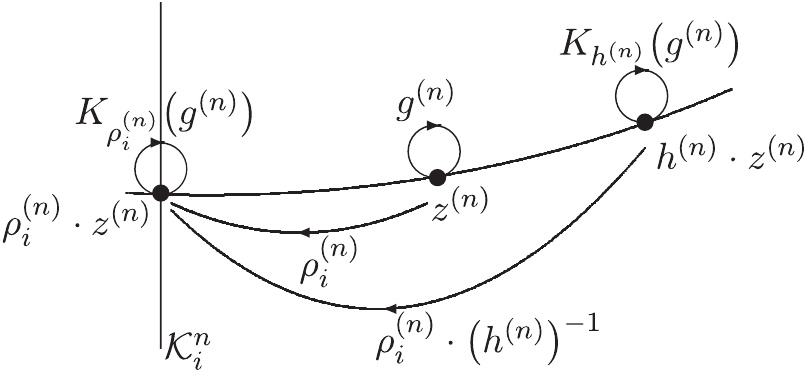}
%\begin{center}
%\begin{picture}
%(200,100) \qbezier(0,50)(95,45)(175,80) \put(10,5){\line(0,1){100}} \put(11,2){$\mathcal{K}^n_i$}
%\put(90,54.5){\circle*{5}} \put(88,42){$z\n$} \put(89,62){\circle{15}} \put(78,73){$g\n$}
%\put(90,69.5){\vector(1,0){1}} \put(150,70.5){\circle*{5}} \put(153,57){$h\n\cdot z\n$} \put(149,78){\circle{15}}
%\put(125,90){$K_{h\n}\big(g\n\big)$} \put(150,85.5){\vector(1,0){1}} \put(-36,38){$\rho\n_i\cdot z\n$} \put(10,50){\circle*{5}}
%\put(10,57){\circle{15}} \put(-15,70){$K_{\rho\n_i}\big(g\n\big)$} \put(10,64.5){\vector(1,0){1}} \qbezier(13,47)(50,30)(87,47)
%\put(50,38.5){\vector(-1,0){1}} \put(50,26){$\rho\n_i$} \qbezier(13,45)(80,-20)(150,63) \put(75,16.5){\vector(-1,0){1}}
%\put(72,4){$\rho\n_i\cdot\big(h\n\big)^{-1}$}
%\end{picture}
\caption{Partial moving frame $\widehat{\rho}^{ (n)}_i\big(z\n,g\n\big)=\big(z\n,\rho\n_i\big(z\n,g\n\big)\big)$, with $g\n\in\G\n_{z\n}$.}
\label{partial moving frame figure}
%\end{center}
\end{figure}

A partial moving frame is constructed by following the algorithm of Section~\ref{subsection regular jets}.
Namely, a~partial moving frame of order $n\geq n_i^\star$ is obtained by choosing a~minimal cross-section
$\mathcal{K}^n_i\subset\mathscr{S}^n_i$ and solving the normalization equations~\eqref{normalization equations}.
By the implicit function theorem, the normalization equations can be solved near the identity jet and the solution will
depend on the submanifold jet $z\n$ and the isotropy group parameters $g\n\in\G\n_{z\n}$.
See Fig.~\ref{partial moving frame figure} for a~suggestive illustration of a~partial moving frame.
\begin{Theorem}
Let $n\geq n_i^\star$ so that $\G$ acts partially freely on $\mathscr{S}^n_i\subset\Jn$ with its orbits forming
a~regular foliation.
Then an $n^{\text{th}}$ order partial moving frame exists in a~neighborhood of $z\n\in\mathscr{S}^n_i$.
\end{Theorem}

A partial moving frame $\widehat{\rho}_i\colon\mathcal{P}\ii_i\to\mathcal{B}\ii_i$ is constructed by choosing a~series
of compatible cross-sections $\mathcal{K}^n_i\subset\mathscr{S}_i^n$ just as in the regular case.
The def\/inition of the invariantization map~\eqref{invariantization map} and the recurrence formula~\eqref{universal
recurrence relations} still hold for partial moving frames with the understanding that these formulas are now def\/ined
on the prolonged bundle $\mathcal{P}\ii_i$.
We note that none of the normalized dif\/ferential invariants $\iota\big(z\ii\big)=\iota(x,u^{(\infty)})$ can depend on the
isotropy group parameters $g^{(\infty)}\in\G\ii_{z\ii}$.
On the other hand, the invariantization of the jet forms $dx^i$, $\theta^\alpha_J$ may involve the isotropy group
parameters.
\begin{Remark}
The concept of partial moving frame def\/ined above is similar to the recent notion of partial moving frame introduced
in~\cite{O-2011}.
Indeed, given a~cross-section $\mathcal{K}^\infty_i\subset\mathscr{S}^\infty_i$ def\/ining a~partial moving frame
$\widehat{\rho}_i$, its inverse image $\big(\bt\ii\big)^{-1}(\mathcal{K}^\infty_i)$ under the restricted target map
$\bt\ii\colon\mathcal{B}\ii_i\to\mathscr{S}^\infty_i$ happens to be equal to the image of the partial moving frame
$\widehat{\rho}_i\big(\mathcal{P}\ii_i\big)=\big(\bt\ii\big)^{-1}(\mathcal{K}^\infty_i)$.
\end{Remark}
\begin{Example}
Continuing Example~\ref{isotropy example}, and referring to~\cite{OP-2009-1, OP-2008} for all the formulas, the
prolonged bundle over the singular submanifold jets~\eqref{pj singular jets} is given by
\begin{gather*}
\mathcal{P}\ii=\big\{\big(x,y,u\ii,f_x\big):\, u_{x^i y^{j+2}}=0\big\}.
\end{gather*}
Up to second order, the expressions for the prolonged action are
\begin{gather*}
U_{X}=\frac{u_x}{f_x}+\frac{e_{xx}-e_x u_y}{f_x^2}-2\frac{f_{xx} e_x}{f_x^3},
\qquad
U_Y=\frac{u_y}{f_x}+\frac{f_{xx}}{f_x^2},
\\
U_{XX}=
\frac{u_{xx}}{f_x^2}\frac{e_{xxx}-e_{xx} u_y-2e_x u_{xy}-f_{xx} u_x}{f_x^3}+\frac{3e_x f_{xx} u_y-4e_{xx} f_{xx}-3e_x f_{xxx}}{f_x^4}+8\frac{e_x f_{xx}^2}{f_x^5},
\\
U_{XY}=\frac{u_{xy}}{f_x^2}+\frac{f_{xxx}-f_{xx} u_y}{f_x^3}-2\frac{f_{xx}^2}{f_x^4},
\qquad
U_{YY}=0.
\end{gather*}
We note that $U_{YY}\equiv0$ when $u_{yy}\equiv0$ and more generally $U_{X^i Y^{j+2}}\equiv0$, $i,j\geq0$.
A cross-section to the pseudo-group orbit is given by
\begin{gather*}
\mathcal{K}^\infty=\{x=y=u_{x^i}=u_{x^i y}=0,\,i\geq0\}\subset\mathscr{S}^\infty.
\end{gather*}
Solving the normalization equations $X=Y=U_{X^i}=U_{X^iY}=0$, we obtain the partial moving frame
\begin{gather*}
f=e=0,
\qquad
e_x=-u f_x,
\qquad
f_{xx}=-u_y f_x,
\qquad
e_{xx}=(u u_y-u_x)f_x,
\qquad
\ldots.
\end{gather*}
Since the action is transitive on~\eqref{pj singular jets} there are no dif\/ferential invariants.
The (partially) invariantized horizontal coframe is
\begin{gather*}
\omega^x=\iota(dx)=f_x dx,
\qquad
\omega^y=f_x(dy-u dx),
\end{gather*}
and their structure equations are
\begin{gather}
\label{3d mc structure equations}
d\omega^x\equiv\mu_X\wedge\omega^x,
\qquad
d\omega^y\equiv\mu_X\wedge\omega^y,
\qquad
d\mu_X\equiv0,
\end{gather}
where
\begin{gather*}
\mu_X=\widehat{\rho}^{\,*}(\mu_X)=\frac{df_x}{f_x}+u_y dx
\end{gather*}
is the only unnormalized Maurer--Cartan form ref\/lecting the fact that the pseudo-group jet $f_x$ parametrizes the
isotropy group of a~submanifold jet $z\ii\in\mathscr{S}^\infty$.
The structure equations~\eqref{3d mc structure equations} are isomorphic to the Maurer--Cartan structure equations of
a~3-dimensional Lie group indicating that the surface $(x,y,a(x) y+b(x))$ is invariant under a~3-dimensional group of
transformations lying inside the pseudo-group~\eqref{pj pseudo-group}.
\end{Example}

\section{Local equivalence}
\label{Section Equivalence Problems}

In this section we review the solution to the local equivalence problem of submanifolds using the equivariant moving
frame machinery.
We then explain how Cartan's approach based on the theory of $G$-structures is related to the equivariant method.

\subsection{Equivalence of submanifolds}

\begin{figure}[!h]\centering
\includegraphics[scale=0.9]{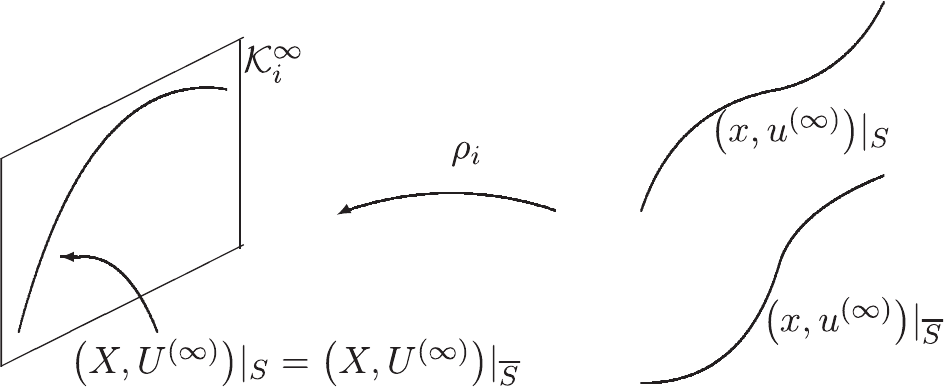}
%\begin{center}
%\begin{picture}
%(250,125) \put(0,15){\line(0,1){60}} \put(0,75){\line(2,1){70}} \put(0,15){\line(2,1){70}} \put(69,49){\line(0,1){60}}
%\put(70,100){$\mathcal{K}^\infty_i$} \qbezier(5,25)(25,101)(65,95)
%\put(20,13){$\big(X,U^{(\infty)}\big)|_S=\big(X,U^{(\infty)}\big)|_{\overline{S}}$} \qbezier(45,25)(35,50)(20,46)
%\put(22,46.5){\vector(-1,0){5}} \qbezier(100,60)(130,70)(160,60) \put(102,61.2){\vector(-2,-1){5}}
%\put(130,75){$\rho_i$} \qbezier(185,60)(195,90)(225,95) \qbezier(225,95)(245,100)(255,120)
%\put(205,80){$\big(x,u^{(\infty)}\big)|_S$} \qbezier(185,10)(215,10)(225,45) \qbezier(225,45)(230,60)(255,70)
%\put(220,25){$\big(x,u^{(\infty)}\big)|_{\overline{S}}$}
%\end{picture}
\caption{Signature of two equivalent submanifolds.}
\label{signature picture}
%\end{center}
\end{figure}

Given two $p$-dimensional submanifolds $S$ and $\overline{S}$ in $M$, the local equivalence problem consists of
determining whether there exists or not a~local dif\/feomorphism $g\in\G$ such that $g\cdot S=\overline{S}$ locally.
In accordance with Cartan's general philosophy, the solution to this problem is determined by the invariants of $\G$.
Let $\mathcal{K}^\infty_i$ be a~cross-section def\/ining a~(partial) moving frame $\widehat{\rho}_i$.
Since the normalized invariants $(X,U\ii)=\iota(x,u\ii)$ give a~local coordinate system on $\mathcal{K}^\infty_i$, two
submanifolds $S$ and $\overline{S}$ are locally equivalent if and only if both submanifolds lie in the domain of
def\/inition of the same (partial) moving frame and their ``projections" $\iota(\j_\infty S)$,
$\iota(\j_\infty\overline{S})$ onto~$\mathcal{K}^\infty_i$ overlap (see Fig.~\ref{signature picture}).
\begin{Definition}
Let $\D_1,\ldots,\D_p$ be the invariant total derivative operators~\eqref{invariant total derivative operators}.
A set of dif\/ferential invariants $\{I_\kappa\}$ is said to be a~\emph{generating set} for the algebra of
dif\/ferential invariants if any invariant can be locally expressed as a~function of the invariants $I_\kappa$ and
their invariant derivatives $\D_J I_\kappa=\D_{j_1}\D_{j_2}\cdots\D_{j_k}I_\kappa$.
\end{Definition}

First proved by Lie,~\cite[p.~760]{LS-1893}, for f\/inite-dimensional Lie group actions and then extended to
inf\/inite-dimensional Lie pseudo-groups by Tresse,~\cite{T-1894}, the \emph{fundamental basis theorem} (also known as
the \emph{Lie--Tresse theorem}) guarantees that, under appropriate technical hypotheses, the algebra of dif\/ferential
invariants is generated by a~f\/inite number of dif\/ferential invariants and~$p$ invariant total derivative operators.
For Lie group actions, recent proofs based on the equivariant moving frame method can be found in~\cite{FO-1999,H-2009}.
For inf\/inite-dimensional Lie pseudo-groups, a~proof also exists,~\cite{OP-2009-1}, but it is much more technical.
It requires the introduction of two modules associated with the prolonged pseudo-group action, and assumes the action
to be free and regular on a~dense open subset of $\J^\infty$.
This more abstruse part of the theory is discussed in Section~\ref{algebra of differential invariants - section}
where we explain how to modify the algebraic constructions introduced in~\cite{OP-2009-1} to cover Lie pseudo-groups
acting regularly and freely on invariant regular subbundles of $\J^\infty$.
Further modif\/ications will also allow us to deal with Lie pseudo-groups acting regularly and non-freely on invariant
regular subbundles of $\J^\infty$.
Other proofs based on Spencer cohomology,~\cite{K-1975}, Weyl algebras,~\cite{MMR-2003}, and homological
methods,~\cite{KL-2006}, are also available.
A global version of the Theorem for transitive algebraic pseudo-group actions was recently proved by Kruglikov and
Lychagin in~\cite{KL-2011}.

Given a~Lie pseudo-group $\G$ with regular submanifold jets $\mathcal{V}^\infty=\cup \mathcal{V}^\infty_i$ and
singular sub\-mani\-fold jets $\mathscr{S}^\infty=\cup \mathscr{S}^\infty_i$, the fundamental basis theorem applies to
each subbundle~$\mathcal{V}^\infty_i$ and~$\mathscr{S}^\infty_i$.
Let~$\{I_\kappa\}$ be a~generating set on $\mathcal{V}^\infty_i$ (or $\mathscr{S}^\infty_i$).
Then, the normalized invariants $(X,U\ii)=\iota(x,u\ii)$ can be expressed in terms of these invariants and their
invariant derivatives
\begin{gather}
\label{generating the normalized invariants}
X^i=\iota\big(x^i\big)=F^i( \ldots I_\kappa \ldots \D_K I_\kappa \ldots ),
\qquad
U^\alpha_J=\iota\big(u^\alpha_J\big)=F^\alpha_J( \ldots I_\kappa \ldots \D_K I_\kappa \ldots ).
\end{gather}
Hence, if $\j_\infty S$, $\j_\infty\overline{S}$ $\subset\mathcal{V}^\infty_i$ (or $\mathscr{S}^\infty_i$), it follows
that their invariantization $\iota(\j_\infty S)$, $\iota(\j_\infty\overline{S})$ is locally prescribed by the
generating invariants and their invariant derivatives.
\begin{Definition}
%\label{regular-singular S}
Let $\mathcal{V}^\infty=\cup \mathcal{V}^\infty_i$ and $\mathscr{S}^\infty=\cup \mathscr{S}^\infty_i$ denote the sets
of regular and singular submanifold jets.
A submanifold $S\subset M$ is said to be \emph{regular} (or \emph{singular}) if $\j_\infty
S\subset\mathcal{V}^\infty_i$ (or $\j_\infty S\subset\mathscr{S}^\infty_i$) for some $i$.
\end{Definition}
\begin{Remark}
Globally, it is possible that $\j_\infty S$ does not lie in a~unique subbundle $\mathcal{V}^\infty_i$ or
$\mathscr{S}^\infty_i$.
If so, the submanifold $S$ should be restricted to an open subset where the containment holds.
\end{Remark}
\begin{Definition}
\label{signature map definition}
Let $S$ be a~regular (or singular) submanifold with $\j_\infty S\subset\mathcal{V}_i^\infty$ (or $\j_\infty
S\subset\mathscr{S}^\infty_i$).
Let $\widehat{\rho}_i$ be a~moving frame (or a~partial moving frame) and $\{I_1,\ldots,I_{\ell}\}$ a~generating set for
the algebra of dif\/ferential invariants.
The $n^\text{th}$ \emph{order signature space} $\mathbb{K}\n$ is the Euclidean space of dimension
$\ell(1+p+p^2+\cdots+p^n)$ coordinatized by $w\n=(\ldots,w_{\kappa;J},\ldots)$, where
$(\kappa;J)=(\kappa,j^1,\ldots,j^r)$ with $1\leq\kappa\leq\ell$ and $(j^1,\ldots,j^r)$ ranging through all unordered
multi-index with $1\leq j^i\leq p$ and $0\leq r\leq n$.
The $n^\text{th}$ \emph{order signature map} associated with $\widehat{\rho}_i$ is the map $\mathbf{I}\n_S\colon
S\to\mathbb{K}\n$ whose components are
\begin{gather*}
w_{\kappa;J}=(\D_J I_\kappa)|_{\j_\infty S},
\qquad
\kappa=1,\ldots,\ell,
\qquad
\#J\leq n.
\end{gather*}
\end{Definition}
\begin{Remark}
In Def\/inition~\ref{signature map definition} the multi-index $J$ is not assumed to be symmetric as the invariant
total derivative operators $\D_j$ do not commute in general.
In applications we can reduce the dimension of the $n^\text{th}$ order signature space~$\mathbb{K}\n$ by ordering as
many multi-indices~$J$ as possible using the commutation relations~\eqref{commutation relations} and expressing the
commutator invariants~$Y^k_{ij}$ in terms of the generating invariants and their invariant derivatives:
\begin{gather*}
[\D_i,\D_j]=\sum_{k=1}^p Y^k_{ij}( \ldots I_\kappa \ldots \D_K I_\kappa \ldots ) \D_k,
\qquad
i,j=1,\ldots,p.
\end{gather*}
\end{Remark}
\begin{Definition}
The \emph{rank} of the signature map $\mathbf{I}\n_S$ at a~point $z\in S$ is the dimension of the space spanned by the
dif\/ferentials
\begin{gather*}
d\big[(\D_J I_\kappa)\big|_{\j_\infty S}\big]\big|_z,
\qquad
1\leq\kappa\leq\ell,
\qquad
0\leq\#J\leq n.
\end{gather*}
The signature map $\mathbf{I}\n_S$ is \emph{regular} if its rank is constant on $S$.
\end{Definition}
\begin{Definition}
Let $S$ be a~regular (or singular) submanifold with $\j_\infty S\subset\mathcal{V}_i^\infty$ (or $\j_\infty
S\subset\mathscr{S}_i^\infty$).
The restricted moving frame (or partial moving frame) $\widehat{\rho}_S^{}=\widehat{\rho}_i |_{\j_\infty S}^{}$ is
said to be \emph{fully regular} if for each $n\geq0$ the signature map $\mathbf{I}\n_S\colon S\to\mathbb{K}\n$ is
regular.
\end{Definition}
\begin{Definition}
%\label{definition signature manifold}
Let $\widehat{\rho}_S^{}$ be fully regular.
The image of the $n^\text{th}$ order signature map $\mathbf{I}\n_S$
\begin{gather*}
%\label{classifying manifold}
\mathfrak{S}\n\big(\widehat{\rho}_S^{}\big)=\big\{\mathbf{I}\n_S(z):z\in S\big\}\subset\mathbb{K}\n
\end{gather*}
is called the $n^\text{th}$ \emph{order signature manifold}.
\end{Definition}
\begin{Proposition}
Let $\widehat{\rho}_S^{}$ be fully regular, and let $\varrho_n$ denote the rank of the $n^\text{th}$ order signature
map $\mathbf{I}\n_S$.
Then
\begin{gather*}
0\leq\varrho_0<\varrho_1<\cdots<\varrho_s=\varrho_{s+1}=\cdots=r\leq p=\dim S.
\end{gather*}
The stabilizing rank $r$ is called the \emph{rank} of $\widehat{\rho}_S^{}$, and the smallest $s$ for which
$\varrho_s=\varrho_{s+1}=r$ is called the \emph{order} of $\widehat{\rho}_S^{}$.
\end{Proposition}

Depending if the two submanifolds $S, \overline{S}\subset M$ are regular or singular we are now ready to state the
solution to the local equivalence problem.
\begin{Theorem}
\label{Theorem coframe equivalence}
Let $\G$ be a~Lie pseudo-group acting regularly and freely on $\mathcal{V}^\infty_i\subset\J^\infty(M,p)$, and
$\widehat{\rho}_i\colon\mathcal{V}^\infty_i\to\mathcal{B}^{(\infty)}_i$ a~moving frame.
If $S, \overline{S}\subset M$ are two regular $p$-dimensional submanifolds with $\j_\infty S,
\j_\infty\overline{S}\subset\mathcal{V}^\infty_i$, then there exists a~local diffeomorphism $g\in\G$ mapping $S$
onto $\overline{S}$ locally if and only if $\widehat{\rho}_S^{}$ and $\widehat{\rho}_{\overline{S}}^{}$ are fully
regular, have the same order $s=\overline{s}$, and the $(s{+}1)^\text{st}$ order signature manifolds
$\mathfrak{S}^{(s+1)}(\widehat{\rho}_S^{})$, $\mathfrak{S}^{(s+1)}(\widehat{\rho}_{\overline{S}})^{}$ overlap.
Moreover, if $z\in S$ and $\overline{z}\in\overline{S}$ are mapped to the same point
\begin{gather*}
\mathbf{I}^{(s+1)}_S(z)=\mathbf{I}^{(s+1)}_{\overline{S}}(\overline{z})
\in
\mathfrak{S}^{(s+1)}\big(\widehat{\rho}_S^{}\big)\cap\mathfrak{S}^{(s+1)}\big(\widehat{\rho}_{\overline{S}}\big)
\end{gather*}
on the overlap of the two signature manifolds, then there exists a~unique local diffeomorphism $g\in\G$ mapping $z$
to $\overline{z}=g\cdot z$.
\end{Theorem}
\begin{proof}
Let $J_1,\ldots,J_r$ be a~set of invariants parametrizing the $s^\text{th}$ order signature manifold
$\mathfrak{S}^{(s)}(\widehat{\rho}_S)$.
By assumption, there exist \emph{signature functions} $F_{K;\kappa}(w_1,\ldots,w_r)$, such that
\begin{gather*}
\mathcal{D}_K I_\kappa=F_{K;\kappa}(J_1,\ldots,J_r),
\qquad
\kappa=1,\ldots,\ell,
\qquad
\#K\leq s.
\end{gather*}
Furthermore, the invariants $J_1,{\ldots},J_r$ also parametrize the $(s{+}1)^\text{th}$ signature mani\-fold
$\mathfrak{S}^{(s{+}1)}\!( \widehat{\rho}_S^{}\!).\!$
Hence, there exist signature functions $\widetilde{F}_{i;\upsilon}(w_1,\ldots,w_r)$ such that
\begin{gather*}
\D_i J_\upsilon=\widetilde{F}_{i;\upsilon}(J_1,\ldots,J_r),
\qquad
\upsilon=1,\ldots,r,
\qquad
i=1,\ldots,p.
\end{gather*}
By the chain rule
\begin{gather*}
\D_i(\D_K I_\kappa)=\sum_{\upsilon=
1}^r \pp{F_{K;\kappa}}{J_\upsilon}(J_1,\ldots,J_r)\cdot\widetilde{F}_{i;\upsilon}(J_1,\ldots,J_r),
\end{gather*}
and we conclude that once $\mathfrak{S}^{(s+1)}(\widehat{\rho}_S^{})$ is known, the signature manifold
$\mathfrak{S}^{(s+k)}(\widehat{\rho}_S^{})$ for $k\geq2$ is obtained by invariant dif\/ferentiation.

Using the assumption that the signature manifolds $\mathfrak{S}^{(s+1)}(\widehat{\rho}_S^{})$ and
$\mathfrak{S}^{(s+1)}(\widehat{\rho}_{\overline{S}}^{})$ overlap, it follows that the invariants
$\{I_1,\ldots,I_{\ell}\}$ generating the algebra of dif\/ferential invariants and their invariant derivatives are
locally equal when restricted to $S$ and $\overline{S}$.
By~\eqref{generating the normalized invariants}, the restriction of the the normalized invariants
$X^i|_S=X^i|_{\overline{S}}$ and $U^\alpha_K|_S=U^\alpha_K|_{\overline{S}}$ are also equal.
Hence, given $z\in S$ and $\overline{z}\in\overline{S}$ for which
$\mathbf{I}^{(s+1)}_S(z)=\mathbf{I}^{(s+1)}_{\overline{S}}(\overline{z})$ we conclude that $\iota(\j_\infty
S|_{z})=\iota\big(\j_\infty\overline{S}|_{\overline{z}}\big)$.
Now, let $h$, $\overline{h}\in\G$ such that
\begin{gather*}
h\ii|_{z}=\rho_i\big(\j_\infty S|_{z}\big)
\qquad
\text{and}
\qquad
\overline{h}\ii|_{\overline{z}}=\rho_i\big(\j_\infty\overline{S}|_{\overline{z}}\big).
\end{gather*}
Since
$\widetilde{\bs}\ii\bigl(\bigl(\overline{h}\ii|_{\overline{z}}\bigr)^{-1}\bigr)=
\iota\big(\j_\infty\overline{S}|_{\overline{z}}\big)=\iota\big(\j_\infty S|_{z}\big)=\widetilde{\bt}\ii\big(h\ii|_{z}\big)$, the map $g=\overline{h}^{-1}\comp h\in\G$ is (locally) well def\/ined and by construction $g\ii\cdot\j_\infty S|_{z}=\j_\infty\overline{S}|_{\overline{z}}$.
In the analytic category, this implies that $g\cdot S=\overline{S}$ locally, and in particular $g\cdot z=\overline{z}$.
\end{proof}
\begin{Theorem}
\label{equivalence theorem for singular submanifolds}
Let $\G$ be a~Lie pseudo-group acting regularly on $\mathscr{S}^\infty_i\subset\J^\infty(M,p)$ and
$\widehat{\rho}_i\colon\mathcal{P}\ii_i\!\to\mathcal{B}\ii_i$ a~partial moving frame.
If $S, \overline{S}$ $\subset M$ are two $p$-dimensional singular submanifolds with $\j_\infty S,
\j_\infty\overline{S}\subset\mathscr{S}^\infty_i$, then there exists a~local diffeomorphism $g\in\mathcal{G}$
sending $S$ onto $\overline{S}$ locally if and only if $\widehat{\rho}_S^{}$ and $\widehat{\rho}_{\overline{S}}^{}$ are
fully regular, have the same order $s=\overline{s}$, and the $(s{+}1)^\text{st}$ order signature manifolds
$\mathfrak{S}^{(s+1)}(\widehat{\rho}_S^{})$, $\mathfrak{S}^{(s+1)}(\widehat{\rho}_{\overline{S}}^{})$ overlap.
Moreover, if $z\in S$ and $\overline{z}\in\overline{S}$ are mapped to the same point
\begin{gather*}
\mathbf{I}^{(s+1)}_S(z)=\mathbf{I}^{(s+1)}_{\overline{S}}(\overline{z})
\in
\mathfrak{S}^{(s+1)}\big(\widehat{\rho}_S^{}\big)\cap\mathfrak{S}^{(s+1)}\big(\widehat{\rho}_{\overline{S}}^{}\big)
\end{gather*}
on the overlap of the two signature manifolds, then there is a~family of local equivalence maps sending $z$ to
$\overline{z}$.
Any two equivalence maps $g$, $\overline{g}$ are related by
\begin{gather*}
\overline{g}=\overline{h}\comp g\comp h,
\qquad
\text{with}
\qquad
\overline{h}\ii\big|_{\overline{z}}\in\G\ii_{\j_\infty\overline{S}|_{\overline{z}}}
\qquad
\text{and}
\qquad
h\ii\big|_{z}\in\G\ii_{\j_\infty S|_{z}}.
\end{gather*}
\end{Theorem}

\begin{proof}
The argument is identical to the proof of Theorem~\ref{Theorem coframe equivalence}.
The only dif\/ference is that the local dif\/feomorphism $g\in\G$ mapping $S$ onto $\overline{S}$ in the neighborhoods
of $z\in S$ and $\overline{z}\in\overline{S}$ is not uniquely def\/ined.
The dif\/feomorphism $g\in\G$ can be precomposed by any $h\in\G$ satisfying $h\ii|_z\in\G\ii_{\j_\infty S|_z}$ and
composed by $\overline{h}\in\G$ with $\overline{h}\ii|_{\overline{z}}\in\G\ii_{\j_\infty\overline{S}|_{\overline{z}}}$
to obtain a~new equivalence map $\overline{g}=\overline{h}\comp g\comp h$.
\end{proof}

\subsection{Equivalence of coframes}
\label{section equivalence of coframes}

In Cartan's framework, local equivalence problems are solved by recasting them as equivalence problems of coframes.
We brief\/ly explain how this is related to the equivalence problem of submanifolds.
Let $\mathcal{H}$ be a~Lie pseudo-group acting on a~$p$-dimensional manifold $X$, and let
\begin{gather}
\label{2 coframes}
\bo=\!\left\{\omega^i=\sum_{j=1}^p u^i_j(x) dx^j,\;i=1,\ldots,p\right\},
\qquad
\overline{\bo}=\!\left\{\overline{\omega}^i=\sum_{j=1}^p \overline{u}^i_j(\overline{x}) d\overline{x}^j,\;i=1,\ldots,p\right\}\!\!\!
\end{gather}
be two coframes on $X$ adapted to a~given equivalence problem.
We refer the reader to~\cite{G-1989,O-1995} for a~detailed account as to how to formulate equivalence problems using
coframes.
The local equivalence problem for the coframes~\eqref{2 coframes} consists of determining whether there exists or not
a~local dif\/feomorphism $\varphi\in\mathcal{H}$ such that
\begin{gather}
\label{equivalence criterion}
d\varphi^*(\overline{\omega}^i)=\sum_{j=1}^p h^i_j(x) \omega^j,
\qquad
i=1,\ldots,p.
\end{gather}
The matrix $(h^i_j(x))\in GL(p)$ in~\eqref{equivalence criterion} is contained in a~Lie group $H$ called the
\emph{structure group} of the equivalence problem.
As in the equivariant moving frame approach, the primary goal of Cartan's equivalence method is to reduce the
$H$-structure to an $\{e\}$-structure by normalizing the group coef\/f\/icients $h^i_j$ after which it is possible to
determine if the two coframes~\eqref{2 coframes} are locally equivalent.

Since coframes on $X$ are local sections of the coframe bundle $\mathcal{F}(X)$, the equivalence problem of coframes
can be interpreted as an equivalence problem of $p$-dimensional sections in $\mathcal{F}(X)$.
Let $M\subset\mathcal{F}(X)$ be the subbundle of all coframes~\eqref{2 coframes} adapted to a~given equivalence problem.
The action of $\mathcal{H}$ on $X$ naturally induces a~Lie pseudo-group action $\G\subset\mathcal{H}^{(1)}$ on the
subbundle~$M$ via the equivalence criterion~\eqref{equivalence criterion}.
Two coframes~$\bo$,~$\overline{\bo}$ are then locally equivalent if and only if their corresponding sections~$S,\overline{S}\subset M$ are equivalent up to a~local dif\/feomorphism~$g\in\mathcal{G}$.

There are three possible outcomes to the coframe equivalence problem, each having their counterparts in the equivariant
moving frame method.
Borrowing the terminology used in~\cite{O-1995}, after a~series of normalizations and prolongations, the three outcomes
to Cartan's algorithm are:
\begin{description}\itemsep=0pt
\item[a) Complete normalization:] The $H$-structure reduces to an $\{e\}$-structure.
This occurs when a~section in $M\subset\mathcal{F}(X)$ is regular and lies in the domain of def\/inition of a~moving
frame.
\item[b) Prolongation:] The coframe is prolonged to a~larger space on which the equivalence problem reduces to an
$\{e\}$-structure.
This situation occurs when a~section $S\subset M$ is singular and the isotropy group $\G\ii_{\j_\infty S|_{z}}$ is
f\/inite-dimensional.
\item[c) Involution:] The structure equations of the (possibly prolonged) coframe are in involution with nonzero
(reduced) Cartan characters.
This happens when a~section $S\subset M$ is singular and the isotropy group $\G\ii_{\j_\infty S|_{z}}$ is
inf\/inite-dimensional.
\end{description}

We now illustrate each outcome with an example.
At the same time we indicate some links between Cartan's approach and the equivariant moving frame method.

\subsubsection{Complete normalization}
%\label{complete normalization section}

Let $\mathcal{H}$ be a~Lie pseudo-group acting on a~manifold $X$ with local coordinates $x=(x^1,\ldots,x^p)$ and
$\G\subset\mathcal{H}^{(1)}$ the induced pseudo-group action on the subbundle of coframes $M\subset\mathcal{F}(X)$
adapted to an equivalence problem.
Given a~moving frame $\widehat{\rho}_i\colon\mathcal{V}^\infty_i\to\mathcal{B}\ii_i$ and a~regular section $S\subset M$
with $\j_\infty S\subset\mathcal{V}^\infty_i$, an invariant coframe on $X$ is obtained by restricting the invariant
horizontal 1-forms $\omega^i=\iota(dx^i)$ to $S$:
\begin{gather}
\label{normalized invariant horizontal coframe}
\bo=\bo|_S=\big\{\omega^j\equiv\widehat{\rho}^{\,*}_i\big(d_J X^j\big)=\iota\big(dx^j\big)\big\}\big|_S.
\end{gather}
The \emph{coframe derivatives} are then given by the invariant derivative operators
\begin{gather*}
\D_j=\D_j|_S=\pp{}{\omega^j},
\qquad
j=1,\ldots,p.
\end{gather*}
Provided the cross-section def\/ining the equivariant moving frame $\widehat{\rho}_i$ is compatible with the
normalizations leading to an $\{e\}$-structure in Cartan's coframe method, the coframe~\eqref{normalized invariant
horizontal coframe} will be equivalent to the one obtained via Cartan's procedure.

In Section~\ref{Section Moving Frame}, we saw that the structure equations for the coframe~\eqref{normalized invariant
horizontal coframe} can be deduced from the universal recurrence relation~\eqref{universal recurrence relations}.
In preparation for the next section, it is more convenient to recover these equations from the structure
equations~\eqref{Lie pseudo-group structure equations} of the equivalence pseudo-group.
Let $z=(x,u)$ be local coordinates on $M\simeq X\times U$.
Then, comparing the identities $dZ^a=\sigma^a+\mu^a$, $a=1,\ldots,m$, with the order 0 recurrence relations
\begin{gather*}
dX^i\equiv\omega^i+\mu^i,
\qquad
i=1,\ldots,p,
\qquad
dU^\alpha\equiv\sum_{j=1}^p U^\alpha_j \omega^j+\mu^{p+\alpha},
\qquad
\alpha=1,\ldots,q,
\end{gather*}
we see that
\begin{gather}
\label{pull-back sigma}
\sigma^i\equiv\omega^i
\qquad
\text{while}
\qquad
\sigma^{p+\alpha}\equiv\sum_{j=1}^p U^\alpha_j \omega^j.
\end{gather}
Hence, pulling-back the structure equations for $\sigma^1, \ldots, \sigma^p$ in~\eqref{Lie pseudo-group structure
equations} by the moving frame $\widehat{\rho}_i$ we obtain the structure equations for the coframe~\eqref{normalized
invariant horizontal coframe}.
\begin{Remark}
The pull-back by $\widehat{\rho}_i$ of the remaining structure equations for $\sigma^{p+1},\ldots,\sigma^m$ and the
Maurer--Cartan forms $\mu^a_B$ yields syzygies among the normalized invariants.
As these syzygies can be recovered from the recurrence relations~\eqref{horizontal recurrence relations}, we can ignore
these normalized structure equations.
\end{Remark}
\begin{Example}
\label{Cartan analogy example}
To make some analogies between the equivariant moving frame method and Cartan's coframe approach, the pseudo-group jets
need to be recursively normalized in the equivariant formalism,~\cite{O-2011}.
We illustrate this by considering the generic branch {\bf I)} of Example~\ref{regular moving frame example}.
Based on the cross-section~\eqref{K1} a~moving frame is obtained by solving the normalization equations
\begin{gather}
X=U=P=Q_{P^4}=Q_{P^5}=Q_{P^4U}=Q_{P^4X}=0,
\qquad
Q_{P^2X^2}=1,
\nonumber
\\
Q_{U^j X^k}=Q_{PU^j X^k}=Q_{P^2U^j}=Q_{P^2U^j X}=Q_{P^3U^j}=Q_{P^3U^j X}=0,
\qquad
j,k\geq0.
\label{K1 normalizations}
\end{gather}
In the recursive moving frame construction, the idea is to normalize the submanifold jets (or solve the normalization
equations) in stages.
For example, at order 0 we solve the normalization equations
\begin{gather*}
X=U=P=Q=0
\end{gather*}
to obtain
\begin{gather*}
\chi=0,
\qquad
\psi=0,
\qquad
\psi_x=-p \psi_u,
\qquad
\psi_{xx}=-\big(q \psi_u+p^2 \psi_{uu}+2p \psi_{xu}\big).
\end{gather*}
Then the (partially) normalized coframe~\eqref{normalized invariant horizontal coframe} yields the $G$-structure
\begin{gather}
\label{cartan coframe}
\begin{pmatrix}
\omega^x
\\
\omega^u
\\
\omega^p
\end{pmatrix}
=
\begin{pmatrix}
D_x\chi&\chi_u&0
\\
0&\psi_u&0
\\
0&\frac{D_x\psi_u}{D_x\chi}&\frac{\psi_u}{D_x\chi}
\end{pmatrix}\!\!
\begin{pmatrix}
dx
\\
du-p dx
\\
dp-q dx
\end{pmatrix}
\end{gather}
which is the starting point of Cartan's method.
The next step in Cartan's algorithm is to compute the structure equations for the 1-forms~\eqref{cartan coframe}.
For the problem at hand, the resulting structure equations are then prolonged,~\cite[p.~403]{O-1995}.
To obtain the same structure equations, using the equivariant formalism, we must normalize
\begin{gather*}
Q_{U^j X^k}=Q_{P U^j X^k}=0,
\qquad
j,k\geq0.
\end{gather*}
The recurrence relations~\eqref{jet coordinate recurrence relations} are then used to obtain expressions for the
(partially) normalized Maurer--Cartan forms $\mu$, $\nu$, $\nu_X$, $\mu_{U^j X^{k+2}}$, $\nu_{U^j X^{k+2}}$, $j,k\geq0$.
For example, the low order norma\-li\-zed Maurer--Cartan forms are
\begin{gather*}
\mu\equiv-\omega^x,
\qquad
\nu\equiv-\omega^u,
\qquad
\nu_X\equiv-\omega^p,
\\
\nu_{XX}\equiv\nu_{UXX}\equiv\nu_{XXX}\equiv0,
\qquad
\mu_{XX}\equiv2\nu_{UX}+Q_{PP} \omega^p.
\end{gather*}
Substituting these expressions into the structure equations~\eqref{mc structure equations}, and recalling
equality~\eqref{pull-back sigma}, we obtain
\begin{gather}
d\omega^x\equiv\mu_X\wedge\omega^x+\mu_U\wedge\omega^u,
\nonumber
\\
d\omega^u\equiv\nu_U\wedge\omega^u+\omega^x\wedge\omega^p,
\nonumber
\\
d\omega^p\equiv\nu_{UX}\wedge\omega^u+(\nu_U-\mu_X)\wedge\omega^p,
\nonumber
\\
d\mu_X\equiv-2\nu_{UX}\wedge\omega^x-\mu_{UX}\wedge\omega^u-\mu_U\wedge\omega^p+Q_{PP} \omega^x\wedge\omega^p,
\nonumber
\\
d\mu_U\equiv-\mu_{UU}\wedge\omega^u-\mu_{UX}\wedge\omega^x+\mu_U\wedge(\nu_U-\mu_X),
\nonumber
\\
d\nu_U\equiv-\nu_{UX}\wedge\omega^x+\mu_U\wedge\omega^p-\nu_{UU}\wedge\omega^u,
\nonumber
\\
d\nu_{UX}\equiv-\nu_{UUX}\wedge\mu^u+(\mu_{UX}-\nu_{UU})\wedge\omega^p+\nu_{UX}\wedge\mu_X.
\label{Cartan prolongation 1}
\end{gather}
These equations are equivalent to the ones obtained using Cartan's method,~\cite[pp.~403--404]{O-1995}.
The correspondence is given by
\begin{gather*}
\theta^1\leftrightarrow\omega^u,
\qquad
\theta^2\leftrightarrow\omega^p,
\qquad
\theta^3\leftrightarrow\omega^x,
\\
\pi^1\leftrightarrow\nu_U,
\qquad
\pi^2\leftrightarrow\nu_{UX},
\qquad
\pi^4\leftrightarrow\mu_U,
\qquad
\pi^6\leftrightarrow\mu_X,
\\
\rho^1\leftrightarrow-\nu_{UU},
\qquad
\rho^2\leftrightarrow-\nu_{UUX},
\qquad
\rho^3\leftrightarrow-\mu_{UX},
\qquad
\rho^4\leftrightarrow-\mu_{UU},
\qquad
T\leftrightarrow-Q_{PP}.
\end{gather*}
The next step in Cartan's algortithm is to normalize $Q_{PP}=0$.
In the equivariant setting, we now set $Q_{PPU^j}=0$, $j\geq0$, which leads to the normalization of the Maurer--Cartan
forms $\mu_{U^{j+1}X}$.
For example, the f\/irst normalized Maurer--Cartan form is
\begin{gather}
\label{muUX normalization}
\mu_{UX}=\frac{1}{4}\big(Q_{PPP} \omega^p+Q_{PPX} \omega^x+2\nu_{UU}\big),
\end{gather}
which we then substitute in~\eqref{Cartan prolongation 1}.
The invariants $Q_{PPP}$ and $Q_{PPX}$ in~\eqref{muUX normalization} are \emph{essential torsion coefficients} of
the resulting structure equations and the next iteration of Cartan's algorithm is to normalize these invariants to zero.
In the equivariant moving frame framework we now set $Q_{P^3U^j}=Q_{P^2U^jX}=Q_{P^3U^jX}=0$, $j\geq0$, and normalize
the Maurer--Cartan forms $\mu_{U^{j+2}}$, $\nu_{U^{j+2}X}$, $\nu_{U^{j+3}}$, $j\geq0$.
At this point we have recovered the universal normalizations~\eqref{intermediate normalizations} and all the
Maurer--Cartan forms are normalized except for~\eqref{unnormalized mc forms}.
From the recurrence relations, the low order normalized Maurer--Cartan forms are
\begin{gather}
\mu\equiv-\omega^x,
\qquad
\nu\equiv-\omega^u,
\qquad
\nu_X\equiv-\omega^p,
\qquad
\nu_{XX}\equiv0,
\nonumber
\\
\mu_{XX}\equiv2\nu_{UX},
\qquad
\mu_{UX}\equiv\frac{1}{2}\nu_{UU},
\qquad
\mu_{UU}\equiv\frac{1}{6}Q_{P^4} \omega^p,
\qquad
\nu_{UXX}\equiv0,
\nonumber
\\
\nu_{UXX}\equiv\frac{1}{6}Q_{P^2X^2} \omega^x,
\qquad
\nu_{UUU}\equiv\frac{1}{3}\big(Q_{P^4X} \omega^p+Q_{P^3X^2} \omega^x\big).
\label{universal mc normalizations}
\end{gather}
Substituting~\eqref{universal mc normalizations} into the structure equations~\eqref{mc structure equations} we obtain
the structure equations for the eight-dimensional invariant coframe
$\{\omega^x,\omega^u,\omega^p,\mu_X,\mu_U,\nu_U,\nu_{UX},\nu_{UU}\}$:
\begin{gather}
d\omega^x\equiv\mu_X\wedge\omega^x+\mu_U\wedge\omega^u,\nonumber \\
d\omega^u\equiv\nu_U\wedge\omega^u+\omega^x\wedge\omega^p,\nonumber \\
d\omega^p\equiv\nu_{UX}\wedge\omega^u+(\nu_U-\mu_X)\wedge\omega^p,\nonumber \\
d\mu_X\equiv-\frac{1}{2}\nu_{UU}\wedge\omega^y-2\nu_{UX}\wedge\omega^x-\mu_U\wedge\omega^p,\nonumber
\\
d\mu_U\equiv-\frac{1}{2}\nu_{UU}\wedge\omega^x+\frac{1}{6}Q_{P^4} \omega^u\wedge\omega^p+(\mu_X-\nu_U)\wedge\mu_U,\nonumber \\
d\nu_U\equiv-\nu_{UU}\wedge\omega^u-\nu_{UX}\wedge\omega^x+\mu_U\wedge\omega^p,\nonumber \\
d\nu_{UX}\equiv-\frac{1}{2}\nu_{UU}\wedge\omega^p+\nu_{UX}\wedge\mu_X+\frac{1}{6}Q_{P^2X^2} \omega^u\wedge\omega^x,\nonumber \\
d\nu_{UU}\equiv2\nu_{UX}\wedge\mu_U+\nu_{UU}\wedge\nu_U+\frac{1}{3}Q_{P^4X} \omega^u\wedge\omega^p+\frac{1}{3}Q_{P^3X^2} \omega^u\wedge\omega^x.
\label{universal point structure equations}
\end{gather}
These equations are equivalent to Cartan's structure equations~\cite[equation~(12.72)]{O-1995}.
The correspondence is given by
\begin{gather*}
\omega^x\leftrightarrow\theta^3,
\qquad
\omega^u\leftrightarrow\theta^1,
\qquad
\omega^p\leftrightarrow\theta^2,
\\
\mu_X\leftrightarrow\pi^6,
\qquad
\mu_U\leftrightarrow\pi^4,
\qquad
\nu_U\leftrightarrow\pi^1,
\qquad
\nu_{UX}\leftrightarrow\pi^2,
\qquad
\nu_{UU}\leftrightarrow-\rho,
\\
Q_{P^4}\leftrightarrow6 K_1,
\qquad
Q_{P^2X^2}\leftrightarrow6 K_2,
\qquad
Q_{P^4X}\leftrightarrow3 K_3,
\qquad
Q_{P^3X^2}\leftrightarrow3 K_4.
\end{gather*}
The structure equations~\eqref{universal point structure equations} and identity $d\comp d=0$ imply that
$Q_{P^3X^2}=\D_P Q_{P^2X^2}$ and $Q_{P^4X}=\D_X Q_{P^4}$.
We note that these dif\/ferential relations can also be deduced from the recurrence relations~\eqref{45 recurrence relations}.

At this stage, the equivalence problem splits into the branches identif\/ied in Example~\ref{regular moving frame example}.
For the generic branch, corresponding to Case~{\bf I)}, we can normalize the remaining pseudo-group parameters by
making the normalizations~\eqref{K1 normalizations}.
To obtain the structure equations of the invariant coframe $\{\omega^x,\omega^u,\omega^p\}$ we use the recurrence
relations to f\/ind the expressions for the low order normalized Maurer--Cartan forms
\begin{gather}
\mu\equiv-\omega^x,
\qquad
\nu\equiv-\omega^u,
\qquad
\nu_X\equiv-\omega^p,\qquad \mu_U\equiv-\frac{1}{5}\big(Q_{P^5X} \omega^x+Q_{P^5U} \omega^u+Q_{P^6} \omega^p\big),
\nonumber
\\
\frac{2}{3}\mu_X\equiv\nu_U\equiv\frac{1}{4}\big(Q_{P^2X^3} \omega^x+Q_{P^2UX^2} \omega^u+Q_{P^3X^2} \omega^p\big),
\qquad
\nu_{XX}\equiv0,
\nonumber
\\
\nu_{UX}\equiv-\big(Q_{P^4X^2} \omega^x+Q_{P^4UX} \omega^u+Q_{P^5X} \omega^p\big).
\label{mc normalizations - regular case}
\end{gather}
Substituting~\eqref{mc normalizations - regular case} in the f\/irst three structure equations of~\eqref{universal
point structure equations} we obtain
\begin{gather*}
d\omega^x\equiv\left(\frac{3}{8}Q_{P^2UX^2}+\frac{1}{5}Q_{P^5X}\right)
\omega^u\wedge\omega^x+\frac{3}{8}Q_{P^3X^2} \omega^p\wedge\omega^x+\frac{1}{5}Q_{P^6} \omega^u\wedge\omega^p,
\\
d\omega^u\equiv\omega^x\wedge\omega^p+\frac{1}{4}Q_{P^2X^3} \omega^x\wedge\omega^u+\frac{1}{4}Q_{P^3X^2} \omega^p\wedge\omega^u,
\\
d\omega^p\equiv\left(Q_{P^5X}-\frac{1}{8}Q_{P^2UX^2}\right)
\omega^u\wedge\omega^p+Q_{P^4X^2} \omega^u\wedge\omega^x+\frac{1}{8}Q_{P^2X^3} \omega^p\wedge\omega^x.
\end{gather*}
Extracting the invariants from the structure coef\/f\/icients, it follows from Cartan's moving frame theory that the
normalized invariants
\begin{gather}
\label{case 1 generating invariants}
Q_{P^3X^2},
\qquad
Q_{P^2UX^2},
\qquad
Q_{P^2X^3},
\qquad
Q_{P^6},
\qquad
Q_{P^5X},
\qquad
Q_{P^4X^2}
\end{gather}
form a~generating set for the algebra of dif\/ferential invariants.
We note that the generating set~\eqref{case 1 generating invariants} does contain all the normalized invariants of
order $\leq6$.
Based on the normalizations~\eqref{K1 normalizations}, the complete list of normalized invariants of order $\leq6$ is
\begin{gather*}
\begin{split}
& Q_{P^3X^2},
\qquad
Q_{P^2UX^2},
\qquad
Q_{P^2X^3},
\qquad
Q_{P^6},
\\
& Q_{P^5U},
\qquad
Q_{P^5X},
\qquad
Q_{P^4U^2},
\qquad
Q_{P^4UX},
\qquad
Q_{P^4X^2}.
\end{split}
\end{gather*}
Using, as always, the recurrence relations we can check that the extra invariants $Q_{P^5U}$, $Q_{P^4U^2}$, $Q_{P^4UX}$
can be expressed in terms of~\eqref{case 1 generating invariants} and their invariant derivatives:
\begin{gather*}
Q_{P^5U}\!=\!
\frac{5}{6}\bigg[\D_U Q_{P^3X^2}\!-\!\D_PQ_{P^2UX^2}\!-\!\frac{3}{10}Q_{P^3X^2}Q_{P^2UX^2}\!+\!Q_{P^5X}Q_{P^3X^2}\!+\!\frac{1}{5}Q_{P^6}Q_{P^2X^3}\bigg],
\\
Q_{P^4U X}\!=\!
\frac{1}{6}\bigg[\D_UQ_{P^2X^3}\!-\!\D_XQ_{P^2UX^2}\!+\!\frac{1}{8}Q_{P^2UX^2}Q_{P^2X^3}\!+\!Q_{P^4X^2}Q_{P^3X^2}\!+\!\frac{1}{5}Q_{P^5X}Q_{P^2X^3}\bigg],
\\
Q_{P^4U^2}\!=\!
\frac{4}{5}\bigg[\D_UQ_{P^5X}\!-\!\D_XQ_{P^5U}\!+\!\frac{1}{4}Q_{P^5X}Q_{P^2UX^2}\!+\!\frac{1}{5}Q_{P^5X}^2\!+\!Q_{P^6}Q_{P^4X^2}\!-\!\frac{1}{8}Q_{P^5U}Q_{P^2X^3}\bigg].
\end{gather*}

The structure equations for the invariant coframe $\{\omega^x,\omega^u,\omega^p\}$ and the determination of
a~generating set of invariants for the other regular cases of Example~\ref{regular moving frame example} are obtained
in a~similar fashion.
\end{Example}

\subsubsection{Prolongation}
%\label{coframe prolongation section}

We now assume that a~section $S\subset M\subset\mathcal{F}(X)$ is singular with $\j_\infty
S\subset\mathscr{S}^\infty_i$.
Also, we consider the case when the isotropy group $\G\ii_{\j_\infty S|_{z}}$ is f\/inite-dimensional and locally
parametrized by $g=(g_1,\ldots,g_r)$.

Let $\pmb{\mu}=\{\mu^1,\ldots,\mu^r\}$ be the Maurer--Cartan forms associated with the isotropy group parameters
$g=(g_1,\ldots,g_r)$.
Then the invariant horizontal 1-forms~\eqref{normalized invariant horizontal coframe} together with the Maurer--Cartan
forms $\pmb{\mu}=\widehat{\rho}^{\,*}_i(\pmb{\mu})|_S$ constitute an invariant coframe on the prolonged bundle
$\mathcal{P}\ii_i|_S$.
As in the previous section, the structure equations of the prolonged coframe $\{\bo,\pmb{\mu}\}$ are obtained by
pulling-back the structure equations of the equivalence pseudo-group by the partial moving frame $\widehat{\rho}_i$.
\begin{Example}
In Example~\ref{partial moving frame example} we saw that if the fourth order (relative) invariants $Q_{P^4}\equiv
Q_{P^2X^2}\equiv0$ are identically zero then all higher order invariants~\eqref{intermediate remaining invariants} are
also equal to zero.
Hence, in Case~${\bf IV)}$ the only invariants are constant functions, and two singular second order ordinary
dif\/ferential equations are equivalent under point transformations,~\cite{C-1955,G-1989,O-1995,T-1896}.

Setting $Q_{P^4}\equiv Q_{P^2X^2}\equiv Q_{P^4X}\equiv Q_{P^3X^2}\equiv0$ in the structure equations~\eqref{universal
point structure equations} we obtain the Maurer--Cartan structure equations for the eighth-dimensional symmetry group
of point transformations (isomorphic to SL(3)) of a~singular second order ordinary dif\/ferential equation:
\begin{gather}
d\omega^x\equiv\mu_X\wedge\omega^x+\mu_U\wedge\omega^u,\nonumber \\
d\omega^u\equiv\nu_U\wedge\omega^u+\omega^x\wedge\omega^p,\nonumber\\
d\omega^p\equiv\nu_{UX}\wedge\omega^u+(\nu_U-\mu_X)\wedge\omega^p,\nonumber\\
d\mu_U\equiv-\frac{1}{2} \nu_{UU}\wedge\omega^x+(\mu_X-\nu_U)\wedge\mu_U,\nonumber\\
d\nu_U\equiv\mu_U\wedge\omega^p-\nu_{UX}\wedge\omega^x-\nu_{UU}\wedge\omega^u,\nonumber \\
d\mu_X\equiv-2 \nu_{UX}\wedge\omega^x-\frac{1}{2} \nu_{UU}\wedge\omega^u-\mu_U\wedge\omega^p,
\nonumber\\
d\nu_{UU}\equiv2 \nu_{UX}\wedge\mu_U+\nu_{UU}\wedge\nu_U,\nonumber \\
d\nu_{UX}\equiv\nu_{UX}\wedge\mu_X-\frac{1}{2} \nu_{UU}\wedge\omega^p.
\label{point singular structure equations}
\end{gather}

Restricting~\eqref{point singular structure equations} to a~f\/ixed point $(x,u,p)$ we obtain the structure equations
\begin{gather}
d\mu_X\equiv0,
\qquad
d\mu_U\equiv(\mu_X-\nu_U)\wedge\mu_U,
\qquad
d\nu_U\equiv0,
\nonumber
\\
d\nu_{UU}\equiv2 \nu_{UX}\wedge\mu_U+\nu_{UU}\wedge\nu_U,
\qquad
d\nu_{UX}\equiv\nu_{UX}\wedge\mu_X,
\label{isotropy mc equations}
\end{gather}
for the f\/ive-dimensional isotropy group $\G\ii_{z\ii}$.
Identifying the jets of a~local vector f\/ield~\eqref{vector field} with the coef\/f\/icients of its Taylor
expansion,~\cite{OPV-2009, V-2008}:
\begin{gather*}
\j_\infty\vv|_{z_0}\simeq\sum_{a=1}^b\sum_{\#B\geq0} \zeta^a_B(z_0)\frac{(z-z_0)^B}{B!}\pp{}{z^a},
\end{gather*}
the f\/iber $\J^\infty TM|_{z_0}$ inherits a~Lie algebra structure.
The monomial vector f\/ields
\begin{gather*}
\vv^B_a=\frac{(z-z_0)^B}{B!}\pp{}{z^a},
\qquad
\#B\geq0,
\qquad
a=1,\ldots,m,
\end{gather*}
provide a~basis of the vector space $\J^\infty TM|_{z_0}$ and they can be identif\/ied as dual vectors to the
Maurer--Cartan forms $\mu^a_B|_{z_0}$.
Choosing, for simplicity, $z_0=(x_0,u_0,p_0,q_0)=0$, the isotropy Lie algebra at the origin is spanned by the
inf\/initesimal generators\footnote{The inf\/initesimal generators are prolonged to $p=u_x$ and $q=u_{xx}$ using
formula~\eqref{prolonged vector field}.}
\begin{gather*}
\vv^x_x=x\pp{}{x},
\qquad
\vv^x_u=u\pp{}{x},
\qquad
\vv^u_u=u\pp{}{u},
\\
\vv^u_{ux}=x^2\pp{}{x}+xu\pp{}{u},
\qquad
\vv^u_{uu}=\frac{xu}{2}\pp{}{x}+\frac{u^2}{2}\pp{}{u},
\end{gather*}
and satisfy commutation relations dual to~\eqref{isotropy mc equations}.
That is,
\begin{alignat*}{4}
& [\vv^x_x, \vv^u_{ux}]=\vv^u_{ux},
\qquad &&
[\vv^x_u, \vv^u_{ux}]=2\vv_{uu}^u,
\qquad &&
[\vv^u_u, \vv^u_{uu}]=\vv^u_{uu}, &
\\
& [\vv^x_u, \vv^x_x]=\vv^x_u,
\qquad &&
[\vv^u_u, \vv^x_u]=\vv^x_u. &
\end{alignat*}
The corresponding 5-dimensional isotropy group is generated by the local transformations
\begin{gather*}
(X,U)=(\lambda_1x,u),
\qquad
(X,U)=(x+\epsilon_1u,u),
\qquad
(X,U)=(x,\lambda_2u),
\\
(X,U)=\left(\frac{x}{1-\epsilon_2x},\frac{u}{1-\epsilon_2x}\right),
\qquad
(X,U)=\left(\frac{x}{1-\epsilon_3u},\frac{u}{1-\epsilon_3u}\right),
\end{gather*}
which can be prolonged to $p=u_x$, $q=u_{xx}$.
\end{Example}

\subsubsection{Involution}
%\label{involution section}

We now consider the case when the isotropy group is inf\/inite-dimensional.
As in the preceding section, let $g\n=(g_1,\ldots,g_{r_n})$ denote the pseudo-group jets parametrizing $\G\n_{\j_n
S|_{z}}$ and let $\pmb{\mu}\n$ denote the corresponding Maurer--Cartan forms.
Since the isotropy group is inf\/inite-dimensional, the collection of 1-forms $\{\bo,\pmb{\mu}^{(n)}\}$ does not form
a~coframe on the prolonged bundle $\mathcal{P}_i\n|_S$.
Nevertheless, the structure equations for the 1-forms $\{\bo,\pmb{\mu}^{(n)}\}$ are still obtained by pulling-back the
structure equations of the equivalence pseudo-group by the partial moving frame.
\begin{Example}
\label{involutive example}
In this example we consider the local equivalence of second order ordinary dif\/ferential equations~\eqref{ODE} under
the pseudo-group of contact transformations~\eqref{contact}.
Ta\-king into account the linear relations~\eqref{contact lifted equations} among the Maurer--Cartan forms, the
recurrence relations~\eqref{jet coordinate recurrence relations} for the lifted invariants are
\begin{gather}
dX\equiv\omega^x+\mu^x,
\qquad
dU\equiv\omega^u+\mu^u,
\qquad
dX^p\equiv\omega^p+\mu^p,
\nonumber
\\
dQ_J\equiv Q_{J,X} \omega^x+Q_{J,U} \omega^u+Q_{J,P} \omega^p+\mu^p_{J,X}
\nonumber
\\
\phantom{dQ_J\equiv}
{}+\mathbb{D}_J[P\mu^p_U+Q\mu^p_P-Q(\mu^x_X+P\mu^x_U+Q\mu^x_P)
\nonumber
\\
\phantom{dQ_J\equiv}
{}-\mu^xQ_X-\mu^uQ_U-\mu^pQ_P]+\mu^xQ_{J,X}+\mu^uQ_{J,U}+\mu^PQ_{J,P}.
\label{ODE recurrence relations contact transformations}
\end{gather}
Since the Maurer--Cartan forms $\mu^p_{J,X}$ are linearly independent, it follows from~\eqref{ODE recurrence relations
contact transformations} that the equivalence pseudo-group is transitive on $\J^\infty$.
Also, the inequality $\dim \G\n>\dim \J^n$ for all $n$ implies that the prolonged action is nowhere free and
all submanifold jets are singular so that $\mathscr{S}^\infty=\J^\infty$.
Choosing the cross-section
\begin{gather*}
\mathcal{K}^\infty=\{x=u=p=q_J=0,\;\#J\geq0\},
\end{gather*}
the normalized Maurer--Cartan forms are
\begin{gather*}
%\label{ODE contact MC normalization}
\mu^x\equiv-\omega^x,
\qquad
\mu^u\equiv-\omega^u,
\qquad
\mu^p\equiv-\omega^p,
\qquad
\mu^p_{J,X}\equiv0.
\end{gather*}
The structure equations for the invariant 1-forms $\bo=\{\omega^x,\omega^u,\omega^p\}$ are obtained by substituting
$\sigma^i\equiv-\mu^i\equiv\omega^i$, $1\leq i\leq3$, into the structure equations~\eqref{horizontal contact structure
equations}.
Taking into account that~\eqref{contact lifted equations} implies $\mu^p_U=\mu^u_{UX}$ when $P=0$, we obtain
\begin{gather}
d\omega^x\equiv\;\mu^x_X\wedge\omega^x+\mu^x_U\wedge\omega^u+\mu^x_P\wedge\omega^p,
\nonumber
\\
d\omega^u\equiv\;\mu^u_U\wedge\omega^u+\omega^x\wedge\omega^p,
\nonumber
\\
d\omega^p\equiv\;\mu^u_{XU}\wedge\omega^u+\big(\mu^u_U-\mu^x_X\big)\wedge\omega^p.
\label{ODE involutive system}
\end{gather}
These equations are equivalent to the structure equations~\cite[equation~(11.5)]{O-1995} obtained using Cartan's approach.
The correspondence is given by
\begin{gather*}
\theta^1\leftrightarrow\omega^u,
\qquad
\theta^2\leftrightarrow\omega^p,
\qquad
\theta^3\leftrightarrow\omega^x,
\qquad
\pi^1\leftrightarrow\mu^u_U,
\qquad
\pi^2\leftrightarrow\mu^u_{XU},
\\
\pi^3\leftrightarrow\mu^u_U-\mu^x_X,
\qquad
\pi^4\leftrightarrow\mu^x_U,
\qquad
\pi^5\leftrightarrow\mu^x_P.
\end{gather*}

By Cartan's involutivity test,~\cite{BCGGG-1991,K-1989,O-1995}, the structure equations~\eqref{ODE involutive system}
are involutive.
In Section~\ref{involution - section}, following Seiler's book~\cite{S-2009}, Cartan's test based on the algebraic
theory of involution is introduced.
This of\/fers an alternative way of verifying, for example, that the structure equations~\eqref{ODE involutive system}
are involutive.

\end{Example}

\section{Algebra of dif\/ferential invariants}
\label{algebra of differential invariants - section}

As seen in Section~\ref{Section Equivalence Problems}, the fundamental basis theorem is at the heart of the local
equivalence problem solution.
Following~\cite{OP-2009-1}, we now introduce the algebraic constructions used to prove the fundamental basis theorem
for Lie pseudo-groups acting freely and regularly on dense open subsets of $\J^\infty$.
Once this is done we explain how to modify the algebraic constructions to take into account Lie pseudo-groups acting
freely and regularly on invariant regular subbundles of~$\J^\infty$.
Further modif\/ications will allow us to deal with Lie pseudo-groups acting regularly and non-freely on invariant
regular subbundles of $\J^\infty$.
The main conclusion is that, with the appropriate modif\/ications and regularity assumptions, the algebraic proof of
the fundamental basis theorem given in~\cite{OP-2009-1} extends to Lie pseudo-groups acting regularly and freely (or
non-freely) on invariant regular subbundles of $\J^\infty$.

\subsection{Regular submanifold jets}
%\label{section algebra differential invariants}

Let $\G$ be a~Lie pseudo-group acting on an $m$-dimensional manifold $M$.
For the moment we restrict our considerations to a~dense open subset $\mathcal{V}^\infty\subset\J^\infty$ where the
action becomes free at order $n^\star$.

Let $\mathbb{R}[t,T]$ denote the algebra of real polynomials in the variables $t=(t_1,\ldots,t_m)$ and
$T=(T^1,\ldots,T^m)$, and def\/ine
\begin{gather*}
\mathcal{T}=\left\{\eta(t,T)=\sum_{a=
1}^m \eta_a(t)T^a\right\}\simeq\mathbb{R}[t]\otimes\mathbb{R}^m \subset \mathbb{R}[t,T]
\end{gather*}
to be the $\mathbb{R}[t]$ module consisting of homogeneous linear polynomials in the variable $T$.
Let $\mathcal{T}^n\subset\mathcal{T}$ denote the subspace of homogeneous polynomials of degree $n$ in $t$.
The notations $\mathcal{T}^{\leq n}=\oplus_{k=0}^n\mathcal{T}^k$ and $\mathcal{T}^{\geq
n}=\oplus_{k=n}^\infty\mathcal{T}^k$ are used to denote the space of polynomials of degree $\leq n$ and $\geq n$ in $t$.
Let $\bH\colon\mathcal{T}\to\mathcal{T}$ be the \emph{highest order term} operator def\/ined by the requirement that
for $0\neq\eta\in\mathcal{T}^{\leq n}$ with $0\neq\bH(\eta)\in\mathcal{T}^n$ the equality $\eta=\bH(\eta)+\lambda$
holds for some $\lambda\in\mathcal{T}^{\leq n-1}$.
Locally, let $(\J^\infty TM)^*\simeq M\times\mathcal{T}$ via the pairing $\langle\j_\infty\vv;t_B T^a\rangle=\zeta^a_B$.
Then, an $n^\text{th}$ order linear dif\/ferential equation
\begin{gather}
\label{linear differential equation}
L(z,\zeta\n)=\sum_{a=1}^m\sum_{\#B\leq n} h^B_a(z) \zeta^a_B=0
\end{gather}
can be identif\/ied with the parametrized polynomial
\begin{gather}
\label{defining polynomial}
\eta(z;t,T)=\sum_{a=1}^m\sum_{\#B\leq n} h^B_a(z) t_B T^a \in \mathcal{T}^{\leq n}.
\end{gather}
\begin{Definition}
The \emph{symbol} $\bS(L)$ of the $n^\text{th}$ order (non-zero) linear dif\/ferential equation~\eqref{linear differential equation}
consists of the highest order terms in its def\/ining polynomial~\eqref{defining polynomial}:
\begin{gather*}
%\label{linear differential polynomial}
\bS[L(z,\zeta\n)]=\mathbf{H}[\eta(z;t,T)]=\sum_{a=1}^m\sum_{\#B=n} h^B_a(z) t_BT^a \in \mathcal{T}^n.
\end{gather*}
\end{Definition}

Let $\mathcal{L}=(\J^\infty\mathfrak{g})^\bot\subset(\J^\infty TM)^*$ denote the \emph{annihilator subbundle} of the
inf\/initesimal generator jet bundle, and def\/ine
\begin{gather}
\label{I module}
\mathcal{I}=\bH(\mathcal{L})
\end{gather}
to be the span of the highest order terms of the annihilating polynomials at each $z\in M$.
Since the inf\/initesimal determining system~\eqref{infinitesimal determining system} is formally integrable, it
follows that the system is closed under dif\/ferentiation with respect to the total derivative operators
$\mathbb{D}_{z^a}$.
At the symbol level, since total dif\/ferentiation with respect to $\mathbb{D}_{z^a}$ corresponds to multiplication by
$t_a$
\begin{gather*}
\bS\big(\mathbb{D}_{z^a}L\big)=t_a \bS(L),
\qquad
a=1,\ldots,m,
\end{gather*}
where $L$ is of the form~\eqref{linear differential equation}, it follows that at each point $z\in M$ the f\/iber
$\mathcal{I}|_z$ forms a~graded submodule of $\mathcal{T}$.
This submodule is called the \emph{symbol module} of the Lie pseudo-group at the point $z$.

We now introduce a~second module called the \emph{prolonged symbol} submodule of the prolonged inf\/initesimal
generator~\eqref{prolonged vector field}.
Introducing the variables $s=(s_1\ldots,s_p)$ and $S=(S^1,\ldots,S^q)$, let
\begin{gather*}
\widehat{\mathcal{S}}=\left\{\widehat{\sigma}(s,S)=\sum_{\alpha=
1}^q\widehat{\sigma}_\alpha(s)S^\alpha\right\}\simeq\mathbb{R}[s]\otimes\mathbb{R}^q
\end{gather*}
denote the $\mathbb{R}[s]$ module of polynomials that are linear in $S$.
The module $\widehat{\mathcal{S}}$ is extended to
\begin{gather*}
\mathcal{S}=\mathbb{R}^p\oplus\widehat{\mathcal{S}}=\sum_{n=-1}^\infty\mathcal{S}^n
\end{gather*}
by introducing
$\mathcal{S}^{-1}=\{c\cdot\widetilde{s}=c_1\widetilde{s}_1+\cdots+c_p\widetilde{s}_p\}\simeq\mathbb{R}^p$, where
$\widetilde{s}=(\widetilde{s}_1,\ldots,\widetilde{s}_p)\in\mathbb{R}^p$ are extra variables.
The space $\mathcal{S}$ is endowed with the structure of an $\mathbb{R}[s]$ module by taking the usual module structure
on $\widehat{\mathcal{S}}$ and setting
\begin{gather*}
\tau(s) \widetilde{s}_i=\tau(0) \widetilde{s}_i
\qquad
\text{for any polynomial}
\quad
\tau(s)\in\mathbb{R}[s].
\end{gather*}
The space $\mathcal{S}$ is called the \emph{submanifold jet module}.
A highest order term operator $\bH\colon\mathcal{S}\to\mathcal{S}$ is also introduced on $\mathcal{S}$.
For $\sigma(\widetilde{s},s,S)=c\cdot\widetilde{s}+\widehat{\sigma}(s,S)$ with $\widehat{\sigma}(s,S)\neq0$ we require
that
\begin{gather*}
\bH[\sigma(\widetilde{s},s,S)]=\bH[\widehat{\sigma}(s,S)].
\end{gather*}
In other words, the elements of $\mathcal{S}^{-1}$ have zero highest order term.
We also f\/ix a~convenient degree compatible term ordering on the monomials of $\mathcal{S}$.
For example, one could choose the degree lexicographic ordering,~\cite{S-2009}.

Given an arbitrary vector f\/ield
\begin{gather*}
\VV=\sum_{i=1}^p \xi^i\pp{}{x^i}+\sum_{\alpha=1}^q\sum_{\#J\geq0}\phi^J_\alpha\pp{}{u^\alpha_J} \in T\J^\infty,
\end{gather*}
whose coef\/f\/icients do not have to be the coef\/f\/icients of the prolonged vector f\/ield~\eqref{prolonged vector field},
the cotangent bundle $T^*\J^\infty$ is identif\/ied with $\J^\infty\times\mathcal{S}$ via the pairing
\begin{gather}
\label{S pairing}
\langle\VV;\widetilde{s}_i\rangle=\xi^i,
\qquad
\langle\VV;S^\alpha\rangle=Q^\alpha=\phi_\alpha-\sum_{i=1}^p u^\alpha_i \xi^i,
\qquad
\langle\VV;s_JS^\alpha\rangle=\phi^J_\alpha,
\qquad
\#J\geq1.
\end{gather}
Fixing $z\ii\in\mathcal{V}^\infty\subset\J^\infty$ with $\pi^\infty_0\big(z\ii\big)=z$, the prolongation map
$\mathbf{p}=\mathbf{p}\ii_{z\ii}\colon\J^{\infty}TM|_z\to T\J^\infty|_{z\ii}$ given in~\eqref{prolongation map} induces
the \emph{dual prolongation map} $\bp^*\colon\mathcal{S}\to\mathcal{T}$ def\/ined by
\begin{gather*}
\langle\j_\infty\vv;\bp^*(\sigma)\rangle=\langle\bp(\j_\infty\vv);\sigma\rangle=\langle\vv\ii;\sigma\rangle
\qquad
\text{for all}
\quad
\j_\infty\vv\in\J^\infty TM|_z
\quad
\text{and}
\quad
\sigma\in\mathcal{S}.
\end{gather*}

Next, let
\begin{gather}
\label{particular polynomials}
\beta_i(t)=t_i+\sum_{\alpha=1}^q u^\alpha_i t_{p+\alpha},\qquad i=1,\ldots,p,
\\
B^\alpha(T)=T^{p+\alpha}-\sum_{i=1}^p u^\alpha_i T^i,
\qquad
\alpha=1,\ldots,q,
\end{gather}
where $u^\alpha_i=\partial u^\alpha/\partial x^i$ are the f\/irst order jet coordinates of the f\/ixed submanifold jet
$z\ii$.
Geo\-metrically, the polynomial $B^\alpha(T)$ is the symbol of the characteristic component $Q^\alpha$ appearing
in~\eqref{S pairing} while $\beta_i(t)$ represents the symbol of the $i^\text{th}$ total derivative operator $D_{x^i}$:
\begin{gather*}
\bS\big(D_{x^i}L\big)=\beta_i(t) \bS(L),
\end{gather*}
where $L$ is any linear dif\/ferential equation in the vector f\/ield jet coordinates.
For f\/ixed f\/irst order jet coordinates $u^\alpha_i$, the polynomials~\eqref{particular polynomials} def\/ine
a~linear map
\begin{gather}
\label{beta map}
\bb\colon \ \mathbb{R}^{2m}\to\mathbb{R}^m
\qquad
\text{given by}
\qquad
s_i=\beta_i(t),
\qquad
S^\alpha=B^\alpha(T),
\end{gather}
and for $\widehat{\sigma}(s,S)\in\widehat{\mathcal{S}}\subset\mathcal{S}$ the equality
\begin{gather}
\label{p and beta permutation}
\bH[\bp^*(\widehat{\sigma})]=\bb^*[\bH(\widehat{\sigma})]
\end{gather}
holds.
\begin{Definition}
The \emph{prolonged annihilator subbundle} of the prolonged Lie algebra $\mathfrak{g}\ii$ at
$z\ii\in\mathcal{V}^\infty$ is
\begin{gather}
\label{prolonged annihilator subbundle}
\mathcal{Z}=\big(\mathfrak{g}\ii|_{z\ii}\big)^\bot=(\bp^*)^{-1}\mathcal{L}|_z\subset\mathcal{S}.
\end{gather}
Furthermore, let
\begin{gather*}
\mathcal{U}=\bH(\mathcal{Z})\subset\mathcal{S}
\end{gather*}
be the span of the highest order terms of the prolonged annihilator subbundle (in general, $\mathcal{U}$ is not
a~submodule).
\end{Definition}
\begin{Definition}
The \emph{prolonged symbol submodule} is def\/ined as the inverse image of the symbol module~\eqref{I module} under the
polynomial pull-back morphism~\eqref{beta map}:
\begin{gather}
\label{definition J}
\mathcal{J}=\big(\bb^*\big)^{-1}(\mathcal{I}|_z)=
\big\{\widehat{\sigma}(s,S)\in\widehat{\mathcal{S}}: \, \bb^*(\widehat{\sigma}(s,S))=
\widehat{\sigma}(\beta(t),B(T))\in\mathcal{I}|_z\big\}\subset\widehat{\mathcal{S}}.
\end{gather}
\end{Definition}

From~\eqref{p and beta permutation} and~\eqref{prolonged annihilator subbundle}, the containment
$\mathcal{U}\subset\mathcal{J}$ always holds.
When the action is locally free the containment becomes an equality,~\cite{OP-2009-1}.
\begin{Lemma}
\label{lemma equality U and J}
Let $n^\star$ be the order of freeness.
Then for all $n>n^\star$
\begin{gather*}
\mathcal{U}^n|_{z^{(n)}}=\mathcal{J}^n|_{z^{(n)}},
\qquad
z\n\in\mathcal{V}^n.
\end{gather*}
\end{Lemma}

Given a~moving frame $\widehat{\rho}\colon\mathcal{V}^\infty\to\mathcal{B}\ii$, the invariantization map is used to
invariantize the preceding algebraic constructions.
For example, the invariantization of a~section
\begin{gather*}
\eta(x,u;t,T)=\sum_{a=1}^m\sum_{\#B\leq n}h^B_b(x,u) t_B T^a
\end{gather*}
of the annihilator bundle $\mathcal{L}$ is the polynomial
\begin{gather*}
\overline{\eta}(X,U;t,T)=\iota[\eta(x,u;t,T)]=\sum_{a=1}^m\sum_{\#B\leq n}h^B_a(X,U) t_B T^a
\end{gather*}
obtained by replacing the coordinates $(x,u)$ on $M$ by their invariantizations $(X,U)=\iota(x,u)$.
Similarly, the invariantization of a~prolonged symbol polynomial
\begin{gather*}
\widehat{\sigma}\big(x,u^{(k)};s,S\big)=\sum_{\alpha=
1}^q\sum_{\#J\leq n} h^J_\alpha(x,u^{(k)}) s_J S^\alpha \in \widehat{\mathcal{S}}^{\leq n}
\end{gather*}
is the polynomial
\begin{gather}
\label{invariantization sigma}
\overline{\sigma}\big(X,U^{(k)},s,S\big)=\iota\big[\widehat{\sigma}(x,u^{(k)};s,S)\big]
=\sum_{\alpha=1}^q\sum_{\#J\leq n} h^J_\alpha\big(X,U^{(k)}\big) s_J S^\alpha.
\end{gather}

Let $\overline{\mathcal{U}}^n\big|_{z^{(n)}}=\iota(\mathcal{U}^n)\big|_{z\n}$ and
$\overline{\mathcal{J}}^n\big|_{z^{(n)}}=\iota(\mathcal{J}^n)\big|_{z\n}$.
The equality $\overline{\mathcal{U}}^n\big|_{z^{(n)}}=\overline{\mathcal{J}}^n\big|_{z^{(n)}}$ for all $n>n^\star$ is the key
to proving the fundamental basis theorem.
Since $\overline{\mathcal{J}}^{>n^\star}\big|_{z\n}$ is a~polynomial ideal it has a~\emph{Gr\"obner
basis},~\cite{CLO-1996}, which brings algebraic structure into the problem.
After identifying the polynomial~\eqref{invariantization sigma} with the dif\/ferential invariant
\begin{gather*}
I_{\overline{\sigma}}=\sum_{\alpha=1}^q\sum_{\#J\geq0}h^J_\alpha\big(X,U^{(k)}\big) U^\alpha_J
\end{gather*}
the following result was proved in~\cite{OP-2009-1}.
\begin{Theorem}
\label{basis theorem}
Let $\G$ be a~Lie pseudo-group acting freely and regularly on $\mathcal{V}^{n^\star}\subset{\rm J}^{n^\star}$, where
$n^\star$ is the order of freeness.
Then a~finite generating set for the algebra of differential invariants consists of
\begin{itemize}\itemsep=0pt
\item the differential invariants $I_\nu=I_{\overline{\sigma}_\nu}$, where
$\overline{\sigma}_1,\ldots,\overline{\sigma}_l$ form a~Gr\"obner basis for the inva\-riantized prolonged symbol
submodule $\overline{\mathcal{J}}^{>n^\star}$, and, possibly, \item a~finite number of additional differential
invariants of order $\leq n^\star$.
\end{itemize}
\end{Theorem}

\begin{Remark}
\label{minimality remark}
While Theorem~\ref{basis theorem} yields a~generating set for the algebra of dif\/ferential inva\-riants, it certainly
does not imply that this set is minimal.
Unfortunately, there is, to this day, no algorithm for determining whether a~generating set containing more than one
dif\/ferential invariant is minimal.
\end{Remark}

The algebraic considerations introduced above and the proof of Theorem~\ref{basis theorem} assumes the pseudo-group
acts freely and regularly on a~dense open subset $\mathcal{V}^\infty\subset\J^\infty$.
In Example~\ref{regular moving frame example}, this would correspond to the set of regular jets $\mathcal{V}^\infty_1$.
For the other regular subbundles $\mathcal{V}^\infty_2, \ldots, \mathcal{V}^\infty_5$ of Example~\ref{regular
moving frame example}, we must incorporate into the algebraic constructions the dependencies among the submanifold jet
coordinates introduced by the determining equations~\eqref{i-th determining system}.
To achieve this, a~second identif\/ication of the one-forms on $\J^\infty$ with polynomials in $\mathcal{S}$ is
introduced,~\cite{OP-2009-1}.
This identif\/ication is given by
\begin{gather}
\label{new S identification}
dx^i\longleftrightarrow\widetilde{s}_i,
\qquad
du^\alpha_J\longleftrightarrow s_J S^\alpha.
\end{gather}
Under the correspondence~\eqref{new S identification}, the exterior dif\/ferential of an $n^\text{th}$ order
dif\/ferential equation $E\big(x,u\n\big)=0$
\begin{gather*}
dE\big(x,u\n\big)=\sum_{i=1}^p \pp{E}{x^i}\big(x,u\n\big) dx^i+\sum_{\alpha=
1}^q\sum_{\#J\leq n}\pp{E}{u^\alpha_J}\big(x,u\n\big) du^\alpha_J=0
\end{gather*}
at a~submanifold jet $z\n=\big(x,u\n\big)$ can be identif\/ied with the parametrized polynomial
\begin{gather}
\label{symbol E}
\alpha(z\n;\widetilde{s},s,S)=\sum_{i=1}^p \pp{E}{x^i}\big(x,u\n\big) \widetilde{s}_i+\sum_{\alpha=
1}^q\sum_{\#J\leq n}\pp{E}{u^\alpha_J}\big(x,u\n\big) s_JS^\alpha \in \mathcal{S}^{\leq n}.
\end{gather}
\begin{Definition}
The \emph{symbol} $\bS\big(E\big(z\n\big)\big)$ of an $n^\text{th}$ order (non-zero) dif\/ferential equation $E\big(z\n\big)$ $=E\big(x,u\n\big)=0$ at
a~submanifold jet $z\n=\big(x,u\n\big)$ consists of the highest order terms in the parametrized polynomial~\eqref{symbol E}:
\begin{gather*}
\bS\big[E\big(z\n\big)\big]=\mathbf{H}\big[\alpha(z\n;\widetilde{s},s,S)\big]=\sum_{\alpha=1}^q\sum_{\#J=
n} \pp{E}{u^\alpha_J}\big(z\n\big) s_J S^\alpha \in \mathcal{S}^n.
\end{gather*}
\end{Definition}

Now, let $\mathcal{V}^\infty_i\subset\J^\infty$ be an invariant regular subbundle with determining system
\begin{gather}
\label{E-infty}
E\ii_i\big(x,u\ii\big)=0.
\end{gather}
Since the system of dif\/ferential equations~\eqref{E-infty} is formally integrable, its symbol at a~f\/ixed
submanifold jet $z\ii=(x,u\ii)$
\begin{gather*}
\mathcal{E}_i|_{z\ii}=\bS\big[E\ii_i\big(z\ii\big)\big]=\mathbf{H}\big[\big(T_{z\ii}\mathcal{V}^\infty_i\big)^\bot\big]
\end{gather*}
is a~submodule of the submanifold jet module $\mathcal{S}$.
Let $\mathcal{M}_i|_{z\ii}\subset\mathcal{S}$ denote the monomial module generated by the leading monomials (with
respect to a~chosen term ordering on $\mathcal{S}$) of the symbol module $\mathcal{E}_i|_{z\ii}$.
We can assume, possibly by restricting to an open subset and employing $\delta$-regular coordinates (see
Def\/inition~\ref{delta-regular} below) that generically $\mathcal{M}_i|_{z\ii}=\mathcal{M}_i$ does not depend upon
$z\ii\in\mathcal{V}^\infty_i$.
Let
\begin{gather*}
\mathcal{S}_i=\mathcal{S}^{-1}\oplus\text{span}\;\{s_J S^\alpha\notin\mathcal{M}_i\} \subset \mathcal{S}
\end{gather*}
denote the \emph{$i^\text{th}$ restricted submanifold jet module} spanned by all monomials \emph{not} in the monomial
module $\mathcal{M}_i$ to which we add $\mathcal{S}^{-1}$.
Applying standard Gaussian elimination we can construct a~linear basis for the space $\mathcal{E}_i|_{z\ii}$ of the form
\begin{gather}
\label{linear basis}
s_I S^\beta+\sum_{s_J S^\alpha\in \mathcal{S}^n_i} h^J_\alpha\big(z\n\big) s_J S^\alpha
\qquad
\text{for all}
\quad
s_IS^\beta\in\mathcal{M}^n_i,
\quad
n\geq0.
\end{gather}
A similar statement holds for the subspace $(T_{z\ii}\mathcal{V}^\infty_i)^\bot$, where the sum in~\eqref{linear basis}
now runs over all the monomials in $\mathcal{S}^{\leq n}_i$.
Therefore, $\mathcal{S}_i$ is a~f\/ixed complement to the symbol module $\mathcal{E}_i|_{z\ii}$ and the annihilating
subspace $(T_{z\ii}\mathcal{V}^\infty_i)^\bot$:
\begin{gather*}
\mathcal{S}=\mathcal{S}_i\oplus\mathcal{M}_i=\mathcal{S}_i\oplus\mathcal{E}_i|_{z\ii}=
\mathcal{S}_i\oplus\big(T_{z\ii}\mathcal{V}^\infty_i\big)^\bot.
\end{gather*}
We can thus identify the submodule
\begin{gather*}
\mathcal{S}_i\simeq\mathcal{S}/(\mathcal{E}_i|_{z\ii})\simeq\mathcal{S}/\big(T_{z\ii}\mathcal{V}^\infty_i\big)^\bot\simeq T^*_{z\ii}\mathcal{V}^\infty_i
\end{gather*}
with the dual space of $\mathcal{V}^\infty_i$ at the point $z\ii\in\mathcal{V}^\infty_i$.
Under the identif\/ication~\eqref{new S identification}, the monomials $s_J S^\alpha\in\mathcal{S}_i$ indicate the
\emph{parametric jet variables} $u_J^\alpha$ on $\mathcal{V}^\infty_i$.

On each invariant regular subbundle $\mathcal{V}^\infty_i$, the algebraic constructions introduced before
Remark~\ref{minimality remark} still hold provided the sub\-mani\-fold jet module $\mathcal{S}$ is replaced by the
restricted sub\-mani\-fold jet module $\mathcal{S}_i$.
Also, since $\mathcal{V}^\infty_i$ is $\G$-invariant, the determining system~\eqref{E-infty} is invariant and the
invariantization of the algebraic constructions is well def\/ined.
In particular, the fundamental basis Theorem~\ref{basis theorem} still holds (with $n^\star$ replaced by $n^\star_i$
and $\overline{\mathcal{J}}^{>n^\star}$ replaced by $\overline{\mathcal{J}}^{>n^\star_i}_i$, where
$\mathcal{J}_i=(\bb^*)^{-1}(\mathcal{I}|_z)\subset\widehat{\mathcal{S}}_i$).
\begin{Example}
Continuing Examples~\ref{regular moving frame example} and~\ref{Cartan analogy example}, since in the generic case
$\mathcal{V}^\infty_1$ is a~dense open subset of $\J^\infty$ the restricted submanifold jet module is the whole
submanifold jet module: $\mathcal{S}_1=\mathcal{S}$.
On the other hand, on $\mathcal{V}^\infty_2$ the dif\/ferential constraint $Q_{P^4}\equiv0$ implies that
\begin{gather}
\label{case 2 determining system}
Q_{X^i U^j P^{k+4}}\equiv0
\qquad
\text{or equivalently}
\qquad
q_{x^i u^j p^{k+4}}\equiv0,
\qquad
i,j,k\geq0.
\end{gather}
The symbol associated to the determining system~\eqref{case 2 determining system} is
\begin{gather*}
\mathcal{E}_2=\mathcal{M}_2=\text{span}\big\{s_x^i s_u^j s_p^{k+4}S:\;i,j,k\geq0\big\},
\end{gather*}
and the corresponding restricted submanifold jet module is
\begin{gather*}
\mathcal{S}_2=
\text{span}\{\widetilde{s}_x, \widetilde{s}_u, \widetilde{s}_p\}\oplus\text{span}\big\{s_x^i s_u^j s_p^k S:\;0\leq k\leq3
\;
\text{and}
\;
i,j\geq0\big\}.
\end{gather*}
Similar restricted submanifold jet modules $\mathcal{S}_3, \ldots, \mathcal{S}_5$ for the subbundles
$\mathcal{V}^\infty_3, \ldots, \mathcal{V}^\infty_5$ can be obtained.
\end{Example}

\subsection{Singular submanifold jets}
%\label{section algebra of differential invariants singular jets}

We now assume that $z\ii\in\mathscr{S}^\infty_i\subset\J^\infty$ is a~singular submanifold jet where the pseudo-group
$\G$ acts regularly but not freely.
As in the previous section, let $\mathcal{S}_i$ denote the restricted submanifold jet module associated with the
regular invariant subbundle $\mathscr{S}^\infty_i$.
Let
\begin{gather*}
%\label{isotropy algebra}
\j_\infty\mathfrak{g}_{z\ii}=\text{ker}\,\bp|_{z\ii}\cap\J^\infty\mathfrak{g}|_z,
\qquad
z=\pi^\infty_0\big(z\ii\big)
\end{gather*}
be the vector f\/ield jets of the isotropy Lie algebra at $z\ii$.
Identifying $(\J^\infty T_M)^*|_z$ with the symbol module $\mathcal{T}$, we introduce the \emph{isotropy annihilator
vector space}
\begin{gather*}
(\j_\infty\mathfrak{g}_{z\ii})^\bot=\mathcal{T}_{i}\subset\mathcal{T}.
\end{gather*}
To streamline the notation we have suppressed the dependence of $\mathcal{T}_i$ on the submanifold jet
$z\ii\in\mathscr{S}^\infty_i$.
\begin{Proposition}
\label{proposition annihilator}
Let $n^\star_i$ be the order of partial freeness of $\G$ on $\mathscr{S}^\infty_i$.
Then for $n\geq n^\star_i$ and $z\n\in\mathscr{S}^n_i$,
\begin{gather}
\label{dual g freeness}
\bp^*\big(\mathcal{S}^{\leq n}_i\big)+\mathcal{L}^{\leq n}|_z=\mathcal{T}^{\leq n}_i.
\end{gather}
\end{Proposition}
\begin{proof}
For $z\n\in\mathscr{S}^n_i$ and $n\geq n^\star_i$ we have
\begin{gather*}
\j_n\mathfrak{g}_{z\n}=\text{ker}\,\bp\n|_{z\n}\cap\J^n\mathfrak{g}|_z=
\big(\text{rng}\big(\bp\n\big)^*\big)^\bot\cap\big(\mathcal{L}^{\leq n}|_z\big)^\bot=
\big(\bp^*\big(\mathcal{S}^{\leq n}_i\big)+\mathcal{L}^{\leq n}|_z\big)^\bot,
\end{gather*}
from which~\eqref{dual g freeness} follows.
\end{proof}

Lemma~\ref{lemma equality U and J}, essential to the proof of the fundamental basis Theorem~\ref{basis theorem}, also
holds at a~singular submanifold jet.
\begin{Proposition}
\label{proposition equality U and J}
For $n>n_i^\star$ and $z\n\in\mathscr{S}^n_i$,
\begin{gather}
\label{equality of modules}
\mathcal{U}^n_i|_{z\n}=\mathcal{J}^n_i|_{z\n}.
\end{gather}
\end{Proposition}
\begin{proof}
By induction, it suf\/f\/ices to prove~\eqref{equality of modules} when $n=n^\star_i{+}1$.
Let $Q\in\mathcal{J}^{n^\star_i+1}_i|_{z^{(n^\star_i+1)}}$ and $P=\bp^*(Q)$.
By~\eqref{p and beta permutation} and~\eqref{definition J}
\begin{gather*}
\bH(P)=\bH[\bp^*(Q)]=\bb^*[\bH(Q)]=\bb^*(Q) \in \mathcal{I}^{n^\star_i+1}|_z,
\end{gather*}
and we conclude that there exists $Y\in\mathcal{T}^{\leq n^\star_i}$ such that $P+Y\in\mathcal{L}^{n^\star_i+1}|_z$.
Actually, we have that $Y\in\mathcal{T}^{\leq n^\star_i}_i$.
To see this, we f\/irst observe that the formal integrability of the determining equations of the isotropy group
$\G^{(n^\star_i)}_{z^{(n^\star_i)}}$ implies that the projection $\mathcal{T}_i^{\leq
n^\star_i+1}\to\mathcal{T}_i^{\leq n^\star_i}$ is surjective.
Hence, by Proposition~\ref{proposition annihilator}, since $P+Y$ and $P$ are in $\mathcal{T}_i^{\leq n_i^\star+1}$ we
have that
\begin{gather*}
Y=(P+Y)-P \in \mathcal{T}_i^{\leq n^\star_i}.
\end{gather*}
Now, let $U\in\mathcal{S}_i^{\leq n^\star_i}$ and $V\in\mathcal{L}^{\leq n^\star_i}\big|_z$ such that $Y=\bp^*(U)+V$, then
\begin{gather*}
\bp^*(Q+U)=(P+Y)-V \in \mathcal{L}^{\leq n^\star_i+1}\big|_{z}.
\end{gather*}
Finally, equation~\eqref{prolonged annihilator subbundle} implies that $Q+U\in\mathcal{Z}^{\leq
n^\star_i{+}1}_i\big|_{z^{(n^\star_i{+}1)}}$.
\end{proof}
\begin{Remark}
The proof of Proposition~\ref{proposition equality U and J} is essentially the same as~\cite[Lemma 5.5]{OP-2009-1}.
It is included to show that formal integrability of the determining equations of the isotropy group
$\G^{(n^\star_i)}_{z^{(n^\star_i)}}$ is essential for the proof to remain valid at a~submanifold jet where the
pseudo-group does not act freely.
Due to Proposition~\ref{proposition equality U and J}, the constructive proof of the fundamental basis theorem given
in~\cite{OP-2009-1} still holds at a~singular submanifold jet $z\ii\in\mathscr{S}^\infty_i$ (with the necessary
algebraic modif\/ications as in the regular case).
\end{Remark}

\section{Involutivity}
\label{involution - section}

In Section~\ref{section equivalence of coframes} we showed via examples that in the appropriate geometrical setting,
the pull-back of the structure equations of an equivalence pseudo-group by a~(partial) moving frame reproduces Cartan's
moving frame results.
In particular, in Example~\ref{involutive example} we recovered the involutive structure equations~\eqref{ODE
involutive system}.
The aim of this f\/inal section is to complete Section~\ref{section equivalence of coframes} by recasting Cartan's
involutivity test within the algebraic framework of Section~\ref{algebra of differential invariants - section}.
Involutivity plays an essential role in the solution of equivalence problems,~\cite{O-1995}.
It guarantees that the equivalence map constructed by specifying its jets (or Taylor series coef\/f\/icients)
converges, and the Cartan characters give the ``dimensional freedom'' of the equivalence map in
Theorem~\ref{equivalence theorem for singular submanifolds}.
The following exposition follows Seiler's book,~\cite{S-2009}.
\begin{Definition}
Let
\begin{gather*}
\eta_\upsilon\big(z\n;t,T\big)=\sum_{a=1}^m\sum_{\#B=n} h^B_{a;\upsilon}\big(z\n\big) t_B T^a,
\qquad
\upsilon=1,\ldots,\ell,
\qquad
n\geq1,
\end{gather*}
be a~basis for the degree $n$ homogeneous component of $\mathbf{H}(\mathcal{T}_i^{\leq n})$.
The \emph{symbol matrix}
\begin{gather*}
%\label{symbol matrix}
\bT^n_i=\big(h^B_{a;\upsilon}\big(z\n\big)\big)
\end{gather*}
is the $\ell\times m\binom{m+n-1}{n}$ matrix where the entries of the $\upsilon^\text{th}$ row are given by the
coef\/f\/icients $h^B_{a;\upsilon}\big(z\n\big)$ of the polynomial $\eta_\upsilon(z\n;t,T)$.
\end{Definition}

To def\/ine the \emph{class} of a~symmetric multi-index $B=(b^1,\ldots,b^n)$ of order $\#B=n$ we rewrite the
multi-index as $\widetilde{B}=(\widetilde{b}^1,\ldots,\widetilde{b}^m)$, where $\widetilde{b}^a$ is the number of
occurrences of the integer $1\leq a~\leq m$ in $B=(b^1,\ldots,b^n)$.
\begin{Definition}
The \emph{class} of a~multi-index $\widetilde{B}=(\widetilde{b}^1,\ldots,\widetilde{b}^m)$ is
$%\begin{gather*}
\text{cl}\, \widetilde{B}=\min \big\{a:\widetilde{b}^a\neq0\big\}.
$%\end{gather*}
\end{Definition}

Next, the columns of the symbol matrix $\bT^n_i$ are ordered in such a~way that the column
$(h^{\widetilde{B}}_{a;1},\ldots,h^{\widetilde{B}}_{a;\ell})^T$ is always to the left of the column
$(h^{\widetilde{A}}_{b;1},\ldots$, $h^{\widetilde{A}}_{b;\ell})^T$ if $\text{cl}\,\widetilde{B}>\text{cl}\,\widetilde{A}$.
For two multi-indices with same class, the order of the columns does not matter.
Once the columns of the symbol matrix are ordered, the matrix is put in row echelon form without performing any column
permutations.
\begin{Definition}
\label{index definition}
Let $\beta^{(a)}_n$, $a=1,\ldots,m$, be the number of pivots with class $1\leq a~\leq m$ of the row echelon form symbol
matrix $\bT^n_i$.
The numbers $\beta^{(a)}_n$ are called the {\it indices} of $\bT^n_i$.
\end{Definition}

Def\/inition~\ref{index definition} depends on the chosen coordinate system and one must work with $\delta$-regular
coordinate systems.
\begin{Definition}
\label{delta-regular}
A coordinate system is said to be $\delta$-\emph{regular} if the sum $\sum\limits_{a=1}^m a \beta^{(a)}_n$ is
maximal.
\end{Definition}

Any coordinate system can be transformed into a~$\delta$-regular one with a~linear transformation def\/ined by a~matrix
coming from a~Zariski open subset of $\mathbb{R}^{m\times m}$,~\cite{S-2009}.
\begin{Definition}
The degree $n$ homogeneous component of $\mathbf{H}\big(\mathcal{T}^{\leq n}_i\big)$ is said to be \emph{involutive} if the
symbol matrix $\bT^{n+1}_i$ of the degree $n+1$ homogeneous component of $\mathbf{H}\big(\mathcal{T}^{\leq n+1}_i\big)$
satisf\/ies the algebraic equality
\begin{gather}
\label{involutivity criterion}
\text{rank}\,\mathbf{T}^{n+1}_i=\sum_{a=1}^m a \beta^{(a)}_n.
\end{gather}
\end{Definition}
\begin{Definition}
Let $n\geq1$, the isotropy annihilator subbundle $\mathcal{T}_i^{\leq n}$ is in said to be involutive if the degree $n$
homogeneous component of $\mathbf{H}\big(\mathcal{T}_i^{\leq n}\big)$ is involutive and the projection
$\mathcal{T}_i^{\leq{n+1}}\to\mathcal{T}_i^{\leq n}$ is surjective.
\end{Definition}

A standard result from the theory of involutivity guarantees that when $\mathcal{T}_i^{\leq n}$ becomes involutive then
for all $k>n$ the isotropy annihilator subbundle $\mathcal{T}_i^{\leq k}$ remains involutive.
\begin{Definition}
%\label{definition Cartan characters}
Assume $\mathcal{T}^{\leq n}_i$ is involutive with indices $\beta^{(a)}_n$.
The \emph{Cartan characters} $\alpha_n^{(a)}$ of~$\mathcal{T}^{\leq n}_i$ are def\/ined by
\begin{gather*}
%\label{Cartan characters formula}
\alpha^{(a)}_n=m\binom{n+m-a-1}{n-1}-\beta^{(a)}_n,
\qquad
1\leq a\leq m.
\end{gather*}
\end{Definition}
\begin{Theorem}
%\label{theorem isotropy size}
Let $\mathcal{T}^{\leq n}_i$ be involutive with Cartan characters $\alpha^{(a)}_n$.
Then the isotropy group of the submanifold jet $z\ii$ depends on $f_a$ arbitrary functions of $a$ variables where the
numbers $f_a$ are determined by the recursion relation
\begin{gather*}
%\label{dimensional recursion relation}
f_m=\alpha^{(m)}_n,
\\
f_a=\alpha^{(a)}_n+\sum_{b=a+1}^m \frac{(a-1)!}{(m-1)!}\Big(s^{(b-1)}_{b-a}(0) \alpha^{(b)}_n-s^{(b-1)}_{b-a}(n) f_b\Big),
\qquad
1\leq a\leq m-1,
\end{gather*}
provided all $f_a$ are non-negative integers.
The numbers $s^{(a)}_b(c)$ are the modified Stirling numbers defined by the identity
\begin{gather*}
(c+y+1)(c+y+2)\cdots(c+y+a)=\sum_{b=0}^a s^{(a)}_{a-b}(c) y^b
\end{gather*}
for all non-negative integers $a$, $b$, $c$ and $a\geq b$.
Here $y$ is an arbitrary variable.
\end{Theorem}
\begin{Example}\looseness=-1
Continuing Example~\ref{involutive example} we verify that the involutivity test~\eqref{involutivity criterion} holds
and compute the size of the isotropy group at a~submanifold jet.
First, we recall that for the contact pseudo-group~\eqref{contact} the normalized Maurer--Cartan forms, up to order 1,
satisfy the linear relations
\begin{gather}
\mu^x\equiv-\omega^x,
\qquad
\mu^u\equiv-\omega^u,
\qquad
\mu^p\equiv-\omega^p,
\qquad
\mu^u_P\equiv\mu^p_X\equiv\mu^x_Q\equiv\mu^u_Q\equiv\mu^p_Q\equiv0,
\nonumber
\\
\mu^q-\mu^p_X\equiv0,
\qquad
\mu^u_X-\mu^p\equiv0,
\qquad
\mu^p_P+\mu^x_X-\mu^u_U\equiv0,
\qquad
\mu^q_Q-\mu^p_P+\mu^x_X\equiv0.
\label{involution normalized mc}
\end{gather}
To show that the degree 1 homogeneous component of $\mathbf{H}(\mathcal{T}^{\leq1}_1)$ is involutive, we f\/irst note
that via the lift map~\eqref{lift of a vector field jet} the vector f\/ield jet $\zeta^a_B$ can be identif\/ied with
the Maurer--Cartan form $\mu^a_B$ so that we have the identif\/ication
\begin{gather}
\label{mc identification}
\mu^a_B\longleftrightarrow\zeta^a_B\longleftrightarrow t_BT^a.
\end{gather}
Secondly, since the invariantization of the algebraic constructions introduced in Section~\ref{algebra of differential invariants - section}
coincide with their progenitor when restricted to a~cross-section, we can use the
linear relations among the normalized Maurer--Cartan forms~\eqref{involution normalized mc} to verify the involutivity
test~\eqref{involutivity criterion} at a~submanifold jet on the cross-section def\/ining a~partial moving frame.
Using the correspondence~\eqref{mc identification}, the isotropy annihilator polynomials associated
to~\eqref{involution normalized mc} are
\begin{gather*}
T^x,
\qquad
T^u,
\qquad
T^p,
\qquad
t_pT^u,
\qquad
t_xT^p,
\qquad
t_qT^x,
\qquad
t_qT^u,
\qquad
t_qT^p,
\\
T^q-t_xT^p,
\qquad
t_xT^u-T^p,
\qquad
t_pT^p+t_xT^x-t_uT^u,
\qquad
t_qT^q-t_pT^p+t_xT^x.
\end{gather*}
Writing the order 1 symbol matrix $\mathbf{T}^1$ we obtain
\begin{gather*}
\bordermatrix{&\!t_qT^q\!&\!t_qT^p\!&\!t_qT^u\!&\!t_qT^x\!&\!t_xT^p\!&\!t_xT^u\!&\!t_xT^x\!&\!t_pT^p\!&\!t_pT^u\!&
\!t_pT^x\!&\!t_uT^p\!&\!t_uT^u\!&\!t_uT^x\!\cr&1&0&0&0&0&0&-1&1&0&0&0&0&0\cr&0&1&0&0&0&0&0&0&0&0&0&0&0\cr&0&0&1&0&0&0&0&0&0&0&0&0&0\cr&0&0&0&1&0&0&0&0&0&0&0&0&0\cr&0&0&0&0&1&0&0&0&0&0&0&0&0\cr&0&0&0&0&-1&0&0&0&0&0&0&0&0\cr&0&0&0&0&0&1&0&0&0&0&0&0&0\cr&0&0&0&0&0&0&1&1&0&0&0&-1&0\cr&0&0&0&0&0&0&0&0&1&0&0&0&0}
\end{gather*}
and f\/ind that the indices of the matrix are $\beta^{(4)}_1=4$, $\beta^{(3)}_1=3$, $\beta^{(2)}_1=1$ and
$\beta^{(1)}_1=0$.
At order~2, the rank of the symbol matrix $\mathbf{T}^2$ is
$27=4 \beta^{(4)}_1+3 \beta^{(3)}_1+2 \beta^{(2)}_1+\beta^{(1)}_1$ which satisf\/ies the involutivity test.
The corresponding Cartan characters are $\alpha^{(4)}_1=\alpha^{(3)}_1=0$, $\alpha^{(2)}_1=2$, $\alpha^{(1)}_1=3$, and
we conclude that the isotropy group involves two arbitrary analytic functions, each depending on two variables.
In other words, the general contact transformation between two second order ordinary dif\/ferential equations depends
on two arbitrary functions of two variables.
This is in accordance with Cartan's involutivity test based on the theory of exterior dif\/ferential
systems,~\cite[Example 11.10]{O-1995}.
\end{Example}
\begin{Example}
%\label{essential invariant example}
In this f\/inal example we consider the simultaneous equivalence of a~two-form and a~vector f\/ield on
$\mathbb{R}^3$,~\cite{GS-1992}.
This example is interesting as the solution admits a~branch with an inf\/inite-dimensional isotropy group and an
essential invariant.
Let
\begin{gather}
\label{omega}
\Omega=a(x,y,z)dx\wedge dy+b(x,y,z)dx\wedge dz+c(x,y,z)dy\wedge dz,
\qquad
a(x,y,z)\neq0,
\end{gather}
be a~non-vanishing two-form, and
\begin{gather}
\label{vector field last example}
\vv=e(x,y,z)\pp{}{x}+f(x,y,z)\pp{}{y}+g(x,y,z)\pp{}{z},
\qquad
g(x,y,z)\neq0,
\end{gather}
a non-zero vector f\/ield on $\mathbb{R}^3$.
If
\begin{gather*}
\overline{\Omega}=A(X,Y,Z)dX\wedge dY+B(X,Y,Z)dX\wedge dZ+C(X,Y,Z)dY\wedge dZ,
\qquad\!
A(X,Y,Z)\neq0,
\\
\overline{\vv}=E(X,Y,Z)\pp{}{X}+F(X,Y,Z)\pp{}{Y}+G(X,Y,Z)\pp{}{Z},
\qquad
G(X,Y,Z)\neq0,
\end{gather*}
is another set of non-vanishing two-form and vector f\/ield on $\mathbb{R}^3$, the map
\begin{gather*}
\Phi\colon
\
X=\phi(x,y,z),
\,
Y=\beta(x,y,z),
\,
Z=\alpha(x,y,z)
\in
\D\big(\mathbb{R}^3\big),
\end{gather*}
is a~local equivalence if it satisf\/ies
\begin{gather}
\label{form-vector equivalence criterion}
\Phi^*(\overline{\Omega})=\Omega
\qquad
\text{and}
\qquad
d\Phi^{-1}(\overline{\vv})=\vv.
\end{gather}

The equivalence problem splits in two cases; either $\vv\interior\Omega=0$ or $\vv\interior\Omega\neq0$.
In the following we consider the case $\vv\interior\Omega=0$.
This imposes the restrictions
\begin{gather*}
e(x,y,z)=\frac{g(x,y,z) c(x,y,z)}{a(x,y,z)}
\qquad
\text{and}
\qquad
f(x,y,z)=-\frac{g(x,y,z) b(x,y,z)}{a(x,y,z)}
\end{gather*}
on the coef\/f\/icients of vector f\/ield~\eqref{vector field last example} (and similar constraints on the
coef\/f\/icients of $\overline{\vv}$).
In local coordinates, the equivalence criterions~\eqref{form-vector equivalence criterion} yield the transformation
rules
\begin{gather}
A(\phi_x\beta_y-\beta_x\phi_y)+B(\phi_x\alpha_y-\alpha_x\phi_y)+C(\beta_x\alpha_y-\alpha_x\beta_y)=a,
\nonumber
\\
A(\phi_x\beta_z-\beta_x\phi_z)+B(\phi_x\alpha_z-\alpha_x\phi_z)+C(\beta_x\alpha_z-\alpha_x\beta_z)=b,
\nonumber
\\
A(\phi_y\beta_z-\beta_y\phi_z)+B(\phi_y\alpha_z-\alpha_y\phi_z)+C(\beta_y\alpha_z-\alpha_y\beta_z)=c,
\nonumber
\\
G=\frac{g}{a}(c \alpha_x-b \alpha_y+a \alpha_z),
\label{form-vector coordinate transformation rules}
\end{gather}
for the components of the two-form~\eqref{omega} and the vector f\/ield~\eqref{vector field last example}.
The inf\/initesimal generator corresponding to the Lie pseudo-group action~\eqref{form-vector coordinate transformation
rules} is
\begin{gather*}
%\label{form-vector infinitesimal generator}
\ww=\xi(x,y,z)\pp{}{x}+\eta(x,y,z)\pp{}{y}+\tau(x,y,z)\pp{}{z}+[c \tau_x-a(\xi_x+\eta_y)-b \tau_y]\pp{}{a}
\\
\phantom{\ww=}
{}-[a \eta_z+b(\xi_x+\tau_z)+c \eta_x]\pp{}{b}+[a \xi_z-b \xi_y-c(\eta_y+\tau_z)]\pp{}{c}
+\frac{g}{a}[c \tau_x-b \tau_y+a \tau_z]\pp{}{g},
\end{gather*}
where $\xi(x,y,z)$, $\eta(x,y,z)$ and $\tau(x,y,z)$ are arbitrary analytic functions.
The rank of the Lie matrix,~\cite{O-2000}, of the f\/irst prolongation $\ww^{(1)}$ reveals that the orbits of the
f\/irst order prolonged action are of codimension one in $\J^1$.
Hence the equivalence pseudo-group admits a~f\/irst order dif\/ferential invariant:
\begin{gather}
\label{invariant}
I=\frac{g}{a}(a_z-b_y+c_x).
\end{gather}

Implementing the moving frame algorithm, the analysis of the recurrence relations reveals that it is possible to make
the ``universal'' normalizations
\begin{gather}
X=Y=Z=0,
\qquad
A=G=1,
\qquad
B_{X^i Y^j Z^k}=C_{X^i Y^j Z^k}=0,
\qquad
i+j+k\geq0,
\nonumber
\\
G_{X^i Y^j Z^k}=0,
\qquad
i+j+k\geq1,
\qquad
A_{X^iY^j}=0,
\qquad
i+j\geq1,
\label{v-omega universal normalizations}
\end{gather}
leading to the normalization of the Maurer--Cartan forms
\begin{gather*}
\mu=\bl(\xi),
\qquad
\nu=\bl(\eta)
\qquad
\alpha=\bl(\tau),
\qquad
\nu_{X^i Y^{j+1}}=\bl(\eta_{x^i y^{j+1}}),
\\
\nu_{X^i Y^j Z^{k+1}}=\bl(\eta_{x^i y^j z^{k+1}}),
\qquad
\mu_{X^i Y^j Z^{k+1}}=\bl(\xi_{x^i y^j z^{k+1}}),
\qquad
\alpha_{X^i Y^j Z^{k+1}}=\bl(\tau_{x^i y^j z^{k+1}}),
\end{gather*}
$i,j,k\geq0$.
With the normalizations~\eqref{v-omega universal normalizations}, the invariant~\eqref{invariant} corresponds to
\begin{gather}
\label{I=A_Z}
I=A_Z=\iota(a_z),
\end{gather}
and the invariant coframe $\bo=\{\omega^x, \omega^y, \omega^z\}=\iota\{dx, dy, dz\}$ is such that
$\Omega=\omega^x\wedge\omega^y$, and $\vv\interior\omega^z=1$.
Geometrically, the invariant~\eqref{invariant} measures the obstruction of $\Omega$ to being closed:
\begin{gather*}
d\Omega=I \omega^x\wedge\omega^y\wedge\omega^z.
\end{gather*}
After making the normalizations~\eqref{v-omega universal normalizations}, the only remaining partially normalized
invariants are
\begin{gather}
\label{v-omega remaining invariants}
A_{X^i Y^j Z^{k+1}},
\qquad
i,j,k\geq0.
\end{gather}
Up to order 2, the recurrence relations for these invariants are
\begin{gather}
dA_Z\equiv A_{XZ} \omega^x+A_{YZ} \omega^y+(A_{ZZ}-A_Z^2)\omega^z,
\nonumber
\\
dA_{XZ}\equiv A_{XXZ} \omega^x+A_{XYZ} \omega^y+(A_{XZZ}-2A_Z A_{XZ})\omega^z
\nonumber
\\
\phantom{dA_{XZ}\equiv}
{}-A_{YZ} \nu_X-A_{XZ} \mu_X-\big(A_{ZZ}-A_Z^2\big)\alpha_X,
\nonumber
\\
dA_{YZ}\equiv A_{XYZ} \omega^x+A_{YYZ} \omega^y+(A_{YZZ}-3A_Z A_{YZ})\omega^z
\nonumber
\\
\phantom{dA_{YZ}\equiv}
{}+A_{YZ} \mu_X-A_{XZ} \mu_Y-\big(A_{ZZ}-A_Z^2\big)\alpha_Y,
\nonumber
\\
dA_{ZZ}\equiv A_{XZZ} \omega^x+A_{YZZ} \omega^y+(A_{ZZZ}-A_Z A_{ZZ})\omega^z.
\label{v-omega recurrence relations}
\end{gather}
At this stage, the equivalence problem splits in two branches:
\begin{description}\itemsep=0pt
\item[Case 1:] $A_Z$ is constant.
\item[Case 2:] $A_Z$ is not constant.
\end{description}

In Case~1, when $A_Z\equiv c$ is constant, the recurrence relations~\eqref{v-omega recurrence relations}, and the
higher order ones, imply that the remaining partially normalized invariants~\eqref{v-omega remaining invariants} are
constant.
For example, from~\eqref{v-omega recurrence relations} we have that
\begin{gather*}
\begin{split}
& A_Z\equiv c,
\qquad
A_{XZ}\equiv A_{YZ}\equiv0,
\qquad
A_{ZZ}\equiv c^2,
\\
& A_{XXZ}\equiv A_{XYZ}\equiv Z_{XZZ}\equiv A_{YYZ}\equiv A_{YZZ}\equiv0,
\qquad
A_{ZZZ}\equiv c^3.
\end{split}
\end{gather*}
Hence, in this case, there are no further normalizations possible.
Up to order~2, the normalized Maurer--Cartan forms are
\begin{gather}
\mu\equiv-\omega^x,
\qquad
\nu\equiv-\omega^y,
\qquad
\alpha=-\omega^z,
\qquad
\nu_Y\equiv-\mu_X+A_Z \omega^z,
\qquad
\nu_Z\equiv\mu_Z\equiv\alpha_Z\equiv0,\nonumber \\
\nu_{XY}\equiv-\mu_{XX}-A_Z \alpha_X-A_{XZ} \omega^z,
\qquad
\nu_{YY}\equiv-\mu_{XY}-A_Z \alpha_Y+A_{YZ} \omega^z,\nonumber \\
\nu_{ZX}\equiv\nu_{ZY}\equiv\nu_{ZZ}\equiv\mu_{XZ}\equiv\mu_{YZ}\equiv\mu_{ZZ}\equiv\alpha_{XZ}\equiv\alpha_{YZ}\equiv\alpha_{ZZ}\equiv0.
\label{form-vector MC relations}
\end{gather}
The order 1 isotropy annihilator polynomials associated to~\eqref{form-vector MC relations} are
\begin{gather}
\label{reduced isotropy annihilator polynomials}
t_y T^y+t_xT^x,
\qquad
t_z T^y,
\qquad
t_z T^x,
\qquad
t_z T^z,
\end{gather}
and the indices of the reduced\footnote{Notice that in~\eqref{reduced isotropy annihilator polynomials} we have omitted
the polynomials in $T^a$, $T^b$, $T^c$ and $T^g$ associated to the dependent variables of the problem as in this
example these do not play an essential role in the involutivity test.
The corresponding Cartan characters are all zero.} symbol matrix $\mathbf{T}^1$ are $\beta^{(3)}_1=3$,
$\beta^{(2)}_1=1$, $\beta^{(1)}_1=0$.
At order 2, the rank of the reduced symbol matrix $\mathbf{T}^2$ is
$11=3 \beta^{(3)}_1+2 \beta^{(2)}_1+\beta^{(1)}_1$, and Cartan's involutivity test is satisf\/ied.
The corresponding Cartan characters are $\alpha^{(3)}_1=0$, $\alpha^{(2)}_1=2$, and $\alpha^{(1)}_1=3$.

The structure equations for the horizontal coframe $\bo=\{\omega^x$, $\omega^y$, $\omega^z\}$ are
\begin{gather*}
d\omega^x\equiv\mu_X\wedge\omega^x+\mu_Y\wedge\omega^y,
\\
d\omega^y\equiv\nu_X\wedge\omega^x-\mu_X\wedge\omega^y-c \omega^y\wedge\omega^z,
\\
d\omega^z\equiv\alpha_X\wedge\omega^x+\alpha_Y\wedge\omega^y.
\end{gather*}
These equations are equivalent to those obtained with Cartan's method,~\cite[equation~(11.29)]{O-1995}.
The correspondence is given by
\begin{gather}
\theta^1\leftrightarrow\omega^x,
\qquad
\theta^2\leftrightarrow\omega^y,
\qquad
\theta^3\leftrightarrow\omega^z,
\qquad
\alpha^1\leftrightarrow\mu_X,
\qquad
\alpha^2\leftrightarrow\mu_Y,
\qquad
\alpha^3\leftrightarrow\nu_X,
\nonumber
\\
\beta^1\leftrightarrow\alpha_X,
\qquad
\beta^2\leftrightarrow\alpha_Y,
\qquad
T\leftrightarrow A_Z.
\label{v-omega correspondence}
\end{gather}

Moving to Case~2, where $A_Z$ is not constant, two sub-branches must be considered:
\begin{description}\itemsep=0pt
\item[Case 2.1:] $A_{ZZ}\equiv A_Z^2$.
\item[Case 2.2:] $A_{ZZ}\not\equiv A_Z^2$.
\end{description}

In Case~2.1, since $dA_{Z}\not\equiv0$, the recurrence relation for $A_Z$ in~\eqref{v-omega recurrence relations}
implies, that $(A_{XZ},A_{YZ})$ $\not\equiv(0,0)$.
Assuming $A_{XZ}>0$, we set
\begin{gather*}
A_{XZ}=1,
\qquad
A_{YZ}=0,
\qquad
A_{X^i Y^j Z}=0,
\qquad
i+j\geq2,
\end{gather*}
and normalize the Maurer--Cartan forms $\mu_{X^i Y^j}$, $i+j\geq1$.
The f\/irst recurrence relation in~\eqref{v-omega recurrence relations} then reduces to
\begin{gather*}
dA_Z\equiv\omega^x.
\end{gather*}
Substituting the equality $A_{ZZ}=A_Z^2$ into the last equation of~\eqref{v-omega recurrence relations} we obtain
\begin{gather*}
2A_Z \omega^x=A_{XZZ} \omega^x+A_{YZZ} \omega^y+\big(A_{ZZZ}-A_Z^3\big)\omega^z,
\end{gather*}
which means that
\begin{gather*}
A_{XZZ}\equiv2A_Z,
\qquad
A_{YZZ}\equiv0,
\qquad
A_{ZZZ}\equiv A_Z^3.
\end{gather*}
Similarly, all higher order invariants $A_{X^i Y^j Z^{k+2}}=F_{ijk}(A_Z)$ are function of the invariant~\eqref{I=A_Z}.
Hence, no further normalizations are possible.
Up to order 2, the normalized Maurer--Cartan forms are
\begin{gather}
\mu\equiv-\omega^x,
\qquad
\nu\equiv-\omega^y,
\qquad
\alpha=-\omega^z,
\qquad
\nu_Y\equiv A_Z \omega^z,
\nonumber
\\
\nu_Z\equiv\mu_Z\equiv\alpha_Z\equiv\mu_X\equiv\mu_Y\equiv0,
\qquad
\nu_{XY}\equiv-A_Z \alpha_X-\omega^z,
\qquad
\nu_{YY}\equiv-A_Z \alpha_Y,
\nonumber
\\
\nu_{ZX}\!\equiv\!\nu_{ZY}\!\equiv\!\nu_{ZZ}\!\equiv\!\mu_{XZ}\!\equiv\!\mu_{YZ}\!\equiv\!\mu_{ZZ}\!\equiv\!\alpha_{XZ}
\!\equiv\!\alpha_{YZ}\!\equiv\!\alpha_{ZZ}\!\equiv\!\mu_{XX}\!\equiv\!\mu_{XY}\!\equiv\!\mu_{YY}\!\equiv\!0.
\label{form-vector MC relations case 2.1}
\end{gather}
The order 1 isotropy annihilator polynomials associated with~\eqref{form-vector MC relations case 2.1} are
\begin{gather*}
t_y T^y,
\qquad
t_z T^y,
\qquad
t_z T^x,
\qquad
t_z T^z,
\qquad
t_xT^x,
\qquad
t_yT^x
\end{gather*}
and the indices of the reduced symbol matrix $\mathbf{T}^1$ are $\beta^{(3)}_1=3$, $\beta^{(2)}_1=2$, $\beta^{(1)}_1=1$.
At order 2, the rank of the reduced symbol matrix $\mathbf{T}^2$ is
$14=3 \beta^{(3)}_1+2 \beta^{(2)}_1+\beta^{(1)}_1$, which satisf\/ies Cartan's involutivity test.
The corresponding Cartan characters are $\alpha^{(3)}_1=0$, $\alpha^{(2)}_1=1$, and $\alpha^{(1)}_1=2$.
Finally, the structure equations for the horizontal coframe $\bo=\{\omega^x$, $\omega^y$, $\omega^z\}$ are
\begin{gather*}
d\omega^x\equiv0,
\qquad
d\omega^y\equiv\nu_X\wedge\omega^x-A_Z \omega^y\wedge\omega^z,
\qquad
d\omega^z\equiv\alpha_X\wedge\omega^x+\alpha_Y\wedge\omega^y.
\end{gather*}
Once again, using the correspondence~\eqref{v-omega correspondence} we recover the structure equations obtained via
Cartan's method,~\cite[p.~370]{O-1995}.
Case~2.2 is treated in a~similar fashion.
The complete analysis, based on Cartan's equivalence method, can be found in~\cite{GS-1992}.
\end{Example}

\subsection*{Acknowledgments}

I am grateful to the three anonymous referees for their comments and suggestions which helped improve the paper.
I would also like to thank Abraham Smith, Niky Kamran and Robert Milson for their suggestions and stimulating
discussions, and Peter Olver for his comments on an early draft of the manuscript.
This project was supported by an AARMS Postdoctoral Fellowship and a~NSERC of Canada Postdoctoral Fellowship.

\pdfbookmark[1]{References}{ref}
\LastPageEnding

\end{document}